\documentclass{article}
\usepackage{latexsym}
\usepackage{epsfig}

\begin{document}

\title{\huge Family Switching Formula and the $-n$ Exceptional Rational Curves}
\author{Ai-Ko Liu\footnote
{email address: akliu@math.berkeley.edu, Current Address: 
Mathematics Department of U.C. Berkeley}}
\maketitle
\newtheorem{main}{Main Theorem}
\newtheorem{theo}{Theorem}
\newtheorem{lemm}{Lemma}
\newtheorem{prop}{Proposition}
\newtheorem{rem}{Remark}
\newtheorem{cor}{Corollary}
\newtheorem{mem}{Examples}
\newtheorem{defin}{Definition}
\newtheorem{axiom}{Axiom}
\newtheorem{conj}{Conjecture}

\section{preliminary}

In this paper, we prove the switching formula for Family
 Seiberg-Witten invariants. 
 Besides its role in the family Seiberg-Witten theory, the formula is
 also one of the main ingredients in the proof of the G${\ddot o}$ttsche and
 the G$\ddot o$ttsche-Yau-Zaslow conjecture[Liu1].

 Recall that in the ordinary Seiberg-Witten theory, a simple formula
 similar to the blow up formula was proved by R. Stern and
 R. Fintushel[FS] regarding the change of Seiberg-Witten invariants for the
 two different $spin^c$ structures which differ by a multiple of
 classes with self intersection number $=-n$.
 Our formula generalizes theirs to the family Seiberg-Witten theory.
It turns out that through our generalization, the Gromov-Taubes
 aspect of R. Stern and R. Fintushel's formula can be fully explored.

 In an earlier paper[Liu3], the current
 author has derived the family blowup formula for $-1$ spheres. 
The family blowup formula plays a rather crucial role in understanding the
 enumerating question upon the number of nodal curves on an algebraic surface.
 In this paper, we generalize the approach to prove the $-n$ sphere
 switching formula.

In the following, we list all the few main theorems proved in the
paper.

 The following theorem is the family switching formula of $-n$ spheres.

\medskip

\begin{main}\label{main; 1} 
 Let $\pi:{\cal X}\mapsto B$, $\pi:{\cal C}\mapsto B$, ${\cal L}$ and
 ${\cal L}_k$ 
  be the fiber bundle of four-manifolds, the relative 
 $S^2$ fiber bundle with fiberwise self-intersection number $-n$, the
 two fiberwise $spin^c$ structures  .

 Assume that the family moduli space expected dimensions are
 non-negative, then $\exists {\bf V}\mapsto B$, a complex virtual vector 
bundle over $B$ called the relative obstruction virtual bundle and the
 following family switching formula relates the family invariants of $spin^c$
 structures ${\cal L}_k$ and ${\cal L}$, 

 $$FSW_B(c, {\cal L}_k)=\sum_{i\geq 0} FSW_B(c_i({\bf V})\cup c, {\cal L}).$$

\end{main}
 
\medskip

 For the notations and the structure of ${\bf V}$, please consult section 
 \ref{section; setup} and the statement of theorem \ref{theo; 1} on page
 \pageref{theo; 1}.

\medskip

 The following theorem for algebraic family Seiberg-Witten invariants
 in section \ref{section; apfsf}
is the algebraic analogue of the switching formula in which we allow
 the ${\bf P}^1$ fibration to have singular fibers. 

\medskip

\begin{main}\label{main; alsw}
 Let $\pi:{\cal X}\mapsto B$ be an algebraic fiber bundle of algebraic surfaces 
 over a proper and smooth algebraic manifold $B$. Let $C$ be a $(1, 1)$ class
 of ${\cal X}$ which restricts to $(1, 1)$ classes of the fibers.

 Let ${\cal C}\subset {\cal X}$ be a rational curve fibration over $B$ with
 smooth generic fibers and self-intersection number $-n$.

 Under these assumptions  $AF(i).-AF(iv).$ (see page \pageref{theo; alsw} for
the details), the pure 
algebraic family invariants of $C+kPD({\cal C})$
 and the mixed family invariant of $C$ are related by the following formula,

 $${\cal AFSW}_{{\cal X}\mapsto B}(1, C+kPD({\cal C}))=\sum_{0\leq i<\infty}
 {\cal AFSW}_{{\cal X}\mapsto B}(c_i({\cal V}_{1\mapsto k}), C),$$

 where the relative obstruction virtual sheaf ${\cal V}_{1\mapsto k}$ can be
 identified with the virtual sheaf
 $${\cal R}^1\pi_{\ast}{\cal O}_{k{\cal C}}(D_C+k{\cal C})-{\cal R}^0\pi_{\ast}
 {\cal O}_{k{\cal C}}(D_C+k{\cal C}).$$
\end{main}

\medskip

  In section \ref{section; canon} and section \ref{section; example} we 
apply the ideas embedded in the proof of the family switching formula to
some concrete example on the restriction of the universal family 
${\cal X}=M_{n+1}\times_{M_n}Y(\Gamma)\mapsto Y(\Gamma)=B$. In section 
\ref{section; canon} we compare the canonical algebraic Kuranishi models of
 two classes $C-{\bf M}(E)E-\sum_{e_i\cdot (C-{\bf M}(E)E)<0}e_i$ and 
$C-{\bf M}(E)E$. In section \ref{section; example}, we apply the
 technique of localized top Chern classes and identify certain localized 
contribution of ${\cal AFSW}_{M_{n+1}\times T(M)\mapsto 
 M_n\times T(M)}(1, C-{\bf M}(E)E)$ to be
 $${\cal AFSW}_{M_{n+1}\times_{M_n}Y(\Gamma)\times T(M)\mapsto Y(\Gamma)\times 
T(M)}(c_{total}(\tau), C-{\bf M}(E)E-\sum_{e_i\cdot (C-{\bf M}(E)E)<0}e_i).$$

\medskip

\begin{main}\label{main; question1}
 Under the {\bf Simplifying Assumption}, the integer
 ${j_{Y_1}}_{\ast}{\bf Z}_{Y_1}(s_{canon})\cap 
c_1(i_Y^{\ast}{\bf H})^{dim_{\bf C}M_n+{C^2-C\cdot 
c_1({\bf K}_M)\over 2}+p_g-\sum_{i\leq n} {m_i^2+m_i\over 2}}$
 , representing the dominated localized contribution of $Y(\Gamma)$ to 
${\cal AFSW}_{M_{n+1}\times T(M)\mapsto M_n\times T(M)}(1, C-{\bf M}(E)E)$, 
can be identified with the mixed family invariant

$${\cal AFSW}_{M_{n+1}\times_{M_n}Y(\Gamma)\times T(M)\mapsto Y(\Gamma)\times 
T(M)}(c_{total}(\tau), C-{\bf M}(E)E-\sum_{i\leq p} e_{k_i})$$

 for some $\tau\in K_0(Y(\Gamma)\times T(M))$.
\end{main}

 Please consult section \ref{section; example} and
 page \pageref{assum} for the notations and the details of the simplifying
 assumption. 

\bigskip

\section{\bf A Simple Review of Fintushel-Stern's Argument}

\bigskip

 After the discovery of the Seiberg-Witten theory [W], the
 calculation of Seiberg-Witten invariants has become an important subject as
 Seiberg-Witten invariants give rise to the smooth diffeomorphism invariants of the
 four-manifolds. It has been the long term goal of several group of
 people to understand the behavior of the Donaldson or Seiberg-Witten
 invariants under several kinds of surgerical operation.  In the
 present, we do not attempt to add new ingredients into this beautiful
 theory. Instead, we would like to generalize a simple formula of
 R. Fintushel and R. Stern to the family version of Seiberg-Witten
 invariants. As it will be shown in the second section, the proof of
 the new formula does not involve any new trick or technical
 improvement of Seiberg-Witten-Floer theory. Instead, the conjecture
 raised at the end of the paper strongly indicates the simplicity of
 the picture while the corresponding Floer theory could be rather complicated.

 It was R. Fintushel and R. Stern who first noticed the importance of
 the switching formula in the context of the original Seiberg-Witten theory.

  Let $C$ be a cohomology class in $H^2(M, {\bf Z})$ 
with $C^2=-n$, then $C$ determines a complex line bundle $E_C$ over $M$. 
Suppose the class $C$
 is represented by a $-n$ two-sphere in the four-manifold $M$, then one
 considers the determinant bundles of the $spin^c$ structures of the
 following form 

 $${\cal L}_k={\cal L}_0\otimes_{\bf C} E_C^{2k}.$$

 To simplify the notation, we will adopt the following alternative
 additive notation

 $${\cal L}_k={\cal L}_0+2kC, k\in {\bf Z}.$$

\medskip 

 \noindent {\bf Question}: How do the Seiberg-Witten invariants
 of ${\cal L}_0$ and ${\cal L}_k$ relate to each other?

\medskip

 As the invariant is defined to be zero if the dimension of its moduli space
  is negative, it is more interesting to consider the case that both the
 moduli space dimension ${\cal L}_0$ and ${\cal L}_k$ are  non-negative.
  They give explicit lower bounds to $c_1({\cal L}_0)^2$ and $c_1({\cal L}_k)^2$.

 With this convention understood, one may state the following theorem of
 Fintushel-Stern [FS].

\begin{theo}(Fintushel-Stern)
 Let $C\in H^2(M, {\bf Z})$ be represented by a $-n$ two-sphere in the
 smooth four manifold $M$. Let ${\cal L}_0$ and
 ${\cal L}_k={\cal L}_0+2kC$ be the determinant line bundles of the 
two $spin^c$ structures related by the tensoring of 
 $2k$ multiple of the complex line bundle $E_C$ associated with $C$.

 Suppose that the dimensions of the Seiberg-Witten moduli spaces 
 of ${\cal L}_0$, ${\cal L}_k$, ${c_1({\cal L}_0)^2-2\chi-3\sigma\over 4}$,
 ${c_1({\cal L}_k)^2-2\chi-3\sigma\over 4}$ are non-negative, 

 then the Seiberg-Witten invariants of both $spin^c$ structures are equal
 to each other,
 
  $$SW({\cal L}_k)=SW({\cal L}_0).$$
\end{theo}

 As the proof is rather similar to the proof of blowup formula, let
 us indicate only a few key points. 
 
Let us take a
 tubular neighborhood of the embedded
 $S^2$ (poincare dual to $C$) inside $M$ which is diffeomorphic to a
 two dimensional disk bundle $D$ over the $S^2$. As the self-intersection 
number is $C^2=-n$, the normal bundle of the $S^2$ can be
 identified with a complex rank one bundle of negative first chern
 class $-n$. In particular, the boundary of the disk bundle $\partial D$, 
an $S^1$ bundle
 over the embedded $S^2$, is an
 oriented three dimensional manifold, denoted by $L_n$. 
R. Fintushel and R. Stern observe that
 $L_n, n\in {\bf N}$ are lens spaces.
A key property used in the proof of Stern-Fintushel is the fact that
 $L_n$ carry positive scalar curvature metrics. This follows from the
 fact that $S^3$ is the universal covering of all the $L_n$.

 Exactly as was
used in the proof of the blowup formula, we consider the so-called long neck
 metric with a neck isometric to $L_n\times {\bf R}$ bridging the $S^2$
 and $M-D$.
  The vanishing theorem of Seiberg-Witten invariants on manifolds
 with positive scalar curvature metrics [W] simplifies the
 gluing argument dramatically and it implies that the Seiberg-Witten 
moduli spaces of ${\cal L}_0$ and ${\cal L}_k$, ${\cal M}_{{\cal L}_0}$
 and ${\cal M}_{{\cal L}_k}$ are diffemorphic.
 A simple obstruction bundle calculation leads to the proof of of the
 simple equality.

 The reader should consult [FS], [Liu3] for more details.

\bigskip

\section{The Set Up for Family Switching Formula and a Sketch of its
 Proof}\label{section; setup}

\bigskip

 Let us discuss the set up of our main theorem and fix our notations.
 Let $\pi:{\cal X}\mapsto B$ be a fiber bundle over $B$ such that the fibers are 
diffeomorphic to a smooth four manifold $M$ with $b^+_2>0$. 
 Let ${\cal C}\mapsto B$ be a
 ${\bf S}^2$ fiber bundle over $B$ which is embedded into ${\cal X}\mapsto
 B$. Let ${\bf N}_{\cal C}
\mapsto {\cal C}$ be the real rank two normal bundle of ${\cal C}\mapsto B$
 in ${\cal X}$. We assume that the fiberwise degree of
  ${\bf N}_{\cal C}$ along 
${\cal C}\mapsto B$ is negative $=-n, n\in {\bf N}$, our goal is 
 to compare the family Seiberg Witten invariants with
 different multiplicities along ${\cal C}\mapsto B$.

 One assumes additionally that a fiberwise
 almost complex structure has been equipped upon 
some tubular neighborhood of ${\cal C}\subset {\cal X}$ such that the 
 fiber bundle ${\cal C}\mapsto B$ is a relative ${\bf CP}^1$ pseudo-holomorphic
 embedding into ${\cal X}\mapsto B$.
Then
 ${\bf N}_{\cal C}$ is equipped with a complex structure and
 is viewed as a complex line bundle over $B$.
  Let ${\bf K}$ denote the canonical line bundle of the fiberwise
 almost complex structure. Then it follows that 
 ${\bf K}|_{cal C}\cong {\bf N}_{\cal C}^{\ast}\otimes 
{\bf T}_{{\cal C}/B}^{\ast}$.

\medskip

Let ${\bf Q}\mapsto {\cal C}$ be a complex line bundle over the
 ${\bf CP}^1$ bundle ${\cal C}\mapsto B$.
 Then there exists a $\bar{\partial}_{\bf Q}$ operator on ${\bf Q}$

$$\bar{\partial}_{\bf Q}:\Omega_{{\cal C}/B}^{0, 0}\otimes{\bf Q}\mapsto 
\Omega_{{\cal C}/B}^{0, 1}\otimes {\bf Q},$$

making ${\bf Q}$ a relative holomorphic line bundle on ${\cal C}\mapsto B$.
The virtual bundle ${\bf R}^{\cdot}\pi_{\ast}({\bf Q})$ is defined to be 
the family index bundle $IND(\bar{\partial}_{\bf Q})=Ker(\bar{\partial}_{\bf Q})-
Coker(\bar{\partial}_{\bf Q})$.

\medskip

 Let ${\cal L}$ be the determinant line bundle of a fiberwise invariant 
$spin^c$ structure on ${\cal X}$.  Let $D({\bf N}_{\cal C})$ be the
 disc bundle inside the normal bundle of ${\cal C}$ in ${\cal X}$.
 Then the boundary of $D({\bf N}_{\cal C})$
 is a lens space bundle over $B$.
 Let us denote $\partial(D({\bf N}_{\cal C}))$ 
by $L(n)$. Under this convention $L(1)$  has a structure of an
 ${\bf S}^3$ bundle over $B$ and the present set up is reduced to
the family blowup setting discussed in [Liu3].

 Let $PD({\cal C})\in H^2({\cal X}, {\bf Z})$ denote the poincare dual of 
 the sub-manifold ${\cal C}\mapsto {\cal X}$. Then we may follow the
 convention in the ordinary Seiberg-Witten theory and denote 
 ${\cal L}_k={\cal L}+2kPD({\cal C})$ to be the determinant line bundle
 of a new $spin^c$ structure by tensoring ${\cal L}$ with $2k$ power of
 the complex line bundle determined by $PD({\cal C})$.  
Then we want to relate the family Seiberg Witten invariants of the
 $spin^c$ structures ${\cal L}$ and ${\cal L}_k={\cal L}+2kPD({\cal C})$.

 Let the relative 
degree of ${\cal L}$ along ${\cal C}\mapsto B$ be $m\in {\bf N}$.
As ${\cal L}$ is a characteristic element restricted fiberwise over
 $B$, $\int_{{\cal C}/B}\bigl({\cal L}+{\cal C}\bigr)=m+n$ 
must be an even integer. As in [LL1] we introduce the
 pure and mixed family Seiberg-Witten invariants and denote them
 by $FSW_B(1,{\cal L})$ and $FSW_B(c,{\cal L}), c\in H^{\cdot}(B,{\bf Z})$,
respectively.  Then we have the following main theorem on the switching
 process, 

\begin{theo}\label{theo; 1} 
 Let $\pi:{\cal X}\mapsto B$, $\pi:{\cal C}\mapsto B$, $L(n)\mapsto B$, 
 ${\bf N}_{\cal C}\mapsto B$, ${\cal L}$, ${\cal L}_k$ 
 and $m,n$ be as defined above and let the family moduli space dimensions
  
$${\int_{{\cal X}/B}c_1({\cal L})^2-2\chi-3\sigma\over 4}+dim_{\bf R}B,
{\int_{{\cal X}/B}c_1({\cal L}_k)^2-2\chi-3\sigma\over 4}+dim_{\bf R}B$$

be non-negative, then $\exists {\bf V}\mapsto B$, a complex virtual vector 
bundle over $B$ called the relative obstruction virtual bundle and the
 following family switching formula relates the family invariants of $spin^c$
 structures ${\cal L}_k$ and ${\cal L}$, 

 $$FSW_B(c, {\cal L}_k)=\sum_{i\geq 0} FSW_B(c_i({\bf V})\cup c, {\cal L}).$$

\medskip

(i). $rank_{\bf C}{\bf V}={-km+k^2n\over 2}$.

\medskip

(ii). The virtual vector bundle ${\bf V}\mapsto B$ is virtually isomorphic to (in
 $K(B)$)
$$-\oplus_{i\leq k}{\bf R}^{\cdot}\pi_{\ast}({\bf N}_{\cal C}^i\otimes
  {\cal P})\in K(B).$$

\end{theo}

The complex line bundle ${\cal P}=\sqrt{{\cal L}\otimes {\bf K}|_{\cal C}}$
 is well-defined because $m+n\equiv 0\pmod{2}$.

\medskip

\begin{rem}\label{rem; sup}
\noindent (i). In defining the index bundle 
$R^{\cdot}\pi_{\ast}({\bf N}_{\cal C}^i\otimes {\cal P})$ of 
$\bar{\partial}_{{\bf N}_{\cal C}^i\otimes {\cal P}}$, the
 bundle ${\bf N}_{\cal C}$ has been given a
 holomorphic structure.  We have assumed that the
 map ${\cal C}\mapsto {\cal X}$ to be relative
 pseudo-holomorphic over $B$. 
 In terms of Gromov theory,
 the existence locus of a $-n$ pseudo-holomorphic rational curve
 within a generic family 
 is of complex codimension $1-n$. To get a
 relative pseudo-holomorphic embedding like
 ${\cal C}$ over the whole $B$, 
one may restrict the family to be above the existence locus of the
 rational curves. An alternative way
 to define a complex structure on ${\bf N}_{\cal C}$ is to view 
${\bf N}_{\cal C}$ as the complex line bundle induced by 
 its ${\bf S}^1$ circle bundle. 

\medskip

\noindent (ii).
 In the above theorem, we also assume ${\cal C}\mapsto B$ to
be smooth, i.e. there is no singular fibers. In the algebraic proof of
 the family switching formula (see section \ref{section; apfsf}),
 we will drop this assumption and allow
 $\pi:{\cal C}\mapsto B$ to have singular values.

\medskip

\noindent (iii). 
In the above main theorem, the virtual rank
 $rank_{\bf C}{\bf R}^{\cdot}\pi_{\ast}({\bf N}_{\cal C}^i\otimes
  {\cal P})$ is equal to $(-in+{m+n-2\over 2}+1)=-in+{m+n\over 2}$
 may not be positive.  It depends on the sign and the absolute value
 of the number $m$.

Therefore the virtual rank of the virtual bundle ${\bf V}$ is 
 equal to $\sum_{i\leq k}(in-{m+n\over 2})={-km+k^2n\over 2}$, which may not be
 positive, either.

\noindent (iv).
 Whenever the virtual rank $rank_{\bf C}{\bf V}$ is negative, the previous equality
 between the different pure and the mixed family Seiberg-Witten 
invariants could be re-interpreted as

 $$FSW(c,{\cal L})=\sum_j FSW(c_j({\bf W})\cup c,{\cal L}_k),$$
 
where the virtual vector bundle
 ${\bf W}$ is virtually isomorphic to $-{\bf V}$ in the $K$ group
 $K(B)$ and is of
 positive virtual rank. 

Or by using the well known relationship [F] between the Segre
 classes and the Chern classes, we may re-write the previous relationship
 as

 $$FSW(c,{\cal L})=\sum_j FSW(s_j({\bf V})\cup c,{\cal L}_k).$$

\end{rem}

\medskip

 Let us discuss several special cases of the family switching formula before
 giving derivation of the main theorem.

  Suppose that ${\cal X}\mapsto B$ is the fiber bundle constructed
 by blowing up a fiber bundle ${\cal X}_0\mapsto B$ along a cross section
 $s:B\mapsto {\cal X}_0$. 
 Then ${\cal C}$ is nothing but the relative
 $-1$ two-sphere in ${\cal X}$ and the disk bundle $D({\bf N}_{\cal C})\mapsto 
{\cal C}$ can be identified
 with the $\overline{{\bf CP}^2}({\bf R}^0
\pi_{\ast}({\bf N}_{\cal C}^{\ast})\oplus {\bf C})\mapsto B$
 minus the cross section at infinity.
 The family switching formula reduces to the family blowup
 formula studied in [Liu3]. The class $PD({\cal C})$ has been 
 denoted by $E$ in this
 special situation.
 Let ${\cal L}={\cal L}_d+E, {\cal L}_d\perp E$ in
 $H^2({\cal X}, {\bf Z})$ and ${\cal L}+2kE={\cal L}_d+(2k+1)E$
 be the two $spin^c$ structures discussed in the switching formula,
 then the family switching formula asserts that
 $${\bf V}=-\oplus_{1\leq j\leq k}
{\bf R}^{\cdot}\pi_{\ast}({\bf N}_{\cal C}^j \otimes {\cal P}).$$

  The ${\bf P}^1$ bundle ${\cal C}\mapsto B$ can be identified 
 with the projectification 
${\bf CP}^1({\bf R}^0\pi_{\ast}({\bf N}_{\cal C}^{\ast}))\mapsto B$ and
 one may apply the
 family index theorem (essentially the Grothendieck-Riemann-Roch formula 
in the smooth category) 
to the map $\pi:{\cal C}\mapsto B$. Let ${\bf K}_0$ denote the
 restriction of the canonical line bundle around the cross section
 $s:B\mapsto {\cal X}_0$ to $s$.
 By the adjunction formula ${\bf K}|_{\cal C}={\bf K}_0\otimes {\bf N}_{\cal C}$
 and ${\cal L}={\cal L}_d\otimes {\bf N}_{\cal C}$, 
 the complex 
line bundle ${\cal P}$ is isomorphic to ${\bf N}_{\cal C}\otimes \sqrt{{\cal
 L}_d\otimes {\bf K}_0|_{\cal C}}$. Therefore,
 ${\bf R}^{\cdot}\pi_{\ast}({\bf N}_{\cal C}\otimes 
{\cal P})=R^{\cdot}\pi_{\ast}({\bf N}_{\cal C}^{j+1}
\otimes \sqrt{{\cal L}_d|_{\cal C}\otimes {\bf K}_0})$.

 Let $N={\bf R}^0\pi_{\ast}({\bf N}_{\cal}^{\ast})$
 be the complex rank two bundle of $\bar{\partial}_{{\bf N}_{\cal}^{-1}}$
 holomorphic sections along ${\cal C}\mapsto B$. 
Then by the projection 
formula of ${\bf P}^1$ bundles (see exercise 8.3, 8.4 on 
[Har] page 253 for the corresponding statment in
 the algebraic category) the virtual
 vector bundle ${\bf R}^{\cdot}\pi_{\ast}({\bf N}_{\cal C}^{j+1}
\otimes \sqrt{{\cal L}_d|_{\cal C}\otimes {\bf K}_0})$ can be identified with
${\bf S}^{j-1}(N^{\ast})\otimes \sqrt{{\cal L}_d|_{\cal C}\otimes {\bf K}_0}$
when $j>0$. The resulting family switching formula is consistent with
 the family blowup formula derived in [Liu3].

\bigskip

 Next let us derive one simple corollary of the family switching
 formula.  Let us consider three different
 $spin^c$ structures which are related
 from one to another by consecutive switching multiplicities.
 Let ${\cal L}$, $ {\cal L}_{k_1}={\cal L}+2k_1PD({\cal C}),
{\cal L}_{k_1+k_2}={\cal L}+2(k_1+k_2)PD({\cal C})$ be the three fiberwise 
invariant $spin^c$ structures on ${\cal X}$. 

Then we can apply the family switching formula from ${\cal L}_1
 $ to
${\cal L}_2$, ${\cal L}_2$ to ${\cal L}_3 $ and ${\cal L}_1$ to ${\cal
 L}_3$, respectively.

 Let ${\bf V}_{1, 2}$, ${\bf V}_{2, 3}$  and ${\bf V}_{1, 3}$ denote the
 relative virtual
 obstruction bundles constructed by the main theorem 
for the three different switching processes.
 Then we have,

\medskip

\begin{cor}\label{cor; consistent}
 The direct sum of the relative virtual obstruction bundles 
 ${\bf V}_{1, 2}\oplus {\bf V}_{2, 3}$ is
 virtually isomorphic to  ${\bf V}_{1, 3}$ in $K(B)$.
\end{cor}

\noindent Proof of the corollary:
To prove the corollary, one shows that the direct factors in
 ${\bf V}_{1, 3}$ are in one to one correspondence with the
 direct factors in the direct sums ${\bf V}_{1, 2}\oplus {\bf V}_{2, 3}$.

 This can be achieved by showing that ${\bf V}_{2, 3}$ is isomorphic
 to $-\oplus_{k_1<j\leq k_1+k_2}{\bf R}^{\cdot}\pi_{\ast}({\bf N}_{\cal C}^j
\otimes {\cal
 P}_{1, 2})$ while ${\cal P}_{1, 2}\equiv 
\sqrt{{\cal L}\otimes {\bf K}}|_{\cal C}$, 
${\cal P}_{2, 3}\equiv \sqrt{{\cal L}_{k_1}\otimes {\bf K}}|_{\cal C}$,
${\cal P}_{1, 3}\equiv \sqrt{{\cal L}\otimes {\bf K}}|_{\cal C}$.
 It is apparent that ${\cal P}_{1, 2}={\cal P}_{1, 3}$ and 
 ${\bf N}_{\cal C}^i\otimes {\cal P}_{2, 3}={\bf
 N}_{\cal C}^{i+k_1}\otimes {\cal P}_{1, 3}$. 
 Thus ${\bf R}^{\cdot}\pi_{\ast}({\bf N}_{\cal C}^i\otimes
 {\cal P}_{2, 3})={\bf R}^{\cdot}\pi_{\ast}({\bf N}_{\cal C}^{i+k_1}\otimes
 {\cal P}_{1, 3})$  and the equality on ${\bf V}_{1, 3}$ follows easily. $\Box$

\medskip

 \noindent Proof of the main theorem:

 By tubular neighborhood theorem, one may embed the disk bundle 
$D({\bf N}_{\cal C})$  into ${\cal X}$.
 Thus the fiber bundle ${\cal X}\mapsto B$ can be separated by the lens space
 bundle $L_B(n)=\partial D({\bf N}_{\cal C})\mapsto B$ 
into two connected components ${\cal X}-D({\bf N}_{\cal C})$, 
$D({\bf N}_{\cal C})$.

Consider the fiberwise 'long neck' Riemannian metrics on the fiber
 bundle ${\cal X}\mapsto B$
 which has positive scalar curvature on the subset $D({\bf N}_{\cal C})
\cong {\cal D}\subset
{\cal X}$ (viewed as a neighborhood of ${\cal C}$ in ${\cal X}$). 
Let $\mu$ be a fiberwise self-dual two form on ${\cal X}\mapsto B$
 vanishes on $D({\bf N}_{\cal C})\subset {\cal X}$.

 Fix a fiberwise invariant $spin^c$ structure ${\cal L}$ along ${\cal X}\mapsto B$
 and denote the corresponding positive $spin^c$ spinor bundle by 
 ${\cal S}^+_{\cal L}$. The family Seiberg-Witten moduli space of ${\cal L}$
 is defined to be the set of all tuples of $(A, \Psi, b)$ which satisfy

 $$P_+F_A=q(\Psi, \Psi)+i\mu(b),$$

$$D_A(\Psi)=0,$$

  modulo the equivalence relationship of $U(1)$ gauge transformations on
 $(A, \Psi)$.
 For the definition of the quadratic symbol 
$q(\cdot , \cdot):{\cal S}_{\cal L}^+ \otimes {\cal S}_{\cal L}^+\mapsto 
 (\Omega_{{\cal X}_b}^2)_+$, please consult  page 55 of [Mor].

Let $\Omega_{{\cal X}/B}^i$, $0\leq i\leq 2$, denotes
 the vector bundles of smooth differential $i$-forms of the fibers, then
 $(\Omega_{{\cal X}/B}^2)_+$ denote the vector bundle of fiberwise 
self-dual two forms.
 Let $\Gamma({\cal X}/B, {\cal S}^{\pm}_{\cal L})$ denote the
 infinite dimensional vector bundles over $B$ of the fiberwise sections of
  positive (negative) $spin^c$ spinors in ${\cal S}^+_{\cal L}$ 
(${\cal S}^-_{\cal L}$). Let $d$ and $d^{\ast}$ denote the fiberwise 
deRham operator
 and its adjoint and let $\dot{P_+}, \dot{q}, \dot{\mu}, \dot{c(\cdot)}$ 
denote the
 infinitesimal change of the operator $P_+$, $\tau$, the self-dual
 two form $\mu$ and the Clifford multiplication $c(\cdot)$ 
under the infinitesimal change of the fiberwise Riemannian metrics 
 induced along a given direction $\in {\bf T}_bB$.
  Given a Seiberg-Witten solution $(A_0, \Psi_0)$ of a fiber above the point 
$b\in B$, the infinitesimal
 deformation complex of the family Seiberg-Witten equations is given
 by the following Fredholm operator 

\[
 \begin{array}{ccc}
  0 & d^{\ast}   &    0      \\
  \dot{P_+}dA_0-\dot{q}(\Psi_0, \Psi_0)-i\dot{\mu} & P_+d & 
c(\cdot)\cdot \Psi_0\\
\dot{c(A_0)}\cdot \Psi_0 & -2\tau(\Psi_0, \cdot) &     D_{A_0}      \\

 \end{array}
\]

$${\bf T}_bB 
\oplus \Omega_{{\cal X}/B}^1\oplus \Gamma({\cal X}/B, {\cal S}^+_{\cal L})
\longrightarrow \Omega_{{\cal X}/B}^0\oplus (\Omega_{{\cal X}/B}^2)_+\oplus
 \Gamma({\cal X}/B, {\cal S}^-_{\cal L}).$$

\medskip

The kernel of the above Fredholm operator is the Zariski tangent space of the
 family moduli space ${\cal M}_C$ at the equivalence class of $(A, \Psi, b)$ 
and the cokernel of the Fredholm operator
 is the obstruction space of the solution $(A_0, \Psi_0, b)$.

\medskip

  To determine the obstruction semi-bundle of ${\cal M}_{\cal L}$, 
it is reduced to calculate the
 $H^1$ of the deformation complex over all the equivalence classes 
$[A,\Psi, b]\in {\cal M}_{\cal L}$. 

 By our choice of the fiberwise long neck metrics on ${\cal X}\mapsto B$, the
 spinor $\Psi$ of any Seiberg-Witten solutions $[A, \Psi, b]\in {\cal M}_{\cal L}$
 vanishes along the long neck and the whole
$D({\bf N}_{\cal C})\times_B\{b\}\cong {\cal D}\times_B\{b\}\subset 
{\cal X}\times_B\{b\}$. On the other hand, the connection $A$ is anti-self-dual on
 ${\cal D}\times_B\{b\}\cong D({\bf N}_{\cal C})\times_B\{b\}$,
 i.e. $P_+F_A|_{{\cal D}\times_B\{b\}}\equiv 0$ and
 $A$ is gauge equivalent to the trivial connection on the
 boundary $\partial({\cal D}\times_B\{b\})$. Such a connection is 
usually called a reducible connection.

 The sketch of the gluing argument of the solutions of the 
Seiberg-Witten equations (please consult [FS] page 226-227 
for some details) allows us to conclude,

\medskip

\begin{lemm}\label{lemm; id}
 Let ${\cal M}_{{\cal L}_p}$, $p\in {\bf Z}$ be the family Seiberg-Witten
 moduli space of the $spin^c$ structure ${\cal L}+2pPD({\cal C})$ over
 ${\cal X}\mapsto B$, with the chosen long neck fiberwise Riemannian metric.
 
 Then the restriction map $(A, \Psi, b)\mapsto (A|_{{\cal X}\times_B\{b\}-
{\cal D}\times_B\{b\}},
 \Psi|_{{\cal X}\times_B\{b\}-{\cal D}\times_B\{b\}, b})$ 
establishes an inclusion from
 ${\cal M}_{{\cal L}_p}$ to the family Seiberg-Witten 
moduli space ${\cal M}$ of the $spin^c$ structure 
 ${\cal L}_p|_{{\cal X}-{\cal D}}$ on the fiber bundle
 ${\cal X}-{\cal D}$ of four-manifolds with long cylindrical ends.
\end{lemm}
  
\medskip

\noindent A sketch of the argument: The argument of $B=pt$ case (i.e. for
 a single four-manifold $M$) for rational blowdowns has
 been sketched in [FS]. The discussion for the general case is
 parallel. On the fiber bundle of four-manifolds with long cylindrical ends 
${\cal D}\mapsto B$ there is a unique (up to gauge equivalences) 
fiberwise anti-self-dual connection $A_{red}$ 
of the line bundle ${\cal L}_p|_{\cal D}$.

 Given any solution $(A, \Psi, b)\in {\cal M}_{{\cal L}_p}$, the 
 restriction $A|_{{\cal D}\times_B \{b\}}$ to the fiber 
${\cal D}\times_B \{b\}$ must be gauge equivalent to the restriction of the
 unique reducible connection $A_{p, red}\in {\cal A}^1({\cal D}, {\cal L}_p)$ 
on the specific
 fiber ${\cal D}\times_B\{b\}$. 

 On the other hand starting from any solution $\in {\cal M}$ of the Seiberg-Witten
 equations on the four-manifold with a long
cylindrical end, ${\cal X}\times_B\{b\}-{\cal D}\times_B\{b\}$, one may
 glue it with $(A_{p, red}|_{{\cal D}\times_B\{b\}}, 0)$ to get an approximated 
 solution
 on the closed four-manifold ${\cal X}\times_B\{b\}$. 
 By the general perturbation argument to get the exact solution, the 
obstruction to get the exact solution $\in {\cal M}_{{\cal L}_p}$ is 
 determined by a finite rank obstruction bundle over ${\cal M}$.
$\Box$

 By lemma \ref{lemm; id}
 we may view all the ${\cal M}_{{\cal L}_p}$, $1\leq p\leq k$ as subspaces  
in the moduli space ${\cal M}$ over the fiber bundle ${\cal X}-{\cal D}$ 
with cylindrical ends. In the following, ${\cal M}$ will serve as the 
ambient space.

 The global ${\bf S}^1$ gauge transformations induces
 a tautological ${\bf S}^1$ bundle over ${\cal M}$ and we denote by $e$.
 Let ${\bf e}\mapsto {\cal M}$ denote the complex line bundle
 $e\times_{{\bf S}^1}{\bf C}$. The circle bundle $e$ and the line bundle
 ${\bf e}$ are universal in the sense that their pull-back to the various
 ${\cal M}_{{\cal L}_p}$ are the equal to the tautological bundles  
 defined over ${\cal M}_{{\cal L}_p}$.

 Let $D_{A_{p, red}}:{\cal S}^+_{{\cal L}_p|_{\cal D}}\mapsto 
 {\cal S}^-_{{\cal L}_p|_{\cal D}}$ denote the fiberwise 
Dirac operator induced
 by the reducible connection $A_{p, red}$ on ${\cal L}_p|_{\cal D}$.
  
It is not too hard to see by using the infinitesimal deformation complex of 
 the Seiberg-Witten equations and the gluing of Seiberg-Witten solutions 
that the obstruction bundle defining 
 ${\cal M}_{{\cal L}_p}\subset {\cal M}$ is isomorphic to 
 ${\bf e}\otimes_{\bf C}Coker(D_{A_{p, red}})$.

 In the following, we would like to compare the vector bundles
 $Coker(D_{A_{p, red}})\mapsto B$ and $Coker(D_{A_{p-1, red}})\mapsto B$.
 
 In general, there is an equality

$$rank_{\bf C}Coker(D_{A_{p, red}})-rank_{\bf C}Coker(D_{A_{p-1, red}})
 ={1\over 2}(dim_{\bf R}{\cal M}_{{\cal L}_p}
-dim_{\bf R}{\cal M}_{{\cal L}_{p-1}}).$$

\medskip

\begin{prop}\label{prop; realize}
Let $s=dim_{\bf R}{\cal M}_{{\cal L}_p}
-dim_{\bf R}{\cal M}_{{\cal L}_{p-1}}<0$, then the equivalent class
 $[Coker(D_{A_{p, red}})-Coker(D_{A_{p-1, red}})]\in K(B)$ can be
 realized as a complex rank $-s$ vector bundle ${\bf V}_p\mapsto B$.

Let $s=dim_{\bf R}{\cal M}_{{\cal L}_p}
-dim_{\bf R}{\cal M}_{{\cal L}_{p-1}}>0$, then the equivalent class 
 $[Coker(D_{A_{p-1, red}})-Coker(D_{A_{p, red}})]$ can be realized as a
 complex rank $s$ vector bundle ${\bf V}_p\mapsto B$ 
\end{prop}
 
 The proof of this proposition will be postponed to subsection
 \ref{subsection; iden} on \pageref{subsection; iden} 
when we identify ${\bf V}_p$ explicitly.

\medskip

 Based on the existences of such ${\bf V}_p$ and its dependence on the
 difference of family Seiberg-Witten expected dimensions, we derive 
 a schematic form of the family switching formula in the next subsection. 

\bigskip

\subsection{\bf The Comparison of Family Seiberg-Witten Moduli Space Expected 
Dimensions and Family Switching Formula}\label{subsection; comp}

\medskip

 Consider the $k$-tuple of integers 
$\int_{\cal C}(c_1({\cal L})+(2i-1)PD({\cal C}))$
 while $i$ running from $i=1$ to $i=k$. 
If these numbers do not change signs at all, then one may discuss directly.
 If it happens that the numbers change signs, there exists a smallest 
 critical $1\leq k_{cric}\leq k$
 such that 
$\int_{\cal C}\bigl(c_1({\cal L})+(2k_{cric}+1)PD({\cal C})\bigr)\in {\bf Z}$ and
 $int_{\cal C}c_1({\cal L})\in {\bf Z}$ are of different signs.

 In such situation, we may split the original switching process 
into two cases and prove the 
 family switching formula from ${\cal L}$ to ${\cal L}+2k_{cric}PD({\cal C})$ and
 then from ${\cal L}+2k_{cric}PD({\cal C})$ to 
 ${\cal L}+2kPD({\cal C})$. 

 We start from the following lemma on the relative family Seiberg-Witten
 moduli space dimensions,

\medskip

\begin{lemm}\label{lemm; dim}
 Suppose ${\cal L}+2(p-1)PD({\cal C})$ 
and ${\cal L}+2pPD({\cal C})$ are two fiberwise $spin^c$
 structures of a fiber bundle ${\cal X}\mapsto B$ of smooth four-manifolds. Then
 the difference of their family Seiberg-Witten moduli space expected dimensions
 is equal to $\int_{\cal C}\bigl(c_1({\cal L})+(2p-1)PD({\cal C})\bigr)$.
\end{lemm}

\noindent Proof of the lemma:  To simplify our notation, 
define ${\cal L}'={\cal L}+2(p-1)PD({\cal C})$.
 Then the two $spin^c$ structures are ${\cal L}'$ and ${\cal L}'+2PD({\cal C})$,
 respectively.

 The difference of their family dimensions is given by

 $$\{{\bigl(c_1({\cal L}')+2PD({\cal C})\bigr)^2-2\chi-3\sigma\over 4}+
dim_{\bf R}B\}
 -\{{\bigl(c_1({\cal L}')\bigr)^2-2\chi-3\sigma\over 4}+dim_{\bf R}B\}$$

$$=\int_{\cal C}{4c_1({\cal L}')+4PD({\cal C})\over 4}=
\int_{\cal C}\bigl(c_1({\cal L}')+PD({\cal C})\bigr)=
\int_{\cal C}\bigl(c_1({\cal L})+(2p-1)PD({\cal C})\bigr).$$

 The lemma is proved. $\Box$

\medskip

 If $k_{cric}$ exists such that $\int_{\cal C}\bigl(c_1({\cal L})+
(2i-1)PD({\cal C})\bigr)$ change signs,
 the lemma implies that the family Seiberg-Witten moduli space expected
 dimensions of ${\cal L}+2iPD({\cal C})$, $0\leq i\leq k_{cric}$ are monotonically
 decreasing (increasing), while the family Seiberg-Witten moduli 
 space expected 
dimensions of ${\cal L}+2iPD({\cal C})$, $k_{cric}+1\leq i\leq k$ are
 monotonically increasing (decreasing). 

  Because the family Seiberg-Witten moduli spaces of all these $spin^c$ structures 
 defined with
 the positive scalar curvature long neck metrics are all isomorphic, the
monotonicity of the family dimensions implies that there is a 
 real rank $\bigl|\int_{\cal C}\bigl(c_1({\cal L})+(2p-1)PD({\cal C})\bigr)\bigr|$
 relative obstruction bundle relating the ``adjacant'' $spin^c$ 
structures ${\cal L}+2(p-1)PD({\cal C})$ and ${\cal L}+2(p)PD({\cal C})$ for
 different $p$.

\medskip

The proposition \ref{prop; realize} implies that the difference virtual
 bundle of the obstruction bundles defining 
${\cal M}_{{\cal L}_{p-1}}, {\cal M}_{{\cal L}_{p}}$ in ${\cal M}$ 
can be realized as a rank 
${\bigl|\int_{\cal C}\bigl(c_1({\cal L})+(2p-1)PD({\cal C})\bigr)\bigr|\over 2}$ 
complex vector bundle ${\bf e}\times_{\bf C}{\bf V}_p$ over ${\cal M}$.

  Given a complex vector bundle, the Euler class of the underlying real 
vector bundle is equal to its top Chern class. 

If $\int_{\cal C}\bigl(c_1({\cal L})+(2p-1)PD({\cal C})\bigr)<0$, then
 $[{\cal M}_{{\cal L}_p}]\in H_{\ast}({\cal M}, {\bf Z})$ 
is homologous to 

$$c_{top}({\bf e}\otimes_{\bf C}{\bf V}_p)\cap 
[{\cal M}_{{\cal L}_{p-1}}]=\bigl(\sum_i c_1({\bf e})^i\cup c_{top-i}({\bf V}_p)
 \bigr)\cap [{\cal M}_{{\cal L}_{p-1}}]$$

and we expect
 a family switching formula of the following form,

 $$FSW_B(\eta, {\cal L}+2pPD({\cal C}))=
\sum_{i\leq rank_{\bf C}{\bf V}_p}
FSW_B( \eta\cup c_i({\bf V}_p), {\cal L}+2(p-1)PD({\cal C})),$$

 for $\eta\in H^{\ast}(B, {\bf Z})$.

 For $\int_{\cal C}\bigl(c_1({\cal L})+(2p-1)PD({\cal C})\bigr) \geq 0$, we
 expect a family switching formula of the following form,

$$FSW_B(\eta, {\cal L}+2(p-1)PD({\cal C}))=
\sum_{i\leq rank_{\bf C}{\bf V}_p}
FSW_B( \eta\cup c_i({\bf V}_p), {\cal L}+2pPD({\cal C})),$$

or equivalently, if we choose $\eta$ to be a multiple of 
the total Segre class of ${\bf V}_p$ by $\tilde{\eta}\in H^{\ast}(B, {\bf Z})$,
$\eta=s_{total}({\bf V}_p)\cup \tilde{\eta}$, then we have

$$FSW_B(\tilde{\eta}, {\cal L}+2pPD({\cal C}))
=\sum_{i\geq 0} 
FSW_B(\tilde{\eta}\cup s_i({\bf V}_p), {\cal L}+2(p-1)PD({\cal C})),$$

 by using the relationship

$$c_{total}({\bf V}_p)\cup s_{total}({\bf V}_p)=1\in H^0(B, {\bf Z}).$$

\medskip

We have the following proposition combining the switching formulae for the
adjacant $spin^c$ structures,

\begin{prop}\label{prop; summing}
 Define the virtual bundle ${\bf V}_{1\mapsto k}$ by

 $${\bf V}_{1\mapsto k}=
\oplus_{p, \int_{\cal C}(c_1({\cal L})+(2p-1)PD({\cal C}))<0}{\bf V}_p
 -\oplus_{p, \int_{\cal C}(c_1({\cal L})+(2p-1)PD({\cal C}))>0}{\bf V}_p.$$

\medskip

(i). Then there is a family 
switching formula relating the pure (mixed) family invariants
 of ${\cal L}+2kPD({\cal C})$ and ${\cal L}$,

 $$FSW_B(\eta, {\cal L}+2kPD({\cal C}))=\sum_{i\leq \infty} 
FSW_B(\eta \cup c_i({\bf V}_{1\mapsto k}), {\cal L}).$$

\medskip

(ii). The virtual rank of ${\bf V}_{1\mapsto k}$ is equal to half of the
 difference of the family Seiberg-Witten moduli space expected dimensions,

 $$\sum_{p\leq k}\int_{\cal C}(c_1({\cal L})+(2p-1)PD({\cal C}))=
  \int_{\cal C}\bigl(k\cdot c_1({\cal L})+k^2PD({\cal C})\bigr).$$
 
\end{prop}

\medskip

\noindent Proof of the proposition: 

(i). Suppose that the signs of $\int_{\cal C}(c_1({\cal L})+(2p-1)PD({\cal C}))$
 do not change when the index $p\in {\bf N}$ runs from $1$ to $k$. Then 
 by induction we have
$$FSW_B(\eta, {\cal L}+2kPD({\cal C}))=FSW_B(\eta\cup c_{total}({\bf V}_k), 
 {\cal L}+(2k-2)PD({\cal C}))$$
$$=FSW_B(\eta\cup c_{total}({\bf V}_k\oplus {\bf V}_{k-1}), 
{\cal L}+2(k-2)PD({\cal C}))=\cdots=FSW_B(\eta\cup c_{total}({\bf V}_{1\mapsto k}),
 {\cal L}),$$

 by using the multiplicative property of the total Chern class under
 direct sums;

 or

$$FSW_B(\eta, {\cal L}+2kPD({\cal C}))=FSW_B(\eta\cup s_{total}({\bf V}_k), 
 {\cal L}+(2k-2)PD({\cal C}))$$
$$=FSW_B(\eta\cup s_{total}({\bf V}_k\oplus {\bf V}_{k-1}), 
{\cal L}+2(k-2)PD({\cal C}))=\cdots=
FSW_B(\eta\cup s_{total}(-{\bf V}_{1\mapsto k}),
 {\cal L}),$$

 depending on whether the initial value 
$\int_{\cal C}(c_1({\cal L})+PD({\cal C}))$ is
 negative or positive.

 On the other hand, suppose that the signs of 
$\int_{\cal C}(c_1({\cal L})+(2p-1)PD({\cal C}))$ do change, a similar
 inductive calculation from
 ${\cal L}_k$ to ${\cal L}_{k_{cric}+1}$ and then
 from ${\cal L}_{k_{cric}+1}$ shows that

$$FSW_B(\eta, {\cal L}+2kPD({\cal C}))=FSW_B(\eta\cup c_{total}(\oplus_{p\geq 
 k_{cric}+1}{\bf V}_p), {\cal L}+2k_{cric}PD({\cal C}))$$
$$=FSW_B(\eta\cup c_{total}(\oplus_{p\geq 
 k_{cric}+1}{\bf V}_p)\cup s_{total}(\oplus_{p\leq k_{cric}}{\bf V}_p), 
 {\cal L})=FSW_B(\eta\cup c_{total}({\bf V}_{1\mapsto k}), {\cal L}),$$

\medskip

or $$FSW_B(\eta, {\cal L}+2kPD({\cal C}))=FSW_B(\eta\cup s_{total}(\oplus_{p\geq 
 k_{cric}+1}{\bf V}_p), {\cal L}+2k_{cric}PD({\cal C}))$$
$$=FSW_B(\eta\cup s_{total}(\oplus_{p\geq 
 k_{cric}+1}{\bf V}_p)\cup c_{total}(\oplus_{p\leq k_{cric}}{\bf V}_p), 
 {\cal L})=FSW_B(\eta\cup c_{total}({\bf V}_{1\mapsto k}), {\cal L}),$$

depending on the sign of $\int_{\cal C}\bigl(c_1({\cal L})+PD({\cal C})\bigr)$.

\medskip

 In the discussion, we have used the product formula of the total Chern classes

$$c_{total}({\bf U}_1\oplus {\bf U}_2)=
c_{total}({\bf U}_1)\cup c_{total}({\bf U}_2),$$

 and the relationship between the total Segre and the total 
Chern classes of the virtual vector 
bundles,
 
$$s_{total}({\bf U})=c_{total}(-{\bf U}).$$

\medskip

(ii). The formula on the virtual rank of ${\bf V}_{1\mapsto k}$ is 
 verified by a direct calculation and by a usage of the simple formula
 $\sum_{p\leq k} {2p-1\over 2}=k^2$.

This ends the proof of proposition \ref{prop; summing}.
 $\Box$
 
\medskip

  In the following subsection, we prove proposition \ref{prop; realize} and
 identify ${\bf V}_p$ with the index bundle of the fiberwise 
$\bar{\partial}$ operator over
 ${\cal C}\mapsto B$.

\bigskip

\subsection{\bf The Identification of The Bundle ${\bf V}_p$}
\label{subsection; iden}

\medskip
 
  Consider the fiber bundle of four-manifolds with cylindrical ends 
 ${\cal D}\cong D({\bf N}_{\cal C})\mapsto B$ and the anti-self-dual
 connection $A_{p, red}$ of ${\cal L}_p|_{\cal D}$. The obstruction 
bundle defining ${\cal M}_{{\cal L}_p}\subset {\cal M}$ is
 ${\bf e}\otimes_{\bf C}Coker(D_{A_{p, red}})$.

To relate the differential geometrical
 calculation about the family index $IND(D_{A_{p, red}})$ 
to some algebraic geometric datum on ${\cal C}$, 
one considers the
 $n$-th Hirzebruch surface fiber bundle over $B$
 constructed from ${\cal C}$ canonically.

 Recall that $F_n$, the $n$-th Hirzebruch surface, is the rational ruled
 algebraic surface which has a ${\bf P}^1$ bundle structure over ${\bf P}^1$
 with two disjoint cross sections with self-intersection numbers
 $n$ and $-n$ respectively. For all $n\in {\bf N}\cup \{0\}$ the 
 surface $F_n$ can be constructed as
 toric varieties.  To construct the $F_n$ bundle over $B$, we
 consider the projective space bundle over ${\cal C}$, 
${\bf P}({\bf N}_{\cal C}\oplus {\bf C})\mapsto {\cal C}$.  Moreover,
 the fiber bundle ${\cal X}\mapsto B$ can be viewed as the familywise fiber sum
 between ${\cal C}\subset 
{\cal X}\mapsto B$ and ${\cal C}\cong {\bf P}({\bf N}_{\cal C})\subset 
{\bf F}_n\equiv {\bf P}({\bf N}_{\cal C}\oplus {\bf C})$ 
identified along ${\cal C}$.
 The embedded ${\bf P}^1$ bundle ${\cal C}\subset {\cal X}$ is
 the relative $-n$ curve over $B$, while ${\cal C}\cong 
{\bf P}({\bf N}_{\cal C})\subset 
{\bf P}({\bf N}_{\cal C}\oplus {\bf C})$ is the relative 
(self-intersection number=) $n$ curve over $B$.

\medskip

 Let $D({\bf N}_{\cal C})$ denote the unit disk bundle of the
 normal bundle ${\bf N}_{\cal C}$, containing a ${\bf P}^1$ sub-fiber bundle
 of self-intersection number $-n$. Then $\overline{D({\bf N}_{\cal C})}$
 can be viewed as the unit disk bundle of 
 ${\bf N}_{\cal C}$, $D({\bf N}_{\cal C}^{\ast})$, containing a ${\bf P}^1$
 sub-fiber bundle with self-intersection number $n$.
 Then the $n$-th Hirzebruch surface fiber 
 bundle ${\bf P}({\bf N}_{\cal C}\oplus {\bf C})$
 can be decomposed as the union of the
 $D({\bf N}_{\cal C})$ and $\overline{D({\bf N}_{\cal C})}$
 gluing along their common boundary
 $L(n)=\partial(D({\bf N}_{\cal C}))$ and $\overline{L(n)}=
\partial(\overline{D({\bf N}_{\cal C})})$.
 
\medskip

 Putting in a cylindrical end $\cong L_B(n)\times (-r, r)$ (and let $r\mapsto 
\infty$) between
 $D({\bf N}_{\cal C})$ and $\overline{D({\bf N}_{\cal C})}$ by 
streching the fiberwise Riemannian metrics of 
${\bf P}({\bf N}_{\cal C}\oplus {\bf C})$, the fiberwise Riemannian metrics
 defined on ${\bf F}_n\mapsto B$ forms a long neck between
 these two ${\bf P}^1$ fiber bundles. 
 For simplicity let us denote the ${\bf P}^1$
 fiber 
bundles with the self-intersection number $=n$ ($=-n$) by attaching
 the subscript $\pm$, ${\cal C}_+$ and
 ${\cal C}_-$, respectively.

  By using the positive scalar curvature Riemmanian metrics on
 the lens space fiber bundle $L_B(n)\mapsto B$, 
it is well known that the solutions of the Seiberg-Witten equations
 must decay exponentially and asymptotic to the
 reducible solutions along the long neck [FS].

 \medskip

  Up to the tensor factor ${\bf e}$, the obstruction bundle of 
 ${\cal M}_{{\cal L}_p}\subset {\cal M}$ is 
 completely determined by $A_{p, red}$ and it only depends on ${\cal D}$ instead 
of the whole ${\cal X}$.
 Thus one may glue ${\cal D}$ into ${\bf F}_n\mapsto B$ instead or
 equivalently, replace the long-necked fiber bundle
${\cal X}\mapsto B$ by the long necked fiber bundle ${\bf F}_n\mapsto B$.

 One extend ${\cal L}_p|_{\cal D}$, $p\in {\bf Z}$ simutaneously into fiberwise 
$spin^c$ structures $\hat{\cal L}_p$ on ${\bf F}_n\mapsto B$ such that
 they all match on ${\bf F}_n-{\cal D}$.

\medskip

\noindent Proof of proposition \ref{prop; realize}:

 Let $$D_{p}:\Gamma({\bf F}_n/B, {\cal S}^+_{\hat{\cal L}_p})
\mapsto \Gamma({\bf F}_n/B, {\cal S}^-_{\hat{\cal L}_p})$$ 
 and $$D_{p-1}:\Gamma({\bf F}_n/B, {\cal S}^+_{\hat{\cal L}_{p-1}})
\mapsto \Gamma({\bf F}_n/B, {\cal S}^-_{\hat{\cal L}_{p-1}})$$

 be some fiberwise Dirac operators on ${\bf F}_n\mapsto B$ of the 
 fiberwise $spin^c$ structures 
$\hat{\cal L}_p$, and $\hat{\cal L}_{p-1}$.
 By the excision property of the family index of Dirac operators and the
 homotopy equivalence of the index bundles under changes of connections, 
it suffices to show that when 
 $dim_{\bf R}{\cal M}_{{\cal L}_p}-dim_{\bf R}{\cal M}_{{\cal L}_{p-1}}<0$ (>0),
the class of the virtual bundle $[IND(D_{p-1})-IND(D_{p})]$ 
($[IND(D_{p})-IND(D_{p-1})]$) in
 $K(B)$ be realized as a complex vector bundle of rank 
${1\over 2}|dim_{\bf R}{\cal M}_{{\cal L}_p}-
dim_{\bf R}{\cal M}_{{\cal L}_{p-1}}|$.

  Apparently if we change the fiberwise 
Riemannian metrics of the fiber bundle ${\bf F}_n\mapsto B$ and the
 connections on the spinor bundles, it does not affect the above class of
 virtual family index bundles.   Thus, we are free to replace the original 
 fiberwise long necked Riemannian metrics on ${\bf F}_n\mapsto B$ to the
 fiberwise Kahler metrics on ${\bf F}_n\mapsto B$ as $F_n$ is an 
algebraic surface. Thus, the family Dirac operator $D_{p}$ or
 $D_{p-1}$ can be identified with the
 $\bar\partial+\bar\partial^{\ast}$ operator.

\medskip

 Let us review briefly about the cohomology ring structure of a Hirzebruch
 surface $F_n$.
 The Hirzebruch surface $F_n\mapsto {\bf P}^1$
 is a rational ruled surface with two cross sections
 $C_+$ and $C_-$. Let $F$ denote the divisor class generated by the fibers 
 ${\bf P}^1$. The 
 middle cohomology
 of $F_n$, $H^2(F_n, {\bf Z})$
 is generated by two classes, $[C_-]$ with $[C_-]^2=-n$,
 and $[F], with [F]^2=0$. The class $[C_+]$ is related to $[C_-]$ by
$ [C_+]=[C_-]+n[F]$.  The cone of the effective classes in $H^2(F_n, {\bf Z})$ 
is generated by the
 primitive generators of the extremal rays, $[F]$ and $[C_-]$.

 The holomorphic sections vanishing at $F, C_+, C_-$ 
define holomorphic line bundles
 on $F_n$, denoted by ${\bf E}_F, {\bf E}_{C_+}$ and ${\bf E}_{C_-}$. 
Then $c_1({\bf E}_F)=[F]$, $c_1({\bf E}_{C_+})=[C_+]$ and 
$c_1({\bf E}_{C_-})=[C_-]$, respectively.

\medskip

\begin{lemm}\label{lemm; canon}
 The first Chern class of the canonical line bundle ${\bf K}_{F_n}$ is
 equal to $-2[F]-[C_+]-[C_-]=-(n+2)[F]-2[C_-]$.
\end{lemm}

\medskip

\noindent Proof of the lemma: This follows from the canonical divisor formula
 of any toric surface (see page 85 section 4.3. of [F2]),

$${\bf K}_{F_n}={\bf E}^{-2}_{F}\otimes {\bf E}^{-1}_{C_+}\otimes 
{\bf E}^{-1}_{C_-}.$$
 $\Box$

 Given a $spin^c$ determinant line bundle $\hat{\cal L}$ 
on $F_n$, there is 
a unique connection (up to gauge equivalence) which gives $\hat{\cal L}$ 
a structure of holomorphic line bundle. Moreover
the spinor bundles ${\cal S}_{\hat{\cal L}}^+\cong \bigl(\bigwedge^{0, 0}F_n\oplus 
 \bigwedge^{0,2}F_n\bigr)\otimes \sqrt{\hat{\cal L}\otimes {\bf K}_{F_n}}$,
 ${\cal S}_{\hat{L}}^-\cong \bigwedge^{0,1}F_n\otimes 
\sqrt{\hat{\cal L}\otimes {\bf K}_{F_n}}$ and the $spin^c$ 
Dirac operator can be identified with 
the  twisted $\bar\partial+\bar\partial^{\ast}$ operator
 by the complex line bundle $\sqrt{\hat{\cal L}\otimes {\bf K}_{F_n}}$ (see 
[Law] page 395-400).

 In our case, suppose that $deg({\cal L}|_{\cal C})=a$, $a\equiv n(mod2)$, then 
 the extension $\hat{\cal L}$ over ${\bf F}_n\mapsto B$ can be identified as 
 ${\bf E}_F^a \otimes {\bf E}_{C_+}^b\otimes {\cal L}_0$ for some
 yet to be chosen even $b\in 2{\bf Z}$ and  a complex line bundle pull-back 
from the base $B$, ${\cal L}_0$. We have made use of
 the fact $[C_+]\cup [C_-]=0$.

 Accordingly the complex line bundle $\sqrt{\hat{\cal L}\otimes {\bf K}_{F_n}}$
 can be re-expressed as ${\bf E}_{F}^{{a+n\over 2}-1}\otimes
 {\bf E}_{C_+}^{{b\over 2}-1}\otimes \sqrt{{\cal L}_0}$.

By the projection formula, 
the index bundle of the ${\bf E}_{F}^{{a+n\over 2}-1}\otimes
 {\bf E}_{C_+}^{{b\over 2}-1}\otimes \sqrt{{\cal L}_0}$ 
twisted $\bar\partial+\bar\partial^{\ast}$ operator
 can be thought as the index bundle of the ${\bf E}_{F}^{{a+n\over 2}+1}\otimes
 {\bf E}_{C_+}^{{b\over 2}-1}$ twisted $\bar\partial+\bar\partial^{\ast}$ operator
 tensoring with $\sqrt{{\cal L}_0}$.

 The following lemma implies that for a suitable choice of $b\in {\bf N}$,
 the kernel of the $\bar\partial+\bar\partial^{\ast}$ operator vanishes.

\medskip

\begin{lemm}\label{lemm; van}
 Let $F_n$ be the Hirzebruch surface with a toric Kahler metric. Then for 
 all $a\in {\bf Z}, a+n\in 2{\bf Z}$, there exists at least one $b\in 2{\bf Z}$
 such that the ${\bf E}_{F}^{{a+n\over 2}-1}\otimes
 {\bf E}_{C_+}^{{b\over 2}-1}$ twisted $\bar\partial+\bar\partial^{\ast}$ operator

 $$\bar\partial+\bar\partial^{\ast}:(\Omega^{0,0}\oplus \Omega^{0,2})
\otimes {\bf E}_{F}^{{a+n\over 2}-1}\otimes
 {\bf E}_{C_+}^{{b\over 2}-1}
\mapsto
 \Omega^{0, 1}\otimes {\bf E}_{F}^{{a+n\over 2}-1}\otimes
 {\bf E}_{C_+}^{{b\over 2}-1}$$
has a trivial kernel.
\end{lemm}

\medskip

\noindent Proof of the lemma: Given the toric Kahler metric on $F_n$ and the
 holomorphic structure defined by $bar\partial$ with
 $\bar\partial^2=0$, the kernel of the twisted 
$\bar\partial+\bar\partial^{\ast}$ operator can be identified with
 $$H^0_{\bar\partial}(F_n, {\bf E}_{F}^{{a-n\over 2}-1}\otimes
 {\bf E}_{C_+}^{{b\over 2}-1})\oplus  H^2_{\bar\partial}(F_n, 
{\bf E}_{F}^{{a-n\over 2}-1}\otimes
 {\bf E}_{C_+}^{{b\over 2}-1}).$$

 It suffices to find $b\in 2{\bf Z}$ such that both $\bar\partial$ cohomologies
 vanish. It is well known that the zero-th and the second 
$\bar\partial$ cohomologies of
 the holomorphic line bundle ${\bf E}_{F}^{{a-n\over 2}-1}\otimes
 {\bf E}_{C_+}^{{b\over 2}-1}$ are isomorphic to
 the sheaf cohomology 
$H^0(F_n, {\cal O}_{F_n}(({a-n+bn\over 2}-1)F+({b\over 2}-1)C_-))$
 and 
$H^2(F_n, {\cal O}_{F_n}(({a-n+bn\over 2}-1)F+({b\over 2}-1)C_-))$, respectively.

 The second sheaf cohomology is isomorphic to 
$H^0(F_n, {\cal O}_{F_n}((-{a-n+bn\over 2}-1)F+(-{b\over 2}-1)C_-)$, by Serre 
duality on algebraic surfaces and lemma \ref{lemm; canon}.

 We make the following choice on $b$. If ${a-n\over 2}-1$ is negative, we 
 simply set $b=0$.   If ${a-n\over 2}-1\geq 0$, we choose a non-positive
  even $b$ such that $0\leq {a+n+bn\over 2}-1<n$. 

 Because $F$ and $C_-$ generate the cone of effective divisors on 
$F_n$, the dimension 
$h^0(F_n, {\cal O}_{F_n}(AF+BC_-))$ is
 nonzero only when both $A$ and $B$ are non-negative. 
 It is easy to see that our choices of $b$ make both
 ${\bf E}_F^{{a-n+bn\over 2}-1}\otimes {\bf E}_{C_-}^{{b\over 2}-1}$
 and ${\bf E}_F^{-{a-n+bn\over 2}-1}\otimes {\bf E}_{C_-}^{-{b\over 2}-1}$
 non-effective.  Thus, the lemma is proved. $\Box$

  By using lemma \ref{lemm; van} to choose the appropriate $b\in 2{\bf Z}$, the
sheaf of sections of the cokernel bundle $Coker(\bar\partial+\bar\partial^{\ast})$
 can be identified with the 
 first right derived image sheaf of the sheaf of smooth sections
 ${\cal O}_{{\bf F}_n}(({a-n+bn\over 2}+1)F+({b\over 2}-1)C_-)$ holomorphic 
along the fibers of ${\bf F}_n\mapsto B$.

\medskip

We apply lemma \ref{lemm; van} to ${\cal L}_p$ and identify
 ${\cal C}$ with the relative curve ${\cal C}_-\mapsto B$ 
in ${\bf F}_n\mapsto B$.
 Consider the following sheaf short exact sequence,

 $$\hskip -1.3in 
0\mapsto {\cal O}_{{\bf F}_n}(({a-n+bn\over 2}-1)F+({b\over 2}-2)C_-)
\mapsto {\cal O}_{{\bf F}_n}(({a-n+bn\over 2}-1)F+({b\over 2}-1)C_-)\mapsto
 {\cal O}_{{\cal C}_-}(({a-n+bn\over 2}-1)F+({b\over 2}-1)C_-)\mapsto 0.$$

 Taking the right derived image long exact sequence on the short exact sequence
  and we find that if $$\int_{\cal C}(c_1({\cal L})+2pPD({\cal C})+
c_1({\bf K}|_{\cal C}))=
\int_{\cal C}(c_1({\cal L})+(2p-1)PD({\cal C}))-2=a+n-2<0$$ then 
one has the following 
sheaf short exact sequence,
 
$$\hskip -.3in
0\mapsto {\cal R}^1\pi_{\ast}
\bigl({\cal O}_{{\bf F}_n}(({a-n+bn\over 2}-1)F+({b\over 2}-2)C_-)\bigr)\mapsto 
{\cal R}^1\pi_{\ast}
\bigl( {\cal O}_{{\bf F}_n}(({a-n+bn\over 2}-1)F+({b\over 2}-1)C_-)  
 \bigr)$$
$$\mapsto {\cal R}^1\pi_{\ast}
\bigl( {\cal O}_{{\cal C}_-}(({a-n+bn\over 2}-1)F+({b\over 2}-1)C_-) \bigr)
\mapsto 0.$$

  We have used the vanishing of 
${\cal R}^0\pi_{\ast}
\bigl({\cal O}_{{\cal C}_-}(({a-n+bn\over 2}-1)F+({b\over 2}-1)C_-)\bigr)$
 
 for a negative degree invertible sheaf  
${\cal O}_{{\cal C}_-}(({a-n+bn\over 2}+1)F+({b\over 2}-1)C_-)$ on
 ${\cal C}_-$.

This implies that the difference 

$$[{\cal R}^1\pi_{\ast}
\bigl( {\cal O}_{{\bf F}_n}(({a-n+bn\over 2}-1)F+({b\over 2}-1)C_-)  
 \bigr)-{\cal R}^1\pi_{\ast}
\bigl({\cal O}_{{\bf F}_n}(({a-n+bn\over 2}-1)F+({b\over 2}-2)C_-)\bigr)]$$

 is equivalent to the locally free ${\cal R}^1\pi_{\ast}
\bigl({\cal O}_{{\cal C}_-}(({a-n+bn\over 2}-1)F+({b\over 2}-1)C_-)\bigr)$.

 We denote the vector bundle associated with the locally free sheaf as
 ${\bf V}_p$.

 Suppose that 
 $$\int_{\cal C}(c_1({\cal L})+2pPD({\cal C})+
c_1({\bf K}|_{\cal C}))=\int_{\cal C}(c_1({\cal L})+(2p-1)PD({\cal C}))-2=
a+n-2\geq 0,$$ 

 there is another sheaf short exact sequence,

$$\hskip -.3in 
0\mapsto {\cal R}^0\pi_{\ast}
\bigl( {\cal O}_{{\cal C}_-}(({a-n+bn\over 2}-1)F+({b\over 2}-1)C_-) 
\bigr)\mapsto 
{\cal R}^1\pi_{\ast}
\bigl( {\cal O}_{{\bf F}_n}(({a-n+bn\over 2}-1)F+({b\over 2}-2)C_-)\bigr)$$
$$\mapsto 
{\cal R}^1\pi_{\ast}\bigl( {\cal O}_{{\bf F}_n}(({a-n+bn\over 2}-1)F+
({b\over 2}-1)C_-)\bigr)
\mapsto 0.$$

 It implies that 

$$[{\cal R}^1\pi_{\ast}
\bigl( {\cal O}_{{\bf F}_n}(({a-n+bn\over 2}-1)F+({b\over 2}-2)C_-)\bigr)-
{\cal R}^1\pi_{\ast}\bigl( {\cal O}_{{\bf F}_n}(({a-n+bn\over 2}-1)F+
({b\over 2}-1)C_-)\bigr)]$$

is equivalent to the locally free 
${\cal R}^0\pi_{\ast}
\bigl( {\cal O}_{{\cal C}_-}(({a-n+bn\over 2}-1)F+({b\over 2}-1)C_-) 
\bigr)$.

 In this case, we denote ${\bf V}_p$ as the vector bundle associated with
 the locally free sheaf.

 This ends the proof of proposition \ref{prop; realize}. $\Box$

\medskip

\begin{prop}\label{prop; bundle}
 The bundle ${\bf V}_p$
 is isomorphic to 
${\bf R}^1\pi_{\ast}\bigl(\sqrt{({\cal L}\otimes {\bf K})|_{\cal C}
 }\otimes {\bf N}^p_{\cal C}\bigr)$ if
 $\int_{\cal C}(c_1({\cal L})+(2p-1)PD({\cal C}))<0$. If 
 $\int_{\cal C}(c_1({\cal L})+(2p-1)PD({\cal C}))>0$, then 
 ${\bf V}_p$ is isomorphic to 
${\bf R}^0\pi_{\ast}\bigl(\sqrt{({\cal L}\otimes {\bf K})|_{\cal C}
 }\otimes {\bf N}^p_{\cal C}\bigr)$.
\end{prop}

\noindent Proof of proposition \ref{prop; bundle}:
 Based on the identification at the end of the proof of proposition
 \ref{prop; realize}, 
  our task is to identify 
${\cal O}_{{\cal C}_-}(({a-n+bn\over 2}-1)F+
({b\over 2}-1)C_-)$ with the sheaf of smooth sections of 
$\sqrt{({\cal L}\otimes {\bf K})|_{\cal C}
 }\otimes {\bf N}^p_{\cal C}$ holomorphic along the fibers of
 ${\cal C}_-\mapsto B$. 

  We notice that 
 the line bundle ${\cal L}_p|_{\cal C}$ over the curve 
 ${\cal C}\mapsto B$ can be identified with
 $({\bf E}^{a+bn}_{F}\otimes {\bf E}^b_{C_-})|_{{\cal C}_-}$
 on the relative divisor ${\cal C}_-\mapsto B$ inside the relative Hirzebruch
 surface ${\bf F}_n\mapsto B$ for an arbitrary $b\in 2{\bf Z}$.

 On the other hand, the restriction of the canonical line bundle,
 ${\bf K}|_{\cal C}$, can be identified with 
 $${\bf N}^{\ast}_{\cal C}
\otimes {\bf T}^{\ast}_{\cal C}={\bf E}^{-1}_{C_-}\otimes 
{\bf E}^{-2}_{F}|_{{\cal C}_-}={\bf E}_F^{n-2}|_{{\cal C}_-}.$$

 Thus $\sqrt{{\cal L}_p|_{{\cal C}}\otimes {\bf K}_{\cal C}}$ 
can be identified with $\sqrt{{\cal L}_0|_{{\cal C}_-}}\otimes 
({\bf E}_F^{{a+n+bn\over 2}-1}\otimes {\bf E}^{b\over 2}_{C_-})|_{{\cal C}_-}$.

  By using ${\bf E}_F^n|_{{\cal C}_-}\cong{\bf E}_{C_-}^{-1}|_{{\cal C}_-}$, 
the above line bundle is
 isomorphic to $\sqrt{{\cal L}_0|_{{\cal C}_-}}\otimes 
({\bf E}_F^{{a-n+bn\over 2}-1}\otimes {\bf E}^{{b\over 2}-1}_{C_-})|_{{\cal C}_-}$.

 Up to the factor $\sqrt{{\cal L}_0|_{{\cal C}_-}}$ pulled back from the base $B$, 
 the sheaf of smooth sections of the line bundle 
$({\bf E}_F^{{a-n+bn\over 2}-1}\otimes {\bf E}^{{b\over 2}-1}_{C_-})|_{{\cal C}_-}$
 holomorphic along ${\cal C}_-\mapsto B$
 is equal to ${\cal O}_{{\cal C}_-}(({a-n+bn\over 2}-1)F+
 ({b\over 2}-1)C_-)$.

 On the other hand, we have $({\cal L}_p\otimes {\cal L}^{-1})|_{\cal C}
\equiv {\bf N}_{\cal C}^{2p}$.
 So $\sqrt{\bigl({\cal L}_p\otimes {\bf K}\bigr)|_{\cal C}}$
 is isomorphic to $\sqrt{\bigl({\cal L}\otimes {\bf K}\bigr)|_{\cal C}}\otimes
 {\bf N}^p_{\cal C}$.

  Thus, we may identify 
${\cal O}_{{\cal C}_-}(({a-n+bn\over 2}-1)F+
 ({b\over 2}-1)C_-)$ as the sheaf of sections of
$\sqrt{\bigl({\cal L}\otimes 
{\bf K}\bigr)|_{\cal C}}\otimes{\bf N}^p_{\cal C}$. This ends the proof of
the proposition. $\Box$

\begin{cor}\label{cor; explicit}
 Let ${\cal L}|_{\cal C}$ be of degree $a$ over the ${\bf P}^1$ bundle ${\cal C}$.
 Suppose that the ${\bf P}^1$ fiber bundle ${\cal C}\mapsto B$ is the
 projectification of a complex rank two vector bundle ${\bf U}\mapsto B$, 
i.e. ${\cal C}\cong{\bf P}_B({\bf U})$, 
then the relative obstruction virtual bundle
 ${\bf V}_{1\mapsto k}$ of 
the family switching formula can be identified with
$$\oplus_{\int_{\cal C}(c_1({\cal L})+(2p-1)PD({\cal C}))>0}
{\bf S}^{{a+n\over 2}-pn-1}({\bf U})\otimes \sqrt{{\cal L}_0|_{\cal C}} 
\ominus_{\int_{\cal C}(c_1({\cal L})+(2p-1)PD({\cal C}))<0}
{\bf S}^{-{a+n\over 2}+pn-1}({\bf U}^{\ast})\otimes \sqrt{{\cal L}_0|_{\cal C}} 
,$$

 where ${\cal L}_0|_{\cal C}$ is a complex line bundle pulled-back from the
 base $B$.
\end{cor}

\medskip

\noindent Proof: The corollary is a consequence of proposition \ref{prop; summing},
proposition \ref{prop; bundle} and the projection formula of the
 ${\bf P}_B({\bf U})$.   When 
 $\int_{\cal C}(c_1({\cal L})+(2p-1)PD({\cal C}))>0$, it suffices to prove
 that $${\bf R}^0\pi_{\ast}\bigl(\sqrt{{\cal L}\otimes {\bf K}|_{\cal C}}\otimes
 {\bf N}^p_{\cal C}\bigr)\cong {\bf S}^{{a+n\over 2}-pn-1}({\bf U})\otimes 
 \sqrt{{\cal L}_0},$$

 and when $\int_{\cal C}(c_1({\cal L})+(2p-1)PD({\cal C}))<0$, then 
$${\bf R}^1\pi_{\ast}\bigl(\sqrt{{\cal L}\otimes {\bf K}|_{\cal C}}\otimes
 {\bf N}^p_{\cal C}\bigr)\cong {\bf S}^{-{a+n\over 2}+pn-1}({\bf U}^{\ast})\otimes 
 \sqrt{{\cal L}_0}.$$

 These identities 
follow from the projection formula 
(consult [Har], page 253, exercise 8.4.(a), (c) for the formula in the algebraic
 category)
of the ${\bf P}^1$ bundle
 ${\bf P}_B({\bf U})$.
$\Box$

\medskip

 By combining the discussions in subsection \ref{subsection; comp}, 
\ref{subsection; iden}, we have derived the family switching formula 
 of family Seiberg-Witten invariants in the 
 topological category.

\medskip

\subsection{\bf A Brief Remark About the Family Switching Formula and
 Gromov-Taubes theory} \label{subsection; GT}

\bigskip

 At the end of the section, we point out the hidden relationship between
 the obstruction virtual bundle ${\bf V}_{1\mapsto k}$ and the Gromov-Taubes
 theory. We restrict to the special case $deg_{{\cal C}/B}{\cal P}=deg_{{\cal C}/B}
 \sqrt{{\cal L}\otimes {\bf K}|_{\cal C}}<0$.

\medskip

\begin{prop}\label{prop; Q}
 Let $\pi:{\cal C}\mapsto B$ denote the projection map (or its
 restriction to certain subsets) of the ${\bf P}^1$ fiber bundle which admits
 a cross section $s_{\cal C}:B\mapsto {\cal C}$. Let $k$ be a positive
 integer and let $q=deg_{{\cal C}/B}{\cal P}<0$, then
 there exists a smooth complex line bundle ${\bf Q}$ 
on ${\cal C}$ pulled back from $B$ such that
  the obstruction virtual bundle ${\bf V}_{1\mapsto k}$
 can be alternatively represented in the $K$ group $K(B)$ as

 $$\sum_{k\geq i\geq 1}\bigl({\bf R}^1\pi_{\ast}({\bf N}_{\cal C}^i)
+{\bf S}^{-q-1}({\bf C}\oplus {\bf N}_{s_{\cal C}}^{-1})
\otimes {\bf N}_{\cal C}^i|_{s_{\cal C}(B)})\otimes 
 \bigr)\otimes {\bf Q}.$$
\end{prop}

\medskip

 The symbol ${\bf N}_{s_{\cal C}}$ denotes the normal bundle of
 the cross section $s_{\cal C}\subset {\cal C}$.

\medskip

\noindent Proof:   Define $\Delta=-q\cdot s_{\cal C}$ to be the 
non-reduced relative divisor in the relative ${\bf P}^1$ bundle. Then
 there is a short exact sequence relativizing the 
 $0\mapsto {\cal I}_{\Delta\times_B\{b\}}\mapsto {\cal O}_{{\cal C}\times_B\{b\}}
\mapsto {\cal O}_{\Delta\times_B\{b\}}\mapsto 0$.

$$0\mapsto {\cal O}_{\cal C}(-\Delta)\mapsto {\cal O}_{\cal C}\mapsto 
 {\cal O}_{\Delta}\mapsto 0.$$

 Because $deg_{{\cal C}/B}{\cal O}_{\cal C}(-\Delta)=q$ as well, the
 sheaf of fiberwise $\bar\partial$-holomorphic sections of 
 ${\cal P}=\sqrt{{\cal L}\otimes_{\bf C}{\bf K}|_{\cal C}}$ is equivalent to
 ${\cal O}_{\cal C}(-\Delta)\otimes {\bf Q}$ for some ${\cal C}^{\infty}$
 complex line bundle pulled back from $B$.

 By tensoring the short exact sequence with ${\bf N}_{\cal C}^i\otimes {\bf Q}$,
 $1\leq i\leq k$, and then push forward along $\pi:{\cal C}\mapsto B$,
  the equality of the $K$ group is a consequence of the equality

 $${\cal R}^0\pi_{\ast}\bigl({\cal O}_{-q\cdot s_{\cal C}}\bigr)=
\oplus_{0\leq i\leq -q-1}
 {\cal N}_{s_{\cal C}}^{-i}={\bf S}^{-q-1}({\bf O}_B\oplus 
{\cal N}_{s_{\cal C}}^{-1})$$

 for the normal sheaf ${\cal N}_{s_{\cal C}}$. $\Box$

\bigskip

 Up to the factor ${\bf Q}$ which depends on ${\cal L}$, the 
 direct sum 

 $$\sum_{k\geq i\geq 1}{\bf R}^1\pi_{\ast}({\bf N}_{\cal C}^i)$$

 only depends on $k$ and the local data on ${\cal C}\mapsto B$. At the end
 of the subsection, we make a few remark about its implicit link with Gromov-Taubes
 theory of pseudo-holomorphic curves.

\medskip

\noindent (i). The first term ${\bf R}^1\pi_{\ast}({\bf N}_{\cal C})$
 is nothing but the obstruction vector bundle of 
 ${\cal C}\mapsto B$ embedded as a pseudo-holomorphic relative curve
 in a tubular neighborhood inside ${\cal X}\mapsto B$.

 In fact, for all $b\in B$, the ${\bf P}^1$,  ${\cal C}\times_B \{b\}$
 can be viewed as an exceptional curve with self-intersection number $=-n$.
 Then its obstruction vector space 
$H^1({\cal C}\times_B\{b\}, {\bf N}_{\cal C}|_{{\cal C}\times_B\{b\}})$ 
 in Gromov theory is nothing but the fiber of
 ${\bf R}^1\pi_{\ast}({\bf N}_{\cal C})$ at $b\in B$.

\medskip

 \noindent (ii).
Assuming ${\cal X}\mapsto B$ to be symplectic and $FSW_B(1, {\cal L}_k)\not=0$, 
Taubes' hard 
analysis on $SW->Gr$ [T1] implies the existence of pseudo-holomorphic curves
 dual to $c_1(\sqrt{{\cal L}_k\otimes {\bf K}_{{\cal X}/B}})$.
 
 When $k>1$, the multiple terms 
${\bf R}^1\pi_{\ast}({\bf N}_{\cal C}^i)$, $1\leq i\leq k$ reflect that
 in enumerating the family Seiberg-Witten 
 invariant of ${\cal L}_k={\cal L}+2kPD({\cal C})$, a $k$-multiple covering of
 ${\cal C}$ contributes to the family invariant.

 On the other hand, in Taubes' gluing argument ``Gr->SW'' 
from two dimensional vortices
 to solutions of Seiberg-Witten equations (see item 5. on 
page 241 of [T2]), the index bundle of $\Delta_b$ 
 ($N={\bf N}_{\cal C}|_{{\cal C}\times_B\{b\}}, C={\cal C}\times_B\{b\}$),

$$\Delta_b:\oplus_{1\leq i\leq k}N^i\mapsto (\oplus_{1\leq i\leq k}N^i\otimes
 {\bf T}^{0, 1}C$$

with $\Delta_b(\cdot)=\bar\partial(\cdot)+\dots$ equal to 
$\bar\partial$ up to a zeroth order term, can be identified 
 (up to a homotopy of first order elliptic operators) with our 
 $\oplus_{1\leq i\leq k}{\bf R}^1\pi_{\ast}({\bf N}_{\cal C}^i)$.

This establishes an implicit link between two different gluing problems.

\medskip

\noindent (iii). Even though in Gromov-Taubes theory Taubes did develop a concept
 of multiple coverings of exceptional curves. But it differs significantly 
from the usual Ruan-Tian theory [RT]. Suppose that $C\subset M$ is a
 $-n$ exceptional rational curve in a symplectic four-manifold and let
 ${\bf N}_C$ denote the complex normal bundle $C\subset M$. Then
 the obstruction vector space at $[g_k]$ for
 a $k$-multiple covering pseudo-holomorphic map $g_k:{\bf P}^1\mapsto C$
 is $H^1({\bf P}^1, g_k^{\ast}{\bf N}_C)\cong H^1(C, {\bf N}_C^k)$.
 This corresponds to the $i=k$ term of the bundle 
$oplus_{1\leq i\leq k}{\bf R}^1\pi_{\ast}({\bf N}_{\cal C}^i)$.

  In Taubes' setting, a $k$-fold multiple covering of an exceptional curvs is
 better understood as a zero locus defined by $k$-fold multiple of the zero locus 
defining $C\subset M$, through the limiting process of large deformations
 of symplectic forms [T1], [T2], [T3]
 and the identification of zero loci of Dirac spinors
 with pseudo-holomorphic curves. In the algebraic category, this can be
 thought as putting a non-reduced scheme structure on the corresponding 
algebraic curve.
On the other hand, in the usual Gromov-Witten
 theory a $k$-fold covering of an exceptional curve is viewed as a map
 which factors through the $k$-fold covering of ${\bf P}^1\mapsto {\bf P}^1$.

 The reader should be aware of the difference between these two theories and
 the corresponding differences on their dimension formulae.

\medskip

\section{\bf The Algebraic Proof of the Family Switching Formula}
\label{section; apfsf}

\bigskip

 In section 
\ref{section; setup}, we derived the family switching formula in the topological
category.  In this section, we derive the family switching formula of the
 algebraic family invariants ${\cal AFSW}$. In the topological category,
 we have to assume that ${\cal C}\mapsto B$ is a smooth ${\bf P}^1$ bundle
 over $B$. In general, the family invariants in the smooth category are
 difficult to determine directly.  The usage of fiberwise long neck metrics
 allow us to study how do the invariants change under switching of fiberwise
 $spin_c$ structures.  On the other hand, the algebraic family Seiberg-Witten
 invariants can be read off from the topological datum of the 
algebraic Kuranishi models. This fact has two consequences, (1). One may derive
the family switching formula based on the algebraic family Kuranishi models.
 (2). One may weaken the smoothness 
assumption on the relative curve ${\cal C}\mapsto B$.

Let us begin by introducing our setup, and point out the difference from the
case of smooth topology.

 Let $\pi:{\cal X}\mapsto B$ be an algebraic fiber bundle over a smooth algebraic 
 base $B$ such that the fibers are smooth algebraic surfaces.

 Let $C$ be an $(1, 1)$ class of the total space ${\cal X}$ which restricts to
 the fiberwise $(1, 1)$ class on the fibers. Let ${\cal C}\subset {\cal X}$ 
be a relative curve over $B$ satisfying the following conditions:
\label{condition}

\medskip

(A). The generic fibers of ${\cal C}\mapsto B$ are smooth rational curves 
${\bf P}^1$.  

\medskip

(B). The fiberwise self-intersection number 
$\pi_{\ast}({\cal C}\cdot {\cal C})\in {\cal A}_{top}(B)$ is negative, $=-n, n\in
 {\bf N}$.

 Recall the concept of formal excess base dimension, $febd(C, {\cal X}/B)$,
 which is
 an algebraic numerical invariant assoicated to the class $C$ and
 the family ${\cal X}\mapsto B$, giving the formal base dimension to be
 thickened.

 For the two different $(1, 1)$ classes $C$ and $C+kPD({\cal C})$, 
$k\in {\bf Z}$,  their $ febd(C, {\cal X}/B)$ may not be the same.

 By switching the roles of $C$ and $C+kPD({\cal C})$, we may assume that $k>0$.

 We have the following proposition bounding $febd(C, {\cal X}/B)$,

\medskip

\begin{prop}\label{prop; decrease}
 Given two fiberwise $(1, 1)$ classes $C$ and $C+kPD({\cal C})$ with $k>0$,
we always have 

 $$febd(C, {\cal X}/B)\leq febd(C+kPD({\cal C}), {\cal X}/B)$$
\end{prop}

\medskip

\noindent Proof of the proposition: As in [Liu3], denote ${\cal T}_B({\cal X})$ 
 to be the relative $Pic^0$ variety of $\pi: {\cal X}\mapsto B$.
 Let ${\cal E}_C$ denote the locally free sheaf over ${\cal X}\times_B
{\cal T}_B({\cal X})$
 whose first Chern class is equal to the image of $C\in H^2({\cal X}, {\bf Z})$
 under 

$$H^2({\cal X}, {\bf Z})\mapsto H^2({\cal X}
\times_B{\cal T}_B({\cal X}), {\bf Z}).$$

 Then the coherent sheaf 
${\cal R}^0\pi_{\ast}\bigl( {\cal E}_C \bigr)$ determines an algebraic cone
 over ${\cal T}_B({\cal X})$ and the projectified cone, denoted by 
 ${\cal M}_C$, is the algebraic family moduli space associated to $C$.

 Let $D_C$ be the universal 
effective curve over ${\cal M}_C$ poincare dual to $C$.

 Then we have the following short exact sequence,
 
$$0\mapsto {\cal O}_{{\cal X}\times_B {\cal M}_C}(D_C)\mapsto 
{\cal O}_{{\cal X}\times_B {\cal M}_C}
(D_C+k{\cal C})\mapsto {\cal O}_{k{\cal C}\times_B {\cal M}_C}(D_C+
 k{\cal C})\mapsto 0.$$

  Push forward along 
$\pi_{{\cal M}_C}:{\cal X}\times_B {\cal M}_C\mapsto {\cal M}_C$,
and we have the following commutative diagram,

\[
\begin{array}{ccc}
 {\cal R}^2(\pi_{{\cal M}_C})_{\ast} {\cal O}_{{\cal X}\times_B {\cal M}_C}
  & \longrightarrow  & 
  {\cal R}^2(\pi_{{\cal M}_C})_{\ast} {\cal O}_{{\cal X}\times_B {\cal M}_C}(D_C)\\
 \Big\downarrow  &   &  \Big\downarrow  \\
   {\cal R}^2(\pi_{{\cal M}_C})_{\ast} {\cal O}_{{\cal X}\times_B {\cal M}_C}
 &  \longrightarrow & {\cal R}^2(\pi_{{\cal M}_C})_{\ast}{\cal O}_{{\cal X}\times_B
 {\cal M}_C}(D_C+k{\cal C}) \\
  & & \Big\downarrow   \\
  & &  0  \\
\end{array}
\]

  Recall that $febd(C, {\cal X}/B)$ 
is the rank of the maximal trivial invertible 
 subsheaf of 
${\cal R}^2(\pi_{{\cal M}_C})_{\ast} {\cal O}_{{\cal X}\times_B {\cal M}_C}$
 annihilated by the sheaf morphism

$${\cal R}^2(\pi_{{\cal M}_C})_{\ast} {\cal O}_{{\cal X}\times_B {\cal M}_C}
  \longrightarrow 
  {\cal R}^2(\pi_{{\cal M}_C})_{\ast} {\cal O}_{{\cal X}\times_B {\cal M}_C}(D_C),$$

  it is obvious from the commutative diagram that
$febd(C, {\cal X}/B)\leq febd(C+kPD({\cal C}), {\cal X}/B)$. $\Box$

\medskip

\begin{rem}\label{rem; subset}
  From the surjective morphism 
$${\cal R}^2(\pi_{{\cal M}_C})_{\ast} {\cal O}_{{\cal X}\times_B {\cal M}_C}(D_C)
\mapsto 
{\cal R}^2(\pi_{{\cal M}_C})_{\ast} {\cal O}_{{\cal X}\times_B {\cal M}_C}(D_C+
k{\cal C})\mapsto 0,$$

it is apparent that the support of the former coherent 
sheaf contains the support of
 the latter coherent sheaf as a subset.
\end{rem}

Let us state and prove the main theorem of this section,

 Recall the definitions as in [Liu3]
 that the fiber bundle $\pi:{\cal X}\mapsto B$ to be
 one-relatively good if there exists an effective relatively ample $D\subset
 {\cal X}\mapsto B$ such that $D\mapsto B$ is of relative dimension one.
 The fiber bundle $\pi:{\cal X}\mapsto B$ is said to be
 two-relatively good if there exists  effective 
 relative ample $D_1, D_2\mapsto B$
 such that (i). $D_1\mapsto B$ is of relative dimension one. (ii). $D_1\cap D_2
 \mapsto B$ is of relative dimension zero over $B$. The purpose of introducing
 such condition is to ensure the existence of algebraic family Kuranishi models.

\medskip

\begin{theo}\label{theo; alsw}
 Let $\pi:{\cal X}\mapsto B$ be an algebraic fiber bundle of algebraic surfaces 
 over a proper and smooth algebraic manifold $B$. Let $C$ be a $(1, 1)$ class
 of ${\cal X}$ which restricts to $(1, 1)$ classes of the fibers.

 Let ${\cal C}\subset {\cal X}$ be a relative curve over $B$ which satisfies
 condition (A) and (B) on page \pageref{condition}.

 Suppose the following conditions hold for $C$, $C+k{\cal C}$, $k\in {\bf N}$ and
 ${\cal X}\mapsto B$,

\medskip

\noindent $AF(i).$ ${C\cdot C-c_1({\bf K}_{{\cal X}/B})\cdot C\over 2}+dim_{\bf C}B+
febd(C, {\cal X}/B)\geq 0$.

\medskip

\noindent 
$AF(ii).$ ${(C+kPD({\cal C}))\cdot (C+kPD({\cal C}))-c_1({\bf K}_{{\cal X}/B})\cdot 
 (C+kPD({\cal C}))\over 2}+dim_{\bf C}B+febd(C+kPD({\cal C}), {\cal X}/B)\geq 0$.

\medskip

\noindent 
$AF(iii).$ If the fiber geometric genus $p_g>0$, 
the formal excess dimensions of $C$ and $C+kPD({\cal C})$ are equal.
 Namely, $febd(C, {\cal X}/B)=febd(C+kPD({\cal C}), {\cal X}/B)$.

If $p_g=0$, the supports of the coherent sheaves 
 ${\cal R}^2\pi_{\ast}\bigl({\cal E}_C\bigr)$, 
${\cal R}^2\pi_{\ast}\bigl({\cal E}_{C+kPD({\cal C})}\bigr)$ coincide.

\medskip

\noindent $AF(iv).$ Suppose $febd(C, {\cal X}/B)$ is equal to the maximum value, the
 geometric genus of the fibers $\pi^{-1}(b), b\in B$, then the fiber bundle
 $\pi: {\cal X}\mapsto B$ is one-relatively good.

 Otherwise, when $febd(C, {\cal X}/B)<p_g(\pi^{-1}(b)), b\in B$, the
 fiber bundle $\pi:{\cal X}\mapsto B$ is required to be two-relatively good.

\medskip

 Under these assumptions  $AF(i).-AF(iv).$, the pure 
algebraic family invariants of $C+kPD({\cal C})$
 and the mixed family invariant of $C$ are related by the following formula,

 $${\cal AFSW}_{{\cal X}\mapsto B}(1, C+kPD({\cal C}))=\sum_{0\leq i<\infty}
 {\cal AFSW}_{{\cal X}\mapsto B}(c_i({\cal V}_{1\mapsto k}), C),$$

 where the relative obstruction virtual sheaf ${\cal V}_{1\mapsto k}$ can be
 identified with the virtual sheaf
 $${\cal R}^1\pi_{\ast}{\cal O}_{k{\cal C}}(D_C+k{\cal C})-{\cal R}^0\pi_{\ast}
 {\cal O}_{k{\cal C}}(D_C+k{\cal C}).$$
\end{theo}

\medskip

\noindent Proof: As in [Liu3], we consider ${\cal E}_C$ to be the
 locally free sheaf over ${\cal X}_B\times_B {\cal T}_B({\cal X})$ (under the
 convention that ${\cal T}_B({\cal X})=B$ if $q(\pi^{-1}(b))=0, q\in B$)
 with first Chern class equal to the image of $C\in H^2({\cal X}, {\bf Z})$ under
 $H^2({\cal X}, {\bf Z})\mapsto 
H^2({\cal X}\times_B {\cal T}_B({\cal X}), {\bf Z})$. Suppose that 
${\cal R}^0\pi_{\ast}\bigl( {\cal E}_C\bigr)$ is the zero sheaf over
 ${\cal T}_C({\cal X})$, then $C$ is not represented as an effective curve within
 the algebraic family ${\cal X}\mapsto B$. If ${\cal R}^0\pi_{\ast}
\bigl({\cal E}_{C+kPD({\cal C})}\bigr)$ is the zero sheaf as well, then
it is easy to see that the
 algebraic mixed invariants
 ${\cal FASW}_{{\cal X}\mapsto B}(\eta, C)=0$,
 ${\cal FASW}_{{\cal X}\mapsto B}(\eta, C+kPD({\cal C}))=0$, for all 
$\eta\in {\cal A}_{\cdot}(B)$. Under this assumption, the family switching
 formula becomes a trivial identity. Thus, one may assume 
 ${\cal R}^0\pi_{\ast}\bigl({\cal E}_{C+kPD({\cal C})}\bigr)$ is non-trivial.

 In the following we assume ${\cal M}_{C+kPD({\cal C})}$, 
the projectified cone formed
 by ${\cal R}^0\pi_{\ast}\bigl({\cal E}_{C+kPD({\cal C})}\bigr)$ to be non-empty.  
  According to the discussion in [Liu3], one may separate into three cases, 

\medskip

(i). Either $p_g(\pi^{-1}(b))=0, b\in B$ and
 ${\cal R}^2(\pi)_{\ast}\bigl({\cal E}_{C+
kPD({\cal C})}\bigr)$ supports over a Zariski closed 
$Z\subset {\cal T}_B({\cal X})$, $Z\not=\emptyset$. Then there exists
 an algebraic family Kurnanishi model of $C+kPD({\cal C})$ over the
 Zariski open set $U=Z^c$,
$({\cal V}, {\cal W}, \Phi_{\cal VW})$, ${\cal V}$, ${\cal W}$ locally free
 over $U$,

$$\Phi_{\cal VW}: {\cal V}\mapsto {\cal W},$$ a sheaf morphism with

$Ker(\Phi_{\cal VW})\cong {\cal R}^0(\pi)_{\ast}\bigl({\cal E}_{C+kPD({\cal C})}
\bigr)|_U$ and $Coker(\Phi_{\cal VW})\cong 
{\cal R}^1(\pi)_{\ast}\bigl({\cal E}_{C+kPD({\cal C})}\bigr)|_U$.

\medskip

(ii). Or $p_g(\pi^{-1}(b))\geq 0, b\in B$ and 
 ${\cal R}^2(\pi)_{\ast}\bigl({\cal E}_{C+
kPD({\cal C})}\bigr)$ vanishes.  Then there exists an algebraic family
 Kuranshi model $({\cal V}, {\cal W}, \Phi_{\cal VW})$ 
over ${\cal T}_B({\cal X})$, ${\cal V}$, ${\cal W}$ locally
 free and the sheaf morphism $\Phi_{\cal VW}:{\cal V}\mapsto {\cal W}$.

\medskip

(iii). Or $p_g(\pi^{-1}(b))>0, b\in B$, and 
 the coherent sheaf ${\cal R}^2(\pi)_{\ast}\bigl({\cal E}_{C+
kPD({\cal C})}\bigr)$ is non-vanishing. In such situation, there exists
 an algebraic family Kuranishi model $({\cal V}, 
 {\cal W}, \Phi_{\cal VW})$, with algebraic locally free 
 ${\cal V}, {\cal W}$ over ${\cal T}_B({\cal X})$. The locally free
 ${\cal W}$ is constructed by a pair of locally free $\tilde{\cal V},
 \tilde{\cal W}$.

 For the details, please consult [Liu3] section 3.

\medskip

\begin{prop}\label{prop; cal}
 Let $\underline{C}$ be a $(1, 1)$ class of $\pi:{\cal X}\mapsto B$ which
 restricts to fiberwise $(1, 1)$ class. In the case $(i)$ above, 
let $U\subset {\cal T}_B({\cal X})$
 be the Zariski open set which is the complement of the support of 
 ${\cal R}^2\pi_{\ast}\bigl({\cal E}_{\underline{C}}\bigr)$. Otherwise, take
 $U={\cal T}_B({\cal X})$.
Then the pure algebraic family 
 Seiberg-Witten invariant ${\cal AFSW}_{{\cal X}\mapsto B}(1, \underline{C})$
 is equal to $c_{dim_{\bf C}B+q}({\cal W}-{\bf V})\in 
{\cal A}_0(U)$. Likewise, the mixed family
 invariant ${\cal AFSW}_{{\cal X}\mapsto B}(\eta, \underline{C})$, $\eta
 in {\cal A}_{\ast}(B)$ can
be identified with $\eta\cap c_{dim_{\bf C}B+q-deg(\eta)}({\cal W}-{\cal V})
\in {\cal A}_0(U)$.
\end{prop}

\medskip

\noindent Proof: As in [Liu3], we follow the convention that the bold 
 characters ${\bf V}, {\bf W}$ denote the algebraic vector bundle associated 
with the locally free sheaves ${\cal V}$, ${\cal W}$, etc.
 Recall that the family invariant 
${\cal FASW}_{{\cal X}\mapsto B}(1, \underline{C})$ is defined to be

$$c_1^{\cap{\underline{C}^2-\underline{C}\cdot c_1({\bf K}_{{\cal X}/B})\over 2}
+dim_{\bf C}B+p_g}({\bf H})\cap c_{rank_{\bf C}{\bf W}}
({\bf H}\otimes \pi^{\ast}_{{\bf P}_U({\bf V})}
{\bf W})\in{\cal A}_0({\bf P}_U({\bf V}))$$

 in both the (i) and (ii) cases, with $U=Z^c$, 
$Z=supp({\cal R}^2(\pi)_{\ast}\bigl({\cal E}_{C+kPD({\cal C})}\bigr))$,
 $p_g=p_g(\pi^{-1}(b))=0, b\in B$
 in the  case $(i)$, $U={\cal T}_B({\cal X})$ in the case $(ii)$.

\medskip

 In case $(iii)$, the algebraic family Seiberg-Witten invariant is defined to be
 
$$\hskip -1.2in 
{\cal FASW}_{{\cal X}/B}(1, \underline{C})=
c_1^{\cap{\underline{C}^2-\underline{C}\cdot c_1({\bf K}_{{\cal X}/B})\over 2}+
dim_{\bf C}B+febd(\underline{C}, {\cal X}/B)}
({\bf H})\cap c_{rank_{\bf C}{\bf W}}({\bf H}\otimes
 \pi^{\ast}_{{\bf P}_{{\cal T}_B({\cal X})}({\bf V})}{\bf W})\in 
 {\cal A}_0({\bf P}_{{\cal T}_B({\cal X})}({\bf V})).$$

  In this case, we may take $U={\cal T}_B({\cal X})$.

  Recall from [F] page 71 chap. 4 that for a cone ${\bf C}_X$
 over a scheme $X$,
  ${\bf q}:{\bf P}({\bf C}_X\oplus 1)\mapsto X$ denotes the projection map, then
 the Segre class of ${\bf C}_X$ 
is defined by $q_{\ast}(\sum_{i\geq 0} c_1({\cal O}(1))^i
\cap {\bf P}({\bf C}_X\oplus 1))\in {\cal A}_{\ast}(X)$.

 Applying to our situation, we take ${\bf C}_X$ to be the vector bundle cone formed
 by ${\bf V}$ and take ${\bf q}=\pi_{{\bf P}_{U}({\bf V})}
:{\bf P}_{{\cal T}_B({\cal X})}({\bf V})\mapsto U$. We also have
 $c_1({\cal O}(1))=c_1({\bf H})$. 

 By using the top Chern class identity $$c_{top}({\bf H}\otimes q^{\ast}{\bf W})
=\sum_{j\geq 0} c_1({\bf H})^j\cap c_{rank_{\bf C}{\bf W}-j}(q^{\ast}{\bf W})$$
 Then in case (i) and (ii), we have ${\cal AFSW}_{{\cal X}/B}(1, \underline{C})=$

 $${\bf q}_{\ast}\bigl(
c_1({\bf H})^{{\underline{C}^2-\underline{C}\cdot c_1({\bf K}_{{\cal X}/B})\over 2}
+dim_{\bf C}B+p_g}\cap c_{rank_{\bf C}{\bf W}}
({\bf H}\otimes {\bf q}^{\ast}
{\bf W})\cap [{\bf P}_U({\bf V})]\bigr)\in {\cal A}_0(U)$$

$$=\sum_{j\geq 0}{\bf q}_{\ast}\bigl( c_1({\bf H})^{{\underline{C}^2-\underline{C}
\cdot c_1({\bf K}_{{\cal X}/B})\over 2}
+dim_{\bf C}B+p_g}\cap c_1({\bf H})^j\cap 
c_{rank_{\bf C}{\bf W}}(q^{\ast}{\bf W})\cap [{\bf P}_U({\bf V})]\bigr)$$

$$=\sum_{k\geq {\underline{C}^2-\underline{C} \cdot 
c_1({\bf K}_{{\cal X}/B})\over 2}+dim_{\bf C}B+p_g} {\bf q}_{\ast}\bigl( 
 c_1({\bf H})^k\cap [{\bf P}_U({\bf V})]\bigr)
\cap c_{rank_{\bf C}{\bf W}+{\underline{C}^2-\underline{C} \cdot 
c_1({\bf K}_{{\cal X}/B})\over 2}+dim_{\bf C}B+p_g-k}({\bf W})$$

$$=\sum_{l\geq {\underline{C}^2-\underline{C} \cdot 
c_1({\bf K}_{{\cal X}/B})\over 2}+dim_{\bf C}B+p_g-rank_{\bf C}{\bf V}+1}
 s_l({\bf V})\cap c_{rank_{\bf C}{\bf W}+{\underline{C}^2-\underline{C} \cdot 
c_1({\bf K}_{{\cal X}/B})\over 2}+dim_{\bf C}B+p_g-l-rank_{\bf C}{\bf V}+1}
({\bf W}).$$

 By using the identity 
$${\underline{C}^2-\underline{C}\cdot c_1({\bf K}_{{\cal X}/B})\over 2}
+p_g-q+1=rank_{\bf C}{\bf V}-rank_{\bf C}{\bf W},$$

 we may simplify the above summation into 

$$\sum_{l\geq dim_{\bf C}B+q-rank_{\bf C}{\bf W}}s_l({\bf V})\cap
 c_{dim_{\bf C}B+q-l}({\bf W})
 =\sum_{m\leq rank_{\bf C}{\bf W}}s_{dim_{\bf C}B+q-m}({\bf V})c_m({\bf W})$$

$$=c_{dim_{\bf C}B+q}({\bf W}-{\bf V})=c_{dim_{\bf C}B+q}({\cal W}-{\cal V})
\in{\cal A}_0(U).$$

Likewise, the enumeration of the mixed invariant 
${\cal AFSW}_{{\cal X}\mapsto B}(\eta, \underline{C})$ is almost 
 identical. The insertion of $\eta$ effectively drops $dim_{\bf C}B$ to
 $dim_{\bf C}B-deg(\eta)$ and the answer becomes 
 $\eta\cap c_{dim_{\bf C}B+q-deg(\eta)}({\cal W}-{\cal V})\in {\cal A}_0(U)$.

 This finishes the proof of the (i). and (ii). cases.

 For (iii)., the calculation of the pure invariant or the
 $\eta$ inserted mixed invariant is almost identical. The major difference
 is that we replace $p_g$ by the formal excess base dimension 
$febd(C, {\cal X}/B)$ and we use the identity

$${\underline{C}^2-\underline{C}\cdot c_1({\bf K}_{{\cal X}/B})\over 2}
+febd(C, {\cal X}/B)-q+1=rank_{\bf C}{\bf V}-rank_{\bf C}{\bf W}$$

in reducing the intersection number. $\Box$

\medskip

 Consider the invertible sheaves ${\cal E}_{C}$, ${\cal E}_{C+kPD({\cal C})}$ with
 first Chern class $C$, $C+kPD({\cal C})$. By tensoring 
${\cal E}_{C+kPD({\cal C})}$ to the short exact sequence 

$$0\mapsto {\cal O}_{{\cal X}\times_B{\cal T}_B({\cal X})}(-k{\cal C})\mapsto
 {\cal O}_{{\cal X}\times_B{\cal T}_B({\cal X})}\mapsto 
{\cal O}_{k{\cal C}\times_B {\cal T}_B({\cal X})}\mapsto 0,$$

 one gets the following short exact sequence 

$$0\mapsto {\cal E}_C\mapsto {\cal E}_{C+kPD({\cal C})}\mapsto
{\cal O}_{k{\cal C}\times_B {\cal T}_B({\cal X})}\otimes {\cal E}_{C+kPD({\cal C})}
\mapsto 0.$$

By taking the direc images along 
$\pi:{\cal X}\times_B {\cal T}_B({\cal X})\mapsto {\cal T}_B({\cal X})$, 
one gets the following identities in the $K_0$ group of coherent sheaves
 on $ {\cal T}_B({\cal X})$,

$$[\pi_{\ast}\bigl({\cal E}_{C+kPD({\cal C})}\bigr)]=[\pi_{\ast}\bigl(
{\cal E}_C\bigr)]+[\pi_{\ast}\bigl(
{\cal O}_{k{\cal C}\times_B {\cal T}_B({\cal X})}\otimes 
{\cal E}_{C+kPD({\cal C})}\bigr)],$$

where $\pi_{\ast}\bigl({\cal E}_C\bigr)$ denotes 

$$\sum_{0\leq i\leq 2}(-1)^i{\cal R}^i\pi_{\ast}\bigl({\cal E}_C\bigr),$$

etc.

\medskip

\begin{prop}\label{prop; alternating}
 Let $\underline{C}$ be a fiberwise monodromy invariant 
$(1, 1)$ class of ${\cal X}\mapsto B$ and let 
 $({\cal V}, {\cal W}, \Phi_{{\cal V}{\cal W}})$ be an algebraic family
 Kuranishi model of ${\cal M}_{\underline{C}}$ over a Zariski open set
 $U\subset {\cal T}_B({\cal X})$, then there exists an
 equality in the $K$ group of coherent sheaves on $U$ when
 $p_g=0$ or ${\cal R}^2\pi_{\ast}\bigl({\cal E}_{\underline{C}}\bigr)=0$,

$$\hskip -.9in
[{\cal V}-{\cal W}]=[\pi_{\ast}\bigl({\cal E}_{\underline{C}}\bigr)|_U].$$

There exists an equality in the reduced $K$ group of coherent sheaves on
 $U$ when $p_g>0$ and 
${\cal R}^2\pi_{\ast}\bigl({\cal E}_{\underline{C}}\bigr)\not=0$, 

$$\hskip -.9in
[{\cal V}-{\cal W}]=[\pi_{\ast}\bigl({\cal E}_{\underline{C}}\bigr)|_U]
-[{\cal R}^2\pi_{\ast}\bigl({\cal O}_{{\cal X}\times_B{\cal T}_B({\cal X})}
\bigr)|_U].$$
\end{prop}

\medskip

\noindent Proof of the proposition: We discuss the (i)., (ii). cases
 first. Recall that the open set $U\subset {\cal T}_B({\cal X})$ is
 defined to be the complement of the support of the coherent sheaf
 ${\cal R}^2\pi_{\ast}\bigl({\cal E}_{\underline{C}}\bigr)$, then
 $\pi_{\ast}\bigl({\cal E}_{\underline{C}}\bigr)$ is reduced to

$${\cal R}^0\pi_{\ast}\bigl({\cal E}_{\underline{C}}\bigr)
-{\cal R}^1\pi_{\ast}\bigl({\cal E}_{\underline{C}}\bigr).$$

 On the other hand, for both $(i).$ and $(ii).$, we have

$$-\hskip .9in 
0\mapsto {\cal R}^0\pi_{\ast}\bigl({\cal E}_{\underline{C}}\bigr)|_U
{\cal V}\mapsto {\cal W}\mapsto {\cal R}^1\pi_{\ast}
\bigl({\cal E}_{\underline{C}}\bigr)|_U\mapsto 0,$$

 by the defining property of the algebraic family Kuranishi model, then
 we have the identity

$$[{\cal V}-{\cal W}]=[\pi_{\ast}\bigl({\cal E}_{\underline{C}}\bigr)|_C],$$

in the $K_0(U)$.

  For the case when 
${\cal R}^2\pi_{\ast}\bigl({\cal E}_{\underline{C}}\bigr)\not=0$ 
and $p_g>0$, we have the following equalities instead,

$$[{\cal V}]+[\tilde{\cal V}]-[\tilde{\cal W}]=[\pi_{\ast}
\bigl({\cal E}_{\underline{C}}\bigr)].$$
$$[{\cal V}-{\cal W}]=[{\cal V}+\tilde{\cal V}]-[\tilde{\cal W}]-[{\cal F}].$$

where $[{\cal F}]=[{\cal}^2\pi_{\ast}\bigl({\cal O}_{{\cal X}\times_B
 {\cal T}_B({\cal X})}\bigr)]-
 [{\cal O}_{{\cal T}_B({\cal X})}^{febd(\underline{C}, {\cal X}/B)}]$. 
For the details, please consult [Liu3], section 4. 
 
Then up to the trivial factor 
$[{\cal O}_{{\cal T}_B({\cal X})}^{febd(\underline{C}, {\cal X}/B)}]$,
 we have the following equality

$$[{\cal V}-{\cal W}]=[\pi_{\ast}\bigl({\cal E}_{\underline{C}}\bigr)]
-[{\cal}^2\pi_{\ast}\bigl({\cal O}_{{\cal X}\times_B
 {\cal T}_B({\cal X})}\bigr)]$$

in the reduced $K$ group of ${\cal T}_B({\cal X})$. This ends the proof of
 proposition \ref{prop; alternating}.
$\Box$

\medskip

 Apply the proposition \ref{prop; cal} and proposition 
 \ref{prop; alternating} to $\underline{C}=C+kPD({\cal C})$, 
 we find that ${\cal AFSW}_{{\cal X}\mapsto B}(1, C+kPD({\cal C}))$ is 
 equal to
 
 $$=c_{q+dim_{\bf C}B}([\pi_{\ast}\bigl({\cal E}_{C+kPD({\cal C})}\bigr)]-
 [{\cal R}^2\pi_{\ast}\bigl({\cal O}_{{\cal X}\times_B
 {\cal T}_B({\cal X})}\bigr)])$$

 when $febd(C+kPD({\cal C}), {\cal X}/B)<p_g$.

 The above chern class can be re-written as 

$$\hskip -.9in c_{dim_{\bf C}B+q}([\pi_{\ast}\bigl({\cal E}_{C}\bigr)]+
 [\pi_{\ast}\bigl({\cal O}_{k{\cal C}\times_B {\cal X}}
 \times {\cal E}_{C+kPD({\cal C})}\bigr)]-
[{\cal R}^2\pi_{\ast}\bigl({\cal O}_{{\cal X}\times_B
 {\cal T}_B({\cal X})}\bigr)])$$

$$\hskip -.9in \sum_{r\geq 0}
c_r([\pi_{\ast}\bigl({\cal O}_{k{\cal C}\times_B {\cal X}}
 \times {\cal E}_{C+kPD({\cal C})}\bigr)])
\cap c_{dim_{\bf C}B+q-r}([\pi_{\ast}\bigl({\cal E}_{C}\bigr)|_U]-
[{\cal R}^2\pi_{\ast}\bigl({\cal O}_{{\cal X}\times_B
 {\cal T}_B({\cal X})}\bigr)]).$$

 According to the assumption $AF(iii).$, we know that
 $0\leq febd(C, {\cal X}/B)=febd(C+kPD({\cal C}), {\cal X}/B)<p_g$, 
 then apply proposition \ref{prop; cal} and proposition \ref{prop; alternating}
 to $\underline{C}=C$, we find that for 
$\eta_r=c_r([\pi_{\ast}\bigl({\cal O}_{k{\cal C}\times_B {\cal X}}
 \times {\cal E}_{C+kPD({\cal C})}\bigr)])$, the expression

$$c_r([\pi_{\ast}\bigl({\cal O}_{k{\cal C}\times_B {\cal X}}
 \times {\cal E}_{C+kPD({\cal C})}\bigr)])
\cap c_{dim_{\bf C}B+q-r}([\pi_{\ast}\bigl({\cal E}_{C}\bigr)]-
[{\cal R}^2\pi_{\ast}\bigl({\cal O}_{{\cal X}\times_B
 {\cal T}_B({\cal X})}\bigr)])$$

 is nothing but the mixed invariant ${\cal AFSW}_{{\cal X}\mapsto B}(\eta_r, C)$.
By combining all the identities, we get

$$\hskip -.8in {\cal AFSW}_{{\cal X}\mapsto B}(1, C+kPD({\cal C}))=
 \sum_{r\geq 0}{\cal AFSW}_{{\cal X}\mapsto B}(
c_r([\pi_{\ast}\bigl({\cal O}_{k{\cal C}\times_B {\cal X}}
 \times {\cal E}_{C+kPD({\cal C})}\bigr)]), C).$$

 This proves the family switching formula when 
 $0\leq febd(C+kPD({\cal C}))<p_g$.

\medskip

  On the other hand, when $febd(C+kPD({\cal C}), {\cal X}/B)=p_g$,
 either ${\cal R}^2\pi_{\ast}\bigl({\cal E}_{C+kPD({\cal C})}\bigr)=0$ 
 or ${\cal R}^2\pi_{\ast}\bigl({\cal O}_{{\cal X}\times_B
 {\cal T}_B({\cal X})} \bigr)\cong {\cal O}_{{\cal T}_B({\cal X})}^{p_g}$
 and ${\cal F}=0$.  When $p_g=0$, either 
${\cal R}^2\pi_{\ast}\bigl({\cal E}_{C+kPD({\cal C})}\bigr)=0$
 or the algebraic Kuranishi model of $C+kPD({\cal C})$ is built above the
 complement of the coherent sheaf 
${\cal R}^2\pi_{\ast}\bigl({\cal E}_{C+kPD({\cal C})}\bigr)=0$.

 In all these cases, 
${\cal AFSW}_{{\cal X}\mapsto B}(1, C+kPD({\cal C}))$
 has the following schematic form,
 
$$=c_{q+dim_{\bf C}B}([\pi_{\ast}\bigl({\cal E}_{C+kPD({\cal C})}\bigr)|_U])$$

$$=c_{q+dim_{\bf C}B}([\pi_{\ast}\bigl({\cal E}_{C}\bigr)|_U]+
 [\pi_{\ast}\bigl({\cal O}_{k{\cal C}\times_B {\cal X}}
 \times {\cal E}_{C+kPD({\cal C})}\bigr)|_U])$$

$$=\sum_{r\geq 0}c_r([\pi_{\ast}\bigl({\cal O}_{k{\cal C}\times_B {\cal X}}
 \times {\cal E}_{C+kPD({\cal C})}\bigr)|_U])\cap 
 c_{q+dim_{\bf C}B-r}([\pi_{\ast}\bigl({\cal E}_{C}\bigr)|_U]),$$

 When $p_g>0$, $U$ is taken to be the whole space ${\cal T}_B({\cal X})$.
 When $p_g=0$, $U$ is taken to be the complement of the support of the sheaf
 ${\cal R}^2\pi_{\ast}\bigl({\cal E}_{C+kPD({\cal C})}\bigr)$ in
 ${\cal T}_B({\cal X})$.

 By the assumption $AF(iii).$ that $febd(C, {\cal X}/B)=febd(C+kPD({\cal C}), 
 {\cal X}/B)=p_g>0$ or 
the support of ${\cal R}^2\pi_{\ast}\bigl({\cal E}_C\bigr)=U^c\subset
 {\cal T}_B({\cal X})$, the above summation can be recasted into

$$=\sum_{r\geq 0}{\cal AFSW}_{{\cal X}\mapsto B}(
c_r([\pi_{\ast}\bigl({\cal O}_{k{\cal C}\times_B {\cal X}}
 \times {\cal E}_{C+kPD({\cal C})}\bigr)]), C).$$

 This finishes the proof of the theorem.
 $\Box$

\medskip

\begin{rem}\label{rem; clarify}
 The condition $AF(iii)$ has no counterpart in the $C^{\infty}$ proof of the
 family switching formula. The condition on the formal excess base dimension
 when $p_g>0$ is needed as ${\cal AFSW}$ differ from the usual $FSW$ for
 families of $p_g>0$ algebraic surfaces. The condition on the supports of 
 ${\cal R}^2\pi_{\ast}\bigl({\cal E}\bigr)$ is needed when the relative 
curve ${\cal C}\mapsto B$ has singular fibers. When ${\cal C}\mapsto B$ is
 smooth (as was assumed in the $C^{\infty}$ version of 
family switching formula), 
one can show that the supports of ${\cal R}^0\pi_{\ast}\bigl(
{\cal E}_{C+kPD({\cal C})}\bigr)$ and of ${\cal R}^2\pi_{\ast}\bigl(
{\cal E}_{C}\bigr)$ do not overlap. By shrinking the open set $U$ to
 avoid the support of ${\cal R}^2\pi_{\ast}\bigl({\cal E}_C\bigr)$, one
 may prove the family switching formula in the case
 $p_g=0$, ${\cal R}^2\pi_{\ast}\bigl({\cal E}_C\bigr)\not=0$, without assuming the 
supports of  ${\cal R}^2\pi_{\ast}\bigl({\cal E}_C\bigr)$ and 
${\cal R}^2\pi_{\ast}\bigl({\cal E}_{C+kPD({\cal C})}\bigr)$ coincide.
\end{rem}

\medskip

 In the statement of the family switching formula of algebraic family 
Seiberg-Witten invariants, the various Chern classes of the
 relative obstruction virtual 
sheaf $\pi_{\ast}\bigl({\cal O}_{k{\cal C}\times_B 
{\cal T}_B({\cal X})}\otimes {\cal E}_{C+kPD({\cal C})}\bigr)$ 
$$={\cal R}^0\pi_{\ast}\bigl({\cal O}_{k{\cal C}\times_B 
{\cal T}_B({\cal X})}\otimes {\cal E}_{C+kPD({\cal C})}\bigr)-
{\cal R}^1\pi_{\ast}\bigl({\cal O}_{k{\cal C}\times_B 
{\cal T}_B({\cal X})}\otimes {\cal E}_{C+kPD({\cal C})}\bigr)$$

 relates the algebraic family invariants. When the relative curve 
${\cal C}\mapsto B$ is not smooth, neither of
 ${\cal R}^i\pi_{\ast}\bigl({\cal O}_{k{\cal C}\times_B 
{\cal T}_B({\cal X})}\otimes {\cal E}_{C+kPD({\cal C})}\bigr)$, $i=1, 2$
 are locally free and the ranks 
 ${\cal R}^i\pi_{\ast}\bigl({\cal O}_{k{\cal C}\times_B 
{\cal T}_B({\cal X})}\otimes {\cal E}_{C+kPD({\cal C})}\bigr)\otimes 
 k(b)$ jump as $b$ range through the singular values of the map
${\cal C}\mapsto B$.

At the end of the section, we discuss its relationship with the
 relative obstruction bundle appearing in the $C^{\infty}$ version of
 the family switching formula.

Consider the short exact sequences

$$\hskip -.9in 
 0\mapsto {\cal O}_{(n-1){\cal C}}(-{\cal C})\otimes {\cal E}_{C+nPD({\cal C})}
\mapsto {\cal O}_{n{\cal C}}\otimes {\cal E}_{C+nPD({\cal C})}
\mapsto {\cal O}_{\cal C}\otimes {\cal E}_{C+nPD({\cal C})}\mapsto 0$$

for $n= 2,\cdots, k$.
which comes from tensoring ${\cal E}_{C+nPD({\cal C})}$ to

$$\hskip -.9in 
0\mapsto {\cal O}_B(-A)\mapsto {\cal O}_{A\cup B}\mapsto {\cal O}_A\mapsto 0$$

with $A={\cal C}$, $B=(n-1){\cal C}$ being the relative divisors in 
 $\pi:{\cal X}\mapsto B$.

\medskip

\begin{prop}\label{prop; induction}
 In the $K$ group of coherent sheaves on $B$, there is an equality
$$[\pi_{\ast}\bigl({\cal O}_{k{\cal C}\times_B 
{\cal T}_B({\cal X})}\otimes {\cal E}_{C+kPD({\cal C})}\bigr)]
=\sum_{1\leq p\leq k}[\pi_{\ast}\bigl({\cal O}_{\cal C}(p{\cal C})
\otimes {\cal E}_C\bigr)].$$ 
\end{prop}

\medskip

\noindent Proof of the proposition: For $k=1$, the formula is an identity
 once we realize ${\cal E}_{C+PD({\cal C})}={\cal O}({\cal C})\otimes {\cal E}_C$.

Assume that by induction hypothesis the
 equality has been proved for $k-1$, namely,

$$[\pi_{\ast}\bigl({\cal O}_{(k-1){\cal C}\times_B 
{\cal T}_B({\cal X})}\otimes {\cal E}_{C+(k-1)PD({\cal C})}\bigr)]
=\sum_{1\leq p\leq k-1}[\pi_{\ast}\bigl({\cal O}_{\cal C}(p{\cal C})
\otimes {\cal E}_C\bigr)].$$ 

By using the above short exact 
sequences and 
 ${\cal E}_{C+pPD({\cal C})}={\cal E}_C\otimes {\cal O}(p{\cal C})$, 
we get

$$[\pi_{\ast}\bigl({\cal O}_{k{\cal C}\times_B 
{\cal T}_B({\cal X})}\otimes {\cal E}_{C+kPD({\cal C})}\bigr)]
=[\pi_{\ast}\bigl({\cal O}_{(k-1){\cal C}\times_B 
{\cal T}_B({\cal X})}\otimes {\cal E}_{C+(k-1)PD({\cal C})}\bigr)]
+[\pi_{\ast}\bigl({\cal O}_{{\cal C}\times_B 
{\cal T}_B({\cal X})}(k{\cal C})\otimes {\cal E}_C\bigr)]$$

$$=\sum_{1\leq p\leq k-1}[\pi_{\ast}\bigl({\cal O}_{\cal C}(p{\cal C})
\otimes {\cal E}_C\bigr)]+[\pi_{\ast}\bigl({\cal O}_{{\cal C}\times_B 
{\cal T}_B({\cal X})}(k{\cal C})\otimes {\cal E}_C\bigr)]$$

$$=\sum_{1\leq p\leq k}[\pi_{\ast}\bigl({\cal O}_{\cal C}(p{\cal C})
\otimes {\cal E}_C\bigr)].$$ $\Box$

Because equivalent elements in the $K$ group have identical Total Chern
 classes, we may use 
$\sum_{1\leq p\leq k-1}[\pi_{\ast}\bigl({\cal O}_{\cal C}(p{\cal C})
\otimes {\cal E}_C\bigr)]$ in the family switching formula.

\medskip

 If ${\cal C}\mapsto B$ is smooth, then 
$\pi_{\ast}\bigl({\cal O}_{\cal C}(p{\cal C})
\otimes {\cal E}_C\bigr)$ is equal to
${\cal R}^0\pi_{\ast}\bigl({\cal O}_{\cal C}(p{\cal C})\otimes {\cal E}_C\bigr)$

 if $\int_{\cal C} (C+pPD({\cal C})\geq 0$;
it is equal to $0$ if $\int_{\cal C} (C+pPD({\cal C})=-1$;
 it is equal to $-{\cal R}^1\pi_{\ast}\bigl(
{\cal O}_{\cal C}(p{\cal C})\otimes {\cal E}_C\bigr)$ if
 $\int_{\cal C} (C+pPD({\cal C}))<-1$.

  Then under the identifications: 

$${\cal O}_{\cal C}({\cal C})\cong {\cal N}_{\cal C}$$

$${\cal E}^2_C\otimes {\cal K}^{-1}_{{\cal X}/B}\cong {\cal L}$$

 each term $\pi_{\ast}\bigl({\cal O}_{\cal C}(p{\cal C})\otimes {\cal E}_C\bigr)$
 is identified with the locally free sheaf ${\cal V}_p$ assoicated
 to the relative obstruction sheaf (see prop. \ref{prop; bundle}).

\bigskip

\section{\bf The Switching of the Canonical Obstruction Bundles}\label{section; 
 canon}

\bigskip

  As in the previous section let ${\cal X}\mapsto B$ be an algebraic fiber
bundle and let ${\cal X}\supset {\cal C}\mapsto B$ be an ${\bf P}^1$
 fibration with negative fiberwise self-intersection number
 $\int_{{\cal X}/B}PD({\cal C})^2<0$.

  The algebraic construction of the family switching formula suggests that when
 the $(1, 1)$ classes $C$ and $C+kPD({\cal C})$ differ by 
the $k$-multiple of the Poincare dual class, the algebraic family Seiberg-Witten 
invariants of $C$ and $C+kPD({\cal C})$
 are related to each other through the topological data determined
 by $C$, the multiplicity $k$ and ${\cal C}$.

  In fact, this is one of the key observations used in the proof of the 
G${\ddot o}$ttsche-Yau-Zaslow conjecture [Liu1]. 

  In this section, we construct the switching long exact sequence for
the algebraic canonical obstruction bundles of the universal families.

Let $M$ be an algebraic surface over ${\bf C}$. Let $M_n$ denote the
 $n-$th universal space constructed in [V],[Liu1].  Let $\Gamma$ be
 an $n$-vertexes admissible graph. Following the notation of [Liu1],
 let $Y(\Gamma)=\coprod_{\Gamma'<\Gamma} Y_{\Gamma'}\coprod Y_{\Gamma}$ denote the 
 closure of the locally closed subset $Y_{\Gamma}\subset M_n$ called the
 admissible strata.
 Then $B=Y(\Gamma)\subset M_n$ is the base space of the
 fiber bundle $\pi=f_n|_{Y(\Gamma)}:
 Y(\Gamma)\times_{M_n} M_{n+1}\mapsto Y(\Gamma)$.

 As a fiber product of $M_{n+1}\mapsto M_n$ and $Y(\Gamma)\subset M_n$, the
 space $Y(\Gamma)\times_{M_n} M_{n+1}$ can be constructed from $Y(\Gamma)\times M$
 by blowing up $n$ consecutive times.  Let $E_i$, $1\leq i\leq n$ denote the
 $i-th$ exceptional divisor of the blowing ups. Then $E_i$ induces a unique
 cohomology class in $H^{1, 1}(Y(\Gamma)\times_{M_n}M_{n+1}, {\bf Z})$. 
 On the other hand,
 any given $C\in H^{1, 1}(M, {\bf Z})$ also induces a class in 
 $H^{1, 1}(Y(\Gamma)\times_{M_n}M_{n+1}, {\bf Z})$. We
slightly abuse the notation and denote them by the same symbols $E_i$, $1\leq i
\leq n$ and $C$.

 Consider a topological type of algebraic curve singularity, 
let $m_i\in {\bf N}$ be the multiplicities of the minimal resolution of the  
curve singularities.

 Consider $\underline{C}=C-\sum_{1\leq i\leq n}m_i E_i$, also
 denoted as $C-{\bf M}(E)E$, following the convention in [Liu1].

 Over the generic strata $Y_{\Gamma}$ of the smooth space $Y(\Gamma)$, there are
 a finite number of irreducible smooth type $I$ exceptional curves in 
 the fibers of $Y_{\Gamma}\times_{M_n}M_{n+1}\mapsto Y_{\Gamma}$ whose
 cohomology classes 
 are the combinations of the various $E_i$, $e_i=E_i-\sum_{j}E_j$; where in 
the summation $\sum_j$ the 
$j-th$ 
vertexes run through all  
the direct descendents of the $i-th$ vertex in the graph $\Gamma$.
 All such $e_i$ generate a simplicial cone in $H^{1, 1}(Y(\Gamma)\times_{M_n}
M_{n+1}, {\bf Z})$.

 According to the general curve enumeration scheme in section 3 of [Liu4],
 there are a finite
 number of type $I$ exceptional classes $e_{k_i}=\sum E_{k_i}-\sum_{j_{k_i}} 
E_{j_{k_i}}$, $1\leq i\leq p$, 
 characterized by the property that 
 $(C-{\bf M}(E)E)\cdot e_{k_i}<0$ if and only if $1\leq i\leq p$. Such 
$e_{k_i}$, for $1\leq i\leq p$,
 generate a sub-simplicial cone, called the type $I$ exceptional cone
 of $\underline{C}=C-{\bf M}(E)E$ over $Y_{\Gamma}$.

 The complex codimension of $Y(\Gamma)\subset M_n$ is equal to
$-\sum_{1\leq i\leq n} {e_i^2-e_i\cdot c_1({\bf K}_{M_{n+1}/M_n})\over 2}$,
 which is also equal to the number of edges of the admissible graph $\Gamma$.

\medskip
  
  Each $e_{k_i}$ is represented by holomorphic curves above $Y(\Gamma)$ 
and it corresponds to a ${\bf P}^1$ fibration (may contain singular
 fibers) $\Xi_{k_i}\mapsto B$ such that

\medskip

\noindent (i). The generic fibers of $\Xi_{k_i}\mapsto B$
 are smooth ${\bf P}^1$. 

\medskip

\noindent (ii). The fiberwise self intersection number of $\Xi_{k_i}$ is 
equal to $e_{k_i}\cdot e_{k_i}<0$.

\bigskip

  Therefore, each of the curve
 fibrations $\Xi_{k_i}$, $1\leq i\leq p$ are qualified to 
be the ${\cal C}$ in the discussion of the algebraic family Seiberg-Witten 
switching formula in section \ref{section; apfsf}.

 The family switching formula concludes that the family invariants
 of $\underline{C}=C-{\bf M}(E)E$ and of
 $\underline{C}-e_{k_1}$, $\underline{C}-e_{k_1}-e_{k_2}$, $\cdots$,
 or $\underline{C}-\sum_{1\leq i\leq p}e_{k_i}$ are related to each other 
through the insertion of Chern classes of certain relative obstruction bundles.

  In the following, we assume additionally
 that the class $C$ is raised to a higher multiple such that
 $deg_{{\omega}_M}\bigl(C-c_1({\bf K}_M)\bigr)>0$.
 As a consequence of Serre duality, the sheaf 
 ${\cal R}^2\pi_{\ast}\bigl({\cal E}_C\bigr)=0$.

\medskip

\begin{lemm}\label{lemm; vanishing}
 Let ${\cal X}=Y(\Gamma)\times_{M_n}M_{n+1}\mapsto Y(\Gamma)=B$ be the
algebraic fiber bundle.
Suppose that 
$$deg_{{\omega}_{{\cal X}/B}}\bigl(C-c_1({\bf K}_{{\cal X}/B})\bigr)>0,$$

 then ${\cal R}^2\pi_{\ast}\bigl({\cal E}_{C-{\bf M}E(E)-\sum_{1\leq i\leq p}
e_{k_i}}\bigr)=0$, 
${\cal R}^2\pi_{\ast}\bigl({\cal E}_{C-{\bf M}(E)E}\bigr)=0$.
\end{lemm}

\medskip

\noindent Proof of lemma \ref{lemm; vanishing}:
 Both $C-{\bf M}(E)E$ and $C-{\bf M}E(E)-\sum_{1\leq i\leq p}
e_{k_i}$ can be schematically written as 
 $C+\sum_{1\leq i\leq n} d_iE_i$ for some tuples of $d_i\in {\bf Z}$.

 If either of the sheaves
${\cal R}^2\pi_{\ast}\bigl({\cal E}_{C-{\bf M}E(E)-\sum_{1\leq i\leq p}
e_{k_i}} \bigr)$, 
${\cal R}^2\pi_{\ast}\bigl({\cal E}_{C-{\bf M}(E)E}\bigr)$ is not zero, then
 there exists a $b\in T(M)\times B$ such that the base change  
 ${\cal R}^2\pi_{\ast}\bigl({\cal E}_{C+\sum_{1\leq i\leq n}
 d_i E_i}\bigr)\otimes k(b)$
 is nonzero. Then by base change theorem [Ha], the second sheaf cohomology 
 $H^2({\cal X}_b, {\cal E}_{C+\sum_i d_iE_i}|_b)$ is nonzero. 
 By Serre duality, it implies that ${\cal K}_{{\cal X}_b}\otimes
 {\cal E}^{\ast}_{C+\sum d_iE_i}$ has a non-trivial holomorphic section over
 ${\cal X}_b$.
 
 Thus by adjunction formula of the canonical class, 
the class $c_1({\bf K}_M)+\sum_{1\leq i\leq n} 
(1-d_i)E_i-C$ is represented by a holomorphic curve in ${\cal X}_b$.

  Let $\omega_0$ be an ample polarization on $M$.
  Because the smooth fiber algebraic surface ${\cal X}_b$ is blown up from $(M, 
 \omega_0)$ by
 $n$-consecutive blowing-ups (with the blowing up centers determined by
 the image of $b$ in $M_n$), the polarization class 
 on ${\cal X}_b$ can be chosen
 to be $\omega_0-\sum_{1\leq i\leq n}\epsilon_iE_i$ for
 some sequence of sufficiently (and arbitrarily) small 
$\epsilon_i, 1\leq i\leq n$, $0<\epsilon_i\mapsto 0$.

 Then the degree of the effective class $c_1({\bf K}_M)+\sum_{1\leq i\leq n} 
(1-d_i)E_i-C$,

$$(\omega_0-\sum_{1\leq i\leq n}\epsilon_iE_i)\cdot
 (c_1({\bf K}_M)+\sum_{1\leq i\leq n} 
(1-d_i)E_i-C)=(c_1({\bf K}_M)-C)\cdot \omega_0+\sum_{1\leq i\leq n}
 \epsilon_i(d_i-1)>0.$$

 However this implies that 
 $deg_{\omega_0}(c_1({\bf K}_M)-C)\geq lim\sum_{1\leq i\leq n}
 \epsilon_i(d_i-1)=0$ when we take $\epsilon_i\mapsto 0$. Contradicting to the
 assumption $deg_{\omega_0}(c_1({\bf K}_M)-C)=-deg_{\omega_0}(C-c_1({\bf K}_M))<0$.
 $\Box$

 In the long paper [Liu1], the relationship among the family invariants
 was read off by using the long exact sequence

$$\hskip -.9in 
0\mapsto {\cal R}^0\pi_{\ast}\bigl({\cal E}_{C-{\bf M}(E)E-\sum_{i}e_{k_i}}\bigr)
\mapsto {\cal R}^0\pi_{\ast}\bigl({\cal E}_{C-{\bf M}(E)E}\bigr)
\mapsto {\cal R}^0\pi_{\ast}\bigl({\cal O}_{\sum_{1\leq i\leq p} \Xi_{k_i}}
\otimes {\cal E}_{C-{\bf M}(E)E}\bigr)$$

$$\hskip -.9in
 \mapsto {\cal R}^1\pi_{\ast}\bigl({\cal E}_{C-{\bf M}(E)E-\sum_{i}e_{k_i}}\bigr)
\mapsto {\cal R}^1\pi_{\ast}\bigl({\cal E}_{C-{\bf M}(E)E}\bigr)
\mapsto {\cal R}^1\pi_{\ast}\bigl({\cal O}_{\sum_{1\leq i\leq p} \Xi_{k_i}}
\otimes {\cal E}_{C-{\bf M}(E)E}\bigr)\mapsto 0.$$

Following the discussion in section 5.1 of [Liu3], 
 canonical algebraic family Kuranishi models
 $({\bf V}_{canon}, {\bf W}_{canon}, 
\Phi_{{\bf V}_{canon}{\bf W}_{canon}})$ and 
 $({\bf V}_{canon}^{\circ}, {\bf W}_{canon}^{\circ}, 
\Phi_{{\bf V}_{canon}^{\circ}{\bf W}_{canon}^{\circ}})$
can be 
 constructed for $C-{\bf M}(E)E$ and $C-{\bf M}(E)E-\sum_{1\leq i\leq p}e_{k_i}$,
 respectively.

 The following lemma characterizes the 
 ${\bf V}_{canon}, {\bf V}_{canon}^{\circ}$
 for both $C-{\bf M}(E)E$ and $C-{\bf M}(E)E-\sum_{i\leq p} e_{k_i}$.

\medskip

\begin{lemm}\label{lemm; trivial}
Let $C$ be a $(1, 1)$ class on $M$ such that ${\cal E}_{C-c_1({\bf K}_M)}$
is ample on $M$.

Let $({\bf V}_{canon}, {\bf W}_{canon}, 
\Phi_{{\bf V}_{canon}{\bf W}_{canon}})$ and 
 $({\bf V}_{canon}^{\circ}, {\bf W}_{canon}^{\circ}, 
\Phi_{{\bf V}_{canon}^{\circ}{\bf W}_{canon}^{\circ}})$ denote the canonical
algebraic family Kuranishi models 
 of 
$\underline{C}=C-{\bf M}(E)E$ and 
$\underline{C}=C-{\bf M}(E)E-\sum_{i\leq p}e_{k_i}$. Then
${\bf V}_{canon}={\bf V}_{canon}^{\circ}$ and the bundle ${\bf V}_{canon}$
on 
 $T(M)\times M_n$ is the pullback from an algebraic vector bundle over $T(M)$
of rank$={C^2-C\cdot c_1({\bf K}_{M_{n+1}/M_n})\over 2}-q(M)+p_g 
+1$ by the projection map
 $T(M)\times M_n\mapsto T(M)$.
\end{lemm}

\medskip

\noindent Proof of the lemma: Assuming that ${\cal E}_{C-c_1({\bf K}_M)}$ is ample
 on $M$, then by Kodaira vanishing theorem the higher derived image sheaves
 of 
${\cal E}_{C-c_1({\bf K}_{\bf M})}\otimes {\cal K}_M\cong
 {\cal E}_C$ along
 the projection map $\pi_{T(M)}:M\times T(M)\mapsto T(M)$ vanish.

 Thus, ${\cal R}^0(\pi_{T(M)})_{\ast}\bigl({\cal E}_C\bigr)$
 is locally free of rank ${C^2-C\cdot c_1({\bf K}_M)\over 2}-q(M)+p_g(M)+1$. (the 
rank is determined by surface Riemann Roch formula) 

 On the other hand, in the definition 5.3 of [Liu3], both of 
the bundles ${\bf V}_{canon}$, ${\bf V}_{canon}^{\circ}$ are defined to be
the algebraic bundles associated to the locally free sheaf 
${\cal R}^0\pi_{\ast}\bigl({\cal E}_C\bigr)$, where the sheaf ${\cal E}_C$ 
is pulled back from $M\times T(M)$ to $M_{n+1}\times T(M)$ through the
 composition map

$$M_{n+1}\times T(M)\mapsto M_n\times M\times T(M)\mapsto M\times T(M).$$

 Therefore there is a commutative diagram

\[
\begin{array}{ccc}
  M_{n+1}\times T(M) & \longrightarrow & M\times T(M)  \\
   \Big\downarrow &  & \Big\downarrow   \\
  M_n\times T(M)  & \longrightarrow & T(M)  \\
\end{array}
\]

 Thus, ${\cal R}^0\pi_{\ast}\bigl({\cal E}_C\bigr)$ is pulled back from
 $T(M)$. $\Box$

\medskip

 The coherent sheaves ${\cal R}^0\pi_{\ast}\bigl({\cal E}_{C-{\bf M}(E)E-\sum_i
 e_{k_i}}\bigr)$ and
${\cal R}^0\pi_{\ast}\bigl({\cal E}_{C-{\bf M}(E)E}\bigr)$ determine 
 algebraic cones in the total space of
${\bf V}_{canon}^{\circ}={\bf V}_{canon}$ and their 
projectifications are the algebraic family moduli spaces
 ${\cal M}_{C-{\bf M}(E)E-\sum_i e_{k_i}}\mapsto M_n\times T(M)$, 
${\cal M}_{C-{\bf M}(E)E}\mapsto M_n\times T(M)$.

According to the general construction of algebraic family Kuranishi models, 
${\cal M}_{C-{\bf M}(E)E-\sum_i e_{k_i}}$, ${\cal M}_{C-{\bf M}(E)E}$
are sub-schemes of ${\bf P}_{M_n\times T(M)}({\bf V}_{canon})$
 which are the zero loci of the canonical sections of the 
obstruction bundles ${\bf H}\otimes \pi^{\ast}_{{\bf P}({\bf V}_{canon})}
{\bf W}_{canon}^{\circ}$ and ${\bf H}\otimes \pi^{\ast}_{{\bf P}({\bf V}_{canon})}
{\bf W}_{canon}$ .

Let  $\underline{C}=C-{\bf M}(E)E$.
Suppose that $\Phi_{{\bf V}_{canon}{\bf W}_{canon}}:
{\bf V}_{canon}\mapsto {\bf W}_{canon}$ denote the canonical algebraic 
Kuranishi map,
then the canonical section of ${\bf H}\otimes {\bf W}_{canon}$ on 
${\bf P}_{M_n\times T(M)}({\bf V}_{canon})$, $s_{canon}$, 
 is induced by the bundle map $\Phi_{{\bf V}_{canon}{\bf W}_{canon}}$
 as the following.
 Each ray in ${\bf V}_{canon}$ determines a unique image ray in
 $
{\bf W}_{canon}$ through $\Phi_{{\bf V}_{canon}{\bf W}_{canon}}$

 and induces a map

$${\bf H}^{\ast}\longrightarrow \pi^{\ast}_{{\bf P}({\bf V}_{canon})}
{\bf W}_{canon},$$

or equivalently a map

$$ {\bf HOM}({\bf C},\bigl({\bf H}^{\ast}\bigr)^{\ast}\otimes 
\pi^{\ast}_{{\bf P}({\bf V}_{canon})}{\bf W}_{canon})
\cong {\bf H}\otimes 
 \pi^{\ast}_{{\bf P}({\bf V}_{canon})}{\bf W}_{canon}.$$

This map can be viewed as the canonical section $s_{canon}$ of
 ${\bf H}\otimes \pi^{\ast}_{{\bf P}({\bf V}_{canon})}{\bf W}_{canon}$
 defining the algebraic family moduli space ${\cal M}_{\underline{C}}$ as 
 a sub-scheme in ${\bf P}({\bf V}_{canon})$.
 
  The discussion for $\underline{C}=C-{\bf M}(E)E-\sum_i e_{k_i}$ and
 $\Phi_{{\bf V}_{canon}^{\circ}{\bf W}_{canon}^{\circ}}:
{\bf V}_{canon}^{\circ}\mapsto {\bf W}_{canon}^{\circ}$ is parallel and
the corresponding canonical section is denoted by $s_{canon}^{\circ}$.

\bigskip

 When we restrict to the subspace $Y(\Gamma)\subset M_n$ of the
 universal space, the classes 
 $e_{k_i}$, $e_{k_i}\cdot (C-{\bf M}(E)E)<0$, 
$1\leq i\leq p$ become effective exceptional curve classes.
 We would like to  compare the 
 canonical algebraic obstruction bundles ${\bf W}_{canon}$ and 
${\bf W}_{canon}^{\circ}$.

 The main conclusion of the section is the following proposition,

\medskip

\begin{prop}\label{prop; com}
  Let $({\cal V}_{canon}^{\circ}, {\cal W}_{canon}^{\circ}, 
\Phi_{{\cal V}_{canon}^{\circ}{\cal W}_{canon}^{\circ}})$ and
 $({\cal V}_{canon}, {\cal W}_{canon}, 
\Phi_{{\cal V}_{canon}{\cal W}_{canon}})$ denote the sheaf theoretic
 version of algebraic 
canonical family Kuranishi models of $C-{\bf M}(E)E-\sum_{1\leq i\leq p}
e_{k_i}$ and $C-{\bf M}(E)E$, respectively.

  While restricting to the smooth subspace 
 $Y(\Gamma)\times T(M)\subset M_n\times T(M)$, there is a commutative diagram of
 sheaf morphisms between the two algebraic family Kuranishi models.

\[
\hskip -.9in
\begin{array}{ccccc}
 0 & & {\cal R}^0\pi_{\ast}\bigl({\cal O}_{\sum \Xi_{k_i}}\otimes 
{\cal E}_{C-{\bf M}E(E)}\bigr) 
& = & {\cal R}^0\pi_{\ast}\bigl({\cal O}_{\sum \Xi_{k_i}}\otimes 
{\cal E}_{C-{\bf M}E(E)}\bigr)  \\
\Big\downarrow & & \Big\downarrow & & \Big\downarrow \\
{\cal V}_{canon}^{\circ}|_{Y(\Gamma)\times T(M)}  & 
\stackrel{\Phi_{{\cal V}_{canon}^{\circ}{\cal 
W}_{canon}^{\circ}}}{\longrightarrow}& 
{\cal W}_{canon}^{\circ}|_{Y(\Gamma)\times T(M)}& \longrightarrow & 
 {\cal R}^1\pi_{\ast}\bigl({\cal E}_{C-{\bf M}(E)E-\sum_ie_{k_i}}\bigr)\\
\Big\downarrow  & &\Big\downarrow & & \Big\downarrow\\
{\cal V}_{canon}|_{Y(\Gamma)\times T(M)}  & \stackrel{\Phi_{
{\cal V}_{canon}{\cal W}_{canon}}}{\longrightarrow}& 
{\cal W}_{canon}|_{Y(\Gamma)\times T(M)} & 
\longrightarrow & 
{\cal R}^1\pi_{\ast}\bigl({\cal E}_{C-{\bf M}(E)E}\bigr)\\
 \Big\downarrow & & \Big\downarrow & & \Big\downarrow\\
 0 & & {\cal R}^1\pi_{\ast}\bigl({\cal O}_{\sum \Xi_{k_i}}\otimes 
{\cal E}_{C-{\bf M}E(E)}\bigr)  & = & 
{\cal R}^1\pi_{\ast}\bigl({\cal O}_{\sum \Xi_{k_i}}\otimes 
{\cal E}_{C-{\bf M}E(E)}\bigr)\\
\end{array}
\]
\end{prop}
 
\medskip

\noindent Proof of the proposition: The isomorphism of ${\cal V}_{canon}^{\circ}$
 and ${\cal V}_{canon}$ has been addressed in lemma \ref{lemm; trivial}.
The rows are portions of the four
 term exact sequences characterizing the algebraic family Kuranishi model maps.
The last (third) column is a portion of the long exact sequence 
 induced from the short exact sequence
 $$\hskip -.8in 0\mapsto {\cal E}_{C-{\bf M}(E)E-\sum e_{k_i}}\mapsto
 {\cal E}_{C-{\bf M}(E)E}\mapsto {\cal O}_{\sum_i \Xi_{k_i}}\otimes
 {\cal E}_{C-{\bf M}(E)E}\mapsto 0,$$

 as was mentioned earlier.

  The second column in the commutative diagram 
is a four term exact sequence on 
 ${\cal W}_{canon}^{\circ}$ and ${\cal W}_{canon}$.
  By [Liu3], ${\cal W}_{canon}^{\circ}$ and
 ${\cal W}_{canon}$ are by definition 
 ${\cal R}^0\pi_{\ast}\bigl({\cal O}_{\sum_{i\leq n}m_iE_i+\sum_{i\leq p} 
\Xi_{k_i}}\otimes
 {\cal E}_C\bigr)$ and ${\cal R}^0\pi_{\ast}\bigl(
{\cal O}_{\sum_{i\leq n} m_iE_i}\otimes
 {\cal E}_C\bigr)$, respectively.

 We start from the following exact sequence

$$\hskip -.8in 
0\mapsto {\cal O}_{\sum_{i\leq p}\Xi_{k_i}}(-\sum_{i\leq n}m_i E_i)\mapsto 
{\cal O}_{\sum_{i\leq n} m_i E_i+\sum_{i\leq p}\Xi_{k_i}}
\mapsto {\cal O}_{\sum_{i\leq n} m_i E_i}\mapsto
 0$$

 of divisors on $M_{n+1}\times_{M_n}Y(\Gamma)
\times T(M)$ and tensor it by ${\cal E}_C$.
 By using ${\cal E}_C\otimes {\cal O}(-\sum m_iE_i)
={\cal E}_{C-{\bf M}(E)E}$, we take its derived image long exact sequence.

Because 
${\cal R}^1\pi_{\ast}\bigl({\cal O}_{\sum_{i\leq n} m_iE_i+\sum_{i\leq p}
 \Xi_{k_i}}\otimes
 {\cal E}_C\bigr)=0$, the long exact sequence is truncated into a four-term
 exact sequence. The commutativity of the diagram follows from the
 naturality of all the sheaf morphisms involved and the compatibilities of
the connecting homomorphisms.
$\Box$

\bigskip

\begin{prop}\label{prop; coherent}
Let ${\bf H}$ denote the (restriction of the)
 hyperplane bundle of ${\bf P}({\bf V}_{canon})$.
 Let $s_{canon}^{\circ}\in \Gamma({\bf P}({\bf V}_{canon}^{\circ}),
 {\bf H}\otimes \pi^{\ast}_{{\bf P}({\bf V}_{canon}^{\circ})}
 {\bf W}_{canon}^{\circ})$ and 
$s_{canon}\in \Gamma({\bf P}({\bf V}_{canon}),
 {\bf H}\otimes \pi^{\ast}_{{\bf P}({\bf V}_{canon})}
 {\bf W}_{canon}$ denote the 
 the canonical sections defining ${\cal M}_{C-{\bf M}(E)E-\sum e_{k_i}}$,
 ${\cal M}_{C-{\bf M}(E)E}$. 

 Then the bundle map 

$${\bf H}\otimes \pi^{\ast}_{{\bf P}({\bf V}_{canon}^{\circ})}
{\bf W}_{canon}^{\circ}|_{Y(\Gamma)\times T(M)} 
\mapsto {\bf H}\otimes \pi^{\ast}_{{\bf P}(
{\bf V}_{canon})}{\bf W}_{canon}|_{Y(\Gamma)\times T(M)} $$
 induced by the bundle map $${\bf W}_{canon}^{\circ}|_{Y(\Gamma)\times T(M)} 
\mapsto {\bf W}_{canon}|_{Y(\Gamma)\times T(M)} $$
(introduced in the proof of prop. \ref{prop; com}) maps
 $s_{canon}^{\circ}$ to $s_{canon}$.
\end{prop}

\medskip

\noindent Proof of the proposition:
The commutativity of the diagram in prop. \ref{prop; com} implies the
 following commuative square on the corresponding algebraic vector bundles,

\[
\begin{array}{ccc}
 {\bf V}_{canon}^{\circ}|_{Y(\Gamma)\times T(M)} 
& \stackrel{\Phi_{{\bf V}_{canon}^{\circ}
{\bf W}_{canon}^{\circ}}}{\longrightarrow} 
& {\bf W}_{canon}^{\circ}|_{Y(\Gamma)\times T(M)} \\
\Arrowvert & & \Big\downarrow\\
 {\bf V}_{canon}|_{Y(\Gamma)\times T(M)} 
 & \stackrel{\Phi_{{\bf V}_{canon}{\bf W}_{canon}}}{
\longrightarrow} & {\bf W}_{canon}|_{Y(\Gamma)\times T(M)} \\
\end{array}
\]

 Then it implies the following commuative squares of bundle maps,

\[
\begin{array}{ccc}
 {\bf H}^{\ast} & \longrightarrow & \pi^{\ast}_{{\bf P}({\bf V}_{canon}^{\circ})}
{\bf W}_{canon}^{\circ}|_{Y(\Gamma)\times T(M)}  \\
  \Arrowvert & & \Big\downarrow \\
 {\bf H}^{\ast} & \longrightarrow & \pi^{\ast}_{{\bf P}(
{\bf V}_{canon})}{\bf W}_{canon}|_{Y(\Gamma)\times T(M)} \\
\end{array}
\]

  From this commutative square one derives the conclusion of the
 proposition immediately. $\Box$

\medskip

 If ${\cal R}^0\pi_{\ast}\bigl({\cal O}_{\sum_i \Xi_{k_i}}
\otimes {\cal E}_{C-{\bf M}(E)E}\bigr)$ is the zero sheaf over
 $Y(\Gamma)\times T(M)$, then 
 ${\cal R}^1\pi_{\ast}\bigl({\cal O}_{\sum_i \Xi_{k_i}}\bigr)$ is
 locally free and ${\cal W}_{canon}^{\circ}|_{Y(\Gamma)\times T(M)}
\mapsto {\cal W}_{canon}|_{Y(\Gamma)\times T(M)}$
 will be injective. 

 In this case the zero loci of $s_{canon}^{\circ}|_{Y(\Gamma)\times T(M)}$, 
${\cal M}_{C-{\bf M}(E)E-\sum_i e_{k_i}}\times_{M_n}Y(\Gamma)$, and of 
 $s_{canon}|_{Y(\Gamma)\times T(M)}$, 
${\cal M}_{C-{\bf M}(E)E}\times_{M_n}Y(\Gamma)$, coincide as subschemes of
 ${\bf P}({\bf V}_{canon})$.

\medskip

 In general the kernel sheaf 
${\cal R}^0\pi_{\ast}\bigl({\cal O}_{\sum_i \Xi_{k_i}}
\otimes {\cal E}_{C-{\bf M}(E)E}\bigr)$ may be non-zero. 
To analyze the difference of
 ${\cal M}_{C-{\bf M}(E)E-\sum_i e_{k_i}}\times_{M_n}Y(\Gamma)$
 and ${\cal M}_{C-{\bf M}(E)E}\times_{M_n}Y(\Gamma)$, one factorizes the
map 

$$\hskip -1in 
{\cal R}^0\pi_{\ast}\bigl({\cal O}_{\sum_{i\leq n}m_iE_i+\sum_{j\leq p}
 \Xi_{k_j}}\otimes {\cal E}_C\bigr)=
{\cal W}_{canon}^{\circ}|_{Y(\Gamma)\times T(M)}\mapsto 
{\cal W}_{canon}|_{Y(\Gamma)\times T(M)}={\cal R}^0\pi_{\ast}\bigl(
{\cal O}_{\sum_{i\leq n}m_iE_i}\otimes {\cal E}_C\bigr)$$ 

into a sequence of sheaf morphisms,

$${\cal R}^0\pi_{\ast}\bigl({\cal O}_{\sum m_iE_i+\sum_{1\leq i\leq p-r}
 \Xi_{k_i}}\otimes {\cal E}_C\bigr)
\mapsto {\cal R}^0\pi_{\ast}\bigl({\cal O}_{\sum m_iE_i+\sum_{1\leq i\leq p-r-1}
 \Xi_{k_i}}\otimes {\cal E}_C\bigr), 
 0\leq r\leq p-1,$$

 imitating the switching process $C-{\bf M}(E)E\mapsto C-{\bf M}(E)E-e_{k_1}
\mapsto \cdots \mapsto C-{\bf M}(E)E-\sum_{i\leq p}e_{k_i}$ on the cohomology
 classes.

\medskip

 The kernel of each of the above morphisms is isomorphic to 
 ${\cal R}^0\pi_{\ast}\bigl({\cal O}_{\Xi_{k_{p-r}}}\otimes
 {\cal E}_{C-{\bf M}(E)E-\sum_{1\leq i\leq p-r-1}e_{k_i}}\bigr)$, 
 $0\leq r\leq p-1$.

\medskip

\begin{lemm}\label{lemm; ker}
For $0\leq r\leq p-1$, the sheaf ${\cal R}^0\pi_{\ast}\bigl(
{\cal O}_{\sum_{i=r+1}^p \Xi_{k_i}}\otimes 
{\cal E}_{C-{\bf M}(E)E-\sum_{1\leq i\leq r}e_{k_i}}\bigr)$

 fits into the following exact sequence,

$$\hskip -.9in 0\mapsto {\cal R}^0\pi_{\ast}\bigl(
{\cal O}_{\sum_{i=r+2}^p \Xi_{k_i}}\otimes 
{\cal E}_{C-{\bf M}(E)E-\sum_{1\leq i\leq r+1}e_{k_i}}\bigr)
\mapsto {\cal R}^0\pi_{\ast}\bigl(
{\cal O}_{\sum_{i=r+1}^p \Xi_{k_i}}\otimes 
{\cal E}_{C-{\bf M}(E)E-\sum_{1\leq i\leq r}e_{k_i}}\bigr)$$

$$\hskip -1in \mapsto {\cal R}^0\pi_{\ast}\bigl(
{\cal O}_{\Xi_{k_{r+1}}}\otimes 
{\cal E}_{C-{\bf M}(E)E-\sum_{1\leq i\leq r}e_{k_i}}\bigr)\mapsto \cdots.$$

 And for each $r\in {\bf Z}$, $0\leq r\leq p-1$, 
there exists an equality in the $K$ group of coherent sheaves on
 $Y(\Gamma)\times T(M)$,

$$\hskip -1.3in [\pi_{\ast}\bigl(
{\cal O}_{\sum_{i=r+1}^p \Xi_{k_i}}\otimes 
{\cal E}_{C-{\bf M}(E)E-\sum_{1\leq i\leq r}e_{k_i}}\bigr)]
=[\pi_{\ast}\bigl({\cal O}_{\sum_{i=r+2}^p \Xi_{k_i}}\otimes 
{\cal E}_{C-{\bf M}(E)E-\sum_{1\leq i\leq r+1}e_{k_i}}\bigr)]+
[\pi_{\ast}\bigl(
{\cal O}_{\Xi_{k_{r+1}}}\otimes 
{\cal E}_{C-{\bf M}(E)E-\sum_{1\leq i\leq r}e_{k_i}}\bigr)].$$

\end{lemm}

\medskip

\noindent (A sketch of) the 
Proof: The proof of the lemma follows from taking the
 derived long exact sequence of the short exact sequence, 

$$\hskip -1.2in 0\mapsto {\cal O}_{\sum_{i=r+2}^p \Xi_{k_i}}\otimes 
{\cal E}_{C-{\bf M}(E)E-\sum_{1\leq i\leq r+1}e_{k_i}}\mapsto 
{\cal O}_{\sum_{i=r+1}^p \Xi_{k_i}}\otimes 
{\cal E}_{C-{\bf M}(E)E-\sum_{1\leq i\leq r}e_{k_i}}\mapsto 
{\cal O}_{\Xi_{k_{r+1}}}\otimes 
{\cal E}_{C-{\bf M}(E)E-\sum_{1\leq i\leq r}e_{k_i}}\mapsto 0.$$

 We omit the details.  $\Box$

\medskip

The rational curve fibrations dual to
 $e_{k_i}$, $\Xi_{k_i}, 1\leq i\leq r$ form ${\bf P}^1$ 
fibrations on $Y(\Gamma)$. Over the open strata $Y_{\Gamma}\subset Y(\Gamma)$ 
all the fibers of these
 $\Xi_{k_i}$ are smooth and irreducible. When the point 
 specializes to be in $Y(\Gamma)-Y_{\Gamma}$, the fibers of 
some of the $\Xi_{k_i}$ may
 break up into more than one irreducible component and different components smooth
 normal crossing ${\bf P}^1$.

\medskip

 The following proposition constraints the support of the kernel sheaf
${\cal R}^0\pi_{\ast}\bigl({\cal O}_{\sum_{i\leq p} \Xi_{k_i}}
\otimes {\cal E}_{C-{\bf M}(E)E}\bigr)$ of 
${\cal W}_{canon}^{\circ}|_{Y(\Gamma)\times T(M)}\mapsto 
{\cal W}_{canon}|_{Y(\Gamma)\times T(M)}$.

\begin{prop}\label{prop; char}
Let $Z_{k_i}\subset Y(\Gamma)-Y_{\Gamma}\subset M_n$ be the closed subset
 consisting of the singular values of 
 $\Xi_{k_i}\mapsto Y(\Gamma)$ (equivalently, the set over which the fibers
  of $\Xi_{k_i}$ fail to be irreducible). 
 Then the support of the sheaf ${\cal R}^0\pi_{\ast}\bigl({\cal O}_{\sum_{i\leq p}
 \Xi_{k_i}} \otimes {\cal E}_{C-{\bf M}(E)E}\bigr)$ is contained
 in the subset $(\cup_{1\leq i\leq p}Z_{k_i})\times T(M)$.
\end{prop}

\noindent Proof of the proposition: By using lemma \ref{lemm; ker}
 and by using induction, the support of  
${\cal R}^0\pi_{\ast}\bigl({\cal O}_{\sum_{i\leq p}
 \Xi_{k_i}} \otimes {\cal E}_{C-{\bf M}(E)E}\bigr)$, is
contained in the union of the supports

$$\cup_{1\leq i\leq p}supp\bigl({\cal R}^0\pi_{\ast}\bigl(
{\cal O}_{\Xi_{k_i}}\otimes 
{\cal E}_{C-{\bf M}(E)E-\sum_{1\leq j\leq i-1}e_{k_j}}\bigr)\bigr).$$

 It suffices to show that the support of 
${\cal R}^0\pi_{\ast}\bigl(
{\cal O}_{\Xi_{k_i}}\otimes 
{\cal E}_{C-{\bf M}(E)E-\sum_{1\leq j\leq i-1}e_{k_j}}\bigr)\subset Z_{k_i}\times
 T(M)$ and we argue by contradiction.

\medskip

If not, there exists at least an $i$ with  $1\leq i\leq p$ 
and a $t\not\in Z_{k_i}\times T(M)$
 such that
after a base change to $t$ the sheaf cohomology 
$H^0(\Xi_{k_i}|_t, {\cal O}_{\Xi_{k_i}|_t}
\otimes {\cal E}_{C-{\bf M}(E)E-\sum_{1\leq j\leq i-1}e_{k_j}})\not=0$.

Because each pair of distinct $1\leq a, b\leq p$,
 $\Xi_{k_a}|_{Y_{\Gamma}}, \Xi_{k_b}|_{Y_{\Gamma}}$ co-exist
 as ${\bf P}^1$ fiber bundles over $Y_{\Gamma}$, then for all $t\in Y_{\Gamma}$,
 $$e_{k_a}\cdot e_{k_b}=PD(\Xi_{k_a}|_t)\cdot PD(\Xi_{k_b}|_t)
\geq 0, 
  a\not=b.$$

But this implies that for $\Xi_{k_i}|_t\cong {\bf P}^1$, 
$$0\leq deg_{{\bf P}^1}{\cal O}_{{\bf P}^1}\otimes 
{\cal E}_{C-{\bf M}(E)E-\sum_{1\leq j\leq i-1}e_{k_j}}
=e_{k_i}\cdot (C-{\bf M}(E)E-\sum_{1\leq j\leq i-1}e_{k_j})$$
$$=e_{k_i}\cdot (C-{\bf M}(E)E)-e_{k_i}\cdot (\sum_{j\leq i-1}e_{k_j})\leq
 e_{k_i}\cdot (C-{\bf M}(E)E)<0.$$

Contradiction! $\Box$

\medskip

 The proposition singles out the geometric obstruction for 
 $s_{canon}^{\circ}$ and $s_{canon}$ to define the same zero locus over
 $Y(\Gamma)\times T(M)$,
 i.e. ${\cal M}_{C-{\bf M}(E)E-\sum_{i\leq p}e_{k_i}}|_{Y(\Gamma)\times T(M)}
\equiv {\cal M}_{C-{\bf M}(E)E}|_{Y(\Gamma)\times T(M)}$,
 to the possibility that the the fibers of $\Xi_{k_i}, 1\leq i\leq p$ can
 break into more than one irreducible component. This is exactly when 
the so-called higher level
 admissible decompositions pop up within 
the family ${\cal X}=M_{n+1}\times_{M_n}Y(\Gamma)\mapsto Y(\Gamma)=B$. 

 When this happens, the restriction of the algebraic family Kuranishi model
 of $C-{\bf M}(E)E-\sum_{i\leq p}e_{k_i}$, 
$({\cal V}_{canon}^{\circ}, {\cal W}_{canon}^{\circ}, 
\Phi_{{\cal V}_{canon}^{\circ}{\cal W}_{canon}^{\circ}})$, to
 $Y(\Gamma')\times T(M)$ ($\Gamma'<\Gamma$) is not
 an accurate approximation of the original algebraic family Kuranishi model 
$({\cal V}_{canon}, {\cal W}_{canon}, 
\Phi_{{\cal V}_{canon}{\cal W}_{canon}})$ any more and we have to choose a
 better approximation over $Y_{\Gamma'}\times T(M)$, $\Gamma'<\Gamma$, 
through replacing each $e_{k_i}, 1\leq i\leq p$ by
 the effective type $I$ exceptional classes $e'_{k_j}, 
e'_{k_j}\cdot (C-{\bf M}(E)E)<0, 1\leq j\leq p'
$ irreducible on $Y_{\Gamma'}$. (Consult the definition of $>$ on
 page \pageref{>}).

\bigskip
Over $Y(\Gamma)$ the exceptional classes $e_{k_i}, 1\leq i\leq p$
 are all effective. Then we may adjoin 
any effective curve dual to the class $C-{\bf M}(E)E-\sum_{1\leq i
\leq p} e_{k_i}$ with the union of type $I$ exceptional curves dual to 
$\sum_{1\leq i\leq p}e_{k_i}$ to produce 
 a curve dual to $C-{\bf M}(E)E$. This induces an inclusion
 $${\cal M}_{C-{\bf M}(E)E-\sum_{1\leq i \leq p} e_{k_i}}|_{Y(\Gamma)\times
 T(M)}\subset
 {\cal M}_{C-{\bf M}(E)E}|_{Y(\Gamma)\times
 T(M)}.$$

The following proposition characterizes the difference of 
${\cal M}_{C-{\bf M}(E)E-\sum_{1\leq i \leq p} e_{k_i}}$ and
 ${\cal M}_{C-{\bf M}(E)E}$ in terms of
 the intersection of $s_{canon}^{\circ}$ with the algebraic cone associated
 with ${\cal R}^0\pi_{\ast}\bigl({\cal O}_{\sum\Xi_{k_i}}\otimes 
 {\cal E}_{C-{\bf M}(E)E}\bigr)$.

\medskip

\begin{prop}\label{prop; intersect}
Let ${\cal H}$ denote the invertible sheaf associated with the 
hyperplane line bundle on 
 the prjectification of ${\bf V}_{canon}^{\circ}={\bf V}_{canon}$, 
 ${\bf P}_{Y(\Gamma)\times T(M)}({\bf V}_{canon}^{\circ})$,
 Then the locally free sheaf ${\cal H}\otimes 
\pi^{\ast}_{{\bf P}({\bf V}_{canon}^{\circ})}
{\cal W}_{canon}^{\circ}$ determines a vector bundle
 cone over ${\bf P}_{Y(\Gamma)\times T(M)}({\bf V}_{canon}^{\circ})$, the total 
space of ${\bf H}\otimes \pi^{\ast}_{{\bf P}({\bf V}_{canon}^{\circ})}
{\bf W}_{canon}^{\circ}|_{Y(\Gamma)\times T(M)}$.

 The sheaf injection 

$$\hskip -.9in 
0\mapsto {\cal H}\otimes \pi^{\ast}_{{\bf P}({\bf V}_{canon}^{\circ})}
{\cal R}^0\pi_{\ast}\bigl({\cal O}_{\sum_{i\leq p}
 \Xi_{k_i}} \otimes {\cal E}_{C-{\bf M}(E)E}
\bigr)\mapsto {\cal H}\otimes \pi^{\ast}_{{\bf P}({\bf V}_{canon}^{\circ})}
{\cal W}_{canon}^{\circ}|_{Y(\Gamma)\times T(M)}$$

 determines an algebraic subcone ${\bf C}_{\rho}$ of 
${\bf H}\otimes \pi^{\ast}_{{\bf P}({\bf V}_{canon}^{\circ})}
{\bf W}_{canon}^{\circ}|_{Y(\Gamma)\times T(M)}$ over 
${\bf P}_{Y(\Gamma)\times T(M)}({\bf V}_{canon}^{\circ})$.

 The section $s_{canon}^{\circ}$ determines a smooth sub-scheme 
, denoted by the bold character ${\bf s}_{canon}^{\circ}$, in 
 the total space of  ${\bf H}\otimes 
\pi^{\ast}_{{\bf P}_{Y(\Gamma)\times T(M)}({\bf V}_{canon}^{\circ})}
{\bf W}_{canon}^{\circ}|_{Y(\Gamma)\times T(M)}$ isomorphic to the base space 
 ${\bf P}_{Y(\Gamma)\times T(M)}({\bf V}_{canon}^{\circ})$ through the 
projection map $\pi_{{\bf W}_{canon}^{\circ}}:
{\bf H}\otimes \pi^{\ast}_{{\bf P}({\bf V}_{canon}^{\circ})}
{\bf W}_{canon}^{\circ}|_{Y(\Gamma)\times T(M)}\mapsto 
{\bf P}_{Y(\Gamma)\times T(M)}({\bf V}_{canon}^{\circ})$.

 The restriction of the algebraic family moduli space of 
 $C-{\bf M}(E)E$ to $Y(\Gamma)\times T(M)$, 
${\cal M}_{C-{\bf M}(E)E}\times_{M_n}Y(\Gamma)=s^{-1}_{canon}(0)
\subset {\bf P}_{Y(\Gamma)\times T(M)}({\bf V}_{canon})$,
 can be identified with the image 
$\subset {\bf P}_{Y(\Gamma)\times T(M)}
({\bf V}_{canon}^{\circ})\equiv {\bf P}_{Y(\Gamma)\times T(M)}({\bf V}_{canon})$ 
of the scheme theoretical intersection
 ${\bf s}_{canon}^{\circ}\cap {\bf C}_{\rho}$ under the
 projection 
 $\pi_{{\bf H}\otimes 
\pi_{{\bf P}_{Y(\Gamma)\times T(M)}
({\bf V}_{canon}^{\circ})}^{\ast}{\bf W}_{canon}^{\circ}}$.
\end{prop}

\medskip

\noindent Proof of the proposition: To simplify our notations, we
will drop the subscript $Y(\Gamma)\times T(M)$ and denote 
${\bf P}_{Y(\Gamma)\times T(M)}
({\bf V}_{canon}^{\circ})\equiv {\bf P}_{Y(\Gamma)\times T(M)}
({\bf V}_{canon})$ by ${\bf P}({\bf V}_{canon}^{\circ})$.

Let us recall some basic facts about 
 the construction of algebraic cones.
By [F] B.5.5. page 434, the algebraic vector bundle cone determined
 by a locally free sheaf ${\cal D}$ over a scheme $X$
 is $Spec(Sym({\cal D}^{\ast}))$, where 
${\cal D}^{\ast}={\cal HOM}_{{\cal O}_X}({\cal D}, 
{\cal O}_X)$.

Moreover, for a coherent sheaf ${\cal R}$ over a scheme
 $X$, one may construct an algebraic cone
 over $X$ by the recipe ${\bf C}({\cal R})=Spec({\bf S}^{\bullet}({\cal R}))$,
 where ${\bf S}^{\bullet}=Sym({\cal R})$ is the sheaf of 
 graded ${\cal O}_X$ algebra generated by
 ${\bf S}^1={\cal R}$ (of grade one).

  In the current context, we 
take $X={\bf P}_{Y(\Gamma)\times T(M)}({\bf V}_{canon}^{\circ})$.
 Take ${\cal P}$ to be the image of the sheaf morphism 
 $g:{\cal H}\otimes \pi^{\ast}_{{\bf P}({\bf V}_{canon}^{\circ})}
{\cal W}_{canon}^{\circ}|_{Y(\Gamma)\times Y(\Gamma)}\mapsto 
{\cal H}\otimes \pi^{\ast}_{{\bf P}({\bf V}_{canon}^{\circ})}
{\cal W}_{canon}|_{Y(\Gamma)\times Y(\Gamma)}$. 
Then there is a short exact sequence 

$$\hskip -.4in 0\mapsto 
{\cal H}\otimes \pi^{\ast}_{{\bf P}({\bf V}_{canon}^{\circ})}
{\cal R}^0\pi_{\ast}\bigl({\cal O}_{\sum_{i\leq p}
 \Xi_{k_i}} \otimes {\cal E}_{C-{\bf M}(E)E}\bigr)
\mapsto {\cal H}\otimes \pi^{\ast}_{{\bf P}({\bf V}_{canon}^{\circ})}
{\cal W}_{canon}^{\circ}|_{Y(\Gamma)\times T(M)}\mapsto {\cal P}\mapsto 0.$$

 Because ${\cal HOM}_{{\cal O}_X}(\cdot, {\cal O}_X)$ is a contra-variant 
left exact functor on the category of sheaves of ${\cal O}_X$ modules, there
 is a short exact sequence over $X={\bf P}({\bf V}_{canon}^{\circ})$,

$$\hskip -.9in 0\mapsto {\cal P}^{\ast}\mapsto 
\bigl({\cal H}\otimes \pi^{\ast}_{{\bf P}({\bf V}_{canon}^{\circ})}
{\cal W}_{canon}^{\circ}|_{Y(\Gamma)\times T(M)}\bigr)^{\ast}\mapsto 
 {\cal R}\mapsto 0,$$

 where ${\cal R}$ is defined to be the cokernel of the injective 
morphism 
 ${\cal P}^{\ast}\mapsto 
\bigl({\cal H}\otimes \pi^{\ast}_{{\bf P}({\bf V}_{canon}^{\circ})}
{\cal W}_{canon}^{\circ}|_{Y(\Gamma)\times T(M)}\bigr)^{\ast}$.

 We define the cone ${\bf C}_{\rho}\equiv 
{\bf C}({\cal R})$ by the above recipe in [F].
 We may take ${\cal D}$ to be the locally free $\bigl({\cal H}\otimes
 \pi^{\ast}_{{\bf P}({\bf V}_{canon}^{\circ})}
{\cal W}_{canon}^{\circ}|_{Y(\Gamma)\times T(M)}\bigr)^{\ast}$ and
 then ${\bf C}({\cal D})$ is the associated algebraic vector bundle cone.

 Then the above sheaf surjection ${\cal D}\mapsto {\cal R}\mapsto 0$ 
 induces a surjection on the corresponding sheaves of ${\cal O}_{X}$ 
algebras and thus an inclusion of cones
${\bf C}({\cal R})\subset {\bf C}({\cal D})$.

 The next lemma relates the algebraic sub-cone 
${\bf C}_{\rho}$ with the coherent sheaf
 ${\cal H}\otimes \pi^{\ast}_{{\bf P}({\bf V}_{canon}^{\circ})}
{\cal R}^0\pi_{\ast}\bigl({\cal O}_{\sum_{i\leq p}
 \Xi_{k_i}} \otimes {\cal E}_{C-{\bf M}(E)E}\bigr)$.

\medskip

\begin{lemm} \label{lemm; cone}
 The dual of the coherent sheaf ${\cal R}$ satisfies 
${\cal R}^{\ast}={\cal H}\otimes \pi^{\ast}_{{\bf P}({\bf V}_{canon}^{\circ})}
{\cal R}^0\pi_{\ast}\bigl({\cal O}_{\sum_{i\leq p}
 \Xi_{k_i}} \otimes {\cal E}_{C-{\bf M}(E)E}\bigr)$.

 For all $t\in {\bf P}({\bf V}_{canon}^{\circ})$, the fiber of the cone, 
${\bf C}_{\rho}|_t={\bf C}_{\rho}\times_{{\bf P}({\bf V}_{canon}^{\circ})}\{t\}$
 at $t$ is canonically isomorphic to 
${\cal H}\otimes \pi^{\ast}_{{\bf P}({\bf V}_{canon}^{\circ})}
{\cal R}^0\pi_{\ast}\bigl({\cal O}_{\sum_{i\leq p}
 \Xi_{k_i}} \otimes {\cal E}_{C-{\bf M}(E)E}\bigr)\otimes k(t)$.
\end{lemm}

\medskip

\noindent Proof:  By dualizing the exact sequence we get

$$0\mapsto {\cal R}^{\ast}\mapsto 
\bigl(\bigl({\cal H}\otimes \pi^{\ast}_{{\bf P}({\bf V}_{canon}^{\circ})}
{\cal W}_{canon}^{\circ}|_{Y(\Gamma)\times T(M)}\bigr)^{\ast}\bigr)^{\ast}
 \mapsto \bigl({\cal P}^{\ast}\bigr)^{\ast}.$$

 Because ${\cal H}\otimes \pi^{\ast}_{{\bf P}({\bf V}_{canon}^{\circ})}
{\cal W}_{canon}^{\circ}|_{Y(\Gamma)\times T(M)}$ is locally free,
$\bigl(\bigl({\cal H}\otimes \pi^{\ast}_{{\bf P}({\bf V}_{canon}^{\circ})}
{\cal W}_{canon}^{\circ}|_{Y(\Gamma)\times T(M)}\bigr)^{\ast}\bigr)^{\ast}
\cong {\cal H}\otimes \pi^{\ast}_{{\bf P}({\bf V}_{canon}^{\circ})}
{\cal W}_{canon}^{\circ}|_{Y(\Gamma)\times T(M)}$ canonically.
 It suffices to show that the map
${\cal H}\otimes \pi^{\ast}_{{\bf P}({\bf V}_{canon}^{\circ})}
{\cal W}_{canon}^{\circ}|_{Y(\Gamma)\times T(M)}\mapsto \bigl({\cal P}^{\ast}
\bigr)^{\ast}$ factors through ${\cal P}\subset \bigl({\cal P}^{\ast}\bigr)^{\ast}$
 and therefore coincides with the original 
$$g:{\cal H}\otimes \pi^{\ast}_{{\bf P}({\bf V}_{canon}^{\circ})}
{\cal W}_{canon}^{\circ}|_{Y(\Gamma)\times T(M)}\mapsto {\cal P}\mapsto 0.$$

 Starting from 

$$g:{\cal H}\otimes \pi^{\ast}_{{\bf P}({\bf V}_{canon}^{\circ})}
{\cal W}_{canon}^{\circ}|_{Y(\Gamma)\times T(M)}\mapsto {\cal P}\subset
{\cal H}\otimes \pi^{\ast}_{{\bf P}({\bf V}_{canon}^{\circ})}
{\cal W}_{canon}|_{Y(\Gamma)\times T(M)}$$

 and take the double dual of $g$. On the one hand, it is easy to see
 that $(g^{\ast})^{\ast}$ factors through $\bigl({\cal P}^{\ast}\bigr)^{\ast}$.
On the other hand, both 
${\cal H}\otimes \pi^{\ast}_{{\bf P}({\bf V}_{canon}^{\circ})}
{\cal W}_{canon}^{\circ}|_{Y(\Gamma)\times Y(\Gamma)}$ and 
${\cal H}\otimes \pi^{\ast}_{{\bf P}({\bf V}_{canon}^{\circ})}
{\cal W}_{canon}|_{Y(\Gamma)\times Y(\Gamma)}$ are locally free,
thus $(g^{\ast})^{\ast}\equiv g$ canonically. As $g$ factors through
 ${\cal P}$, so is $(g^{\ast})^{\ast}$.

Once this is established, ${\cal R}^{\ast}$ must be isomorphic to
${\cal H}\otimes \pi^{\ast}_{{\bf P}({\bf V}_{canon}^{\circ})}
{\cal R}^0\pi_{\ast}\bigl({\cal O}_{\sum_{i\leq p}
 \Xi_{k_i}} \otimes {\cal E}_{C-{\bf M}(E)E}\bigr)$ as both 
 are the kernel of the morphism $g$.

 By the construction of ${\bf C}_{\rho}$ it is apparent that
 $${\bf C}_{\rho}|_t=Spec(Sym({\cal R})\otimes_{{\cal O}_X}
 k(t))=Spec(Sym({\cal R}\otimes_{{\cal O}_X}
 k(t)))=HOM_{k(t)}({\cal R}\otimes_{{\cal O}_X} k(t), k(t))$$
$$={\cal HOM}_{{\cal O}_X}({\cal R}, {\cal O}_X)\otimes_{{\cal O}_X} k(t)
={\cal R}^{\ast}\otimes_{{\cal O}_X} k(t),$$

 and the first part of the lemma already concludes that
${\cal R}^{\ast}={\cal H}\otimes \pi^{\ast}_{{\bf P}({\bf V}_{canon}^{\circ})}
{\cal R}^0\pi_{\ast}\bigl({\cal O}_{\sum_{i\leq p}
 \Xi_{k_i}} \otimes {\cal E}_{C-{\bf M}(E)E}\bigr)$. 

 Therefore ${\bf C}_{\rho}|_t=
{\cal H}\otimes \pi^{\ast}_{{\bf P}({\bf V}_{canon}^{\circ})}
{\cal R}^0\pi_{\ast}\bigl({\cal O}_{\sum_{i\leq p}
 \Xi_{k_i}} \otimes {\cal E}_{C-{\bf M}(E)E}\bigr)\otimes k(t)$.
$\Box$

\medskip

 Let ${\bf s}_{canon}$ and ${\bf O}$ denote the sub-schemes of the total space
of ${\bf H}\otimes \pi^{\ast}_{{\bf P}({\bf V}_{canon}^{\circ})}
{\bf W}_{canon}|_{Y(\Gamma)\times T(M)}$ 
defined by the sections $s_{canon}$ and the
 zero section. Likewise set ${\bf O}^{\circ}$ to be the zero section
 sub-scheme of the total space of 
${\bf H}\otimes \pi^{\ast}_{{\bf P}({\bf V}_{canon}^{\circ})}
{\bf W}_{canon}^{\circ}|_{Y(\Gamma)\times T(M)}$. 
 We notice that the algebraic family moduli space 
${\cal M}_{C-{\bf M}(E)E}\times_{M_n}Y(\Gamma)$, defined to be
 $s^{-1}_{canon}(0)$, can be viewed as the 
intersection ${\bf s}_{canon}\cap{\bf O}$.

 Define $g$ to be the algebraic map on the total spaces of the vector bundles
 induced by ${\bf H}\otimes \pi^{\ast}_{{\bf P}({\bf V}_{canon}^{\circ})}
{\bf W}_{canon}^{\circ}|_{Y(\Gamma)\times T(M)} 
\mapsto {\bf H}\otimes \pi^{\ast}_{{\bf P}(
{\bf V}_{canon})}{\bf W}_{canon}|_{Y(\Gamma)\times T(M)}$.
 The proposition \ref{prop; coherent} implies that 
$g({\bf s}_{canon}^{\circ})={\bf s}_{canon}$.

 Moreover, the vector bundles projection maps induce the following
 commutative diagram of isomorphisms.

\[
\begin{array}{ccc}
 {\bf s}_{canon}^{\circ} & \stackrel{g}{\longrightarrow} & {\bf s}_{canon}\\
 \Big\downarrow & & \Big\downarrow\\
 {\bf O}_{canon}^{\circ}  & \stackrel{g}{=} & {\bf O}_{canon}\\
\end{array}
\]

 On the other hand, lemma \ref{lemm; cone} implies that ${\bf C}_{\rho}=
g^{-1}({\bf O})$. Thus, 

$$g({\bf s}_{canon}^{\circ}\cap {\bf C}_{\rho})
=g({\bf s}_{canon}^{\circ})\cap g({\bf C}_{\rho})
 ={\bf s}_{canon}\cap {\bf O}.$$

 On the scheme theoretic level, the ideal sheaf ${\cal I}'$ defines 
 ${\bf s}_{canon}^{\circ}\cap {\bf C}_{\rho}\subset {\bf s}_{canon}^{\circ}$
 is the inverse image ideal sheaf $g^{-1}{\cal I}\cdot 
{\cal O}_{{\bf s}_{canon}^{\circ}}$ 
(see [Ha] page 163) of the ideal sheaf ${\cal I}$
 defining ${\bf s}_{canon}\cap {\bf O}\subset {\bf s}_{canon}$ under the
 map $g:{\bf s}_{canon}^{\circ}\mapsto {\bf s}_{canon}$. 

  Since the commutative diagram of isomorphisms imply that $g$ induces
 an isomorphism between ${\bf s}_{canon}^{\circ}$ and ${\bf s}_{canon}$,
 and therefore 
 $g^{-1}{\cal O}_{{\bf s}_{canon}}={\cal O}_{{\bf s}_{canon}^{\circ}}$.
 Thus, 
$g^{-1}{\cal I}\subset 
g^{-1}{\cal O}_{{\bf s}_{canon}}={\cal O}_{{\bf s}_{canon}^{\circ}}$ and
 ${\cal I}'=g^{-1}{\cal I}\cdot 
{\cal O}_{{\bf s}_{canon}^{\circ}}=g^{-1}{\cal I}\cong {\cal I}$.

 Thus, the scheme ${\bf s}_{canon}\cap {\bf O}$ and ${\bf s}_{canon}^{\circ}\cap
 {\bf C}_{\rho}$ are isomorphic.

 This ends the proof of proposition \ref{prop; intersect}. $\Box$

\bigskip

 In general a cone can be decomposed into its irreducible components.
Therefore, we may write ${\bf C}_{\rho}$ as 
 $\sum_{\rho_i}{\bf C}_{\rho_i}$, where each ${\bf C}_{\rho_i}$
is an irreducible sub-cone of ${\bf C}_{\rho}$. By lemma \ref{lemm; cone},
 the cone ${\bf C}_{\rho}$ always contains an irreducible sub-cone
 ${\bf O}^{\circ}$, the zero section cone of the vector bundle. Then there
 is a distinquished irreducible
 component among ${\bf C}_{\rho_i}$, called ${\bf C}_{\rho_0}=
{\bf O}^{\circ}$.

 We have the following immediate corollary of proposition \ref{prop; intersect}.

\medskip

\begin{cor}\label{cor; decom}
Over the closure of the admissible strata, $Y(\Gamma)$, the algebraic 
 family moduli space of curves in $C-{\bf M}(E)E$, 
 ${\cal M}_{C-{\bf M}(E)E}\times_{M_n}Y(\Gamma)$ admits a decomposition
 into $${\cal M}_{C-{\bf M}(E)E}\times_{M_n}Y(\Gamma)=
{\cal M}_{C-{\bf M}(E)E-\sum e_{k_i}}\times_{M_n}Y(\Gamma)\cup
\bigcup_{i\not=0} 
\pi_{{\bf H}\otimes {\bf W}_{canon}}({\bf s}_{canon}^{\circ}\cap 
{\bf C}_{\rho_i}).$$
\end{cor}

 Each of the term 
$\pi_{{\bf H}\otimes {\bf W}_{canon}}({\bf s}_{canon}^{\circ}\cap 
{\bf C}_{\rho_i})$ is a closed sub-scheme of $X={\bf P}({\bf V}_{canon}^{\circ})$.

\bigskip

\section{\bf The Localized Top Chern Class and The Contribution to
the Family Invariants}\label{section; example}

\bigskip

 At the end of this paper, we would like to discuss some 
implication of the discussion presented above. Let us begin by
 reviewing the concept of localized top Chern class of a vector bundle (see
 [F] p244 for the details).

 Let $E\mapsto X$ be a rank $e$ vector bundle on a purely $m$ dimensional 
scheme $X$ and $s:X\mapsto E$ be a section with zero scheme $Z(s)$.
 Let $s_E$ denote the zero section. Define ${\bf Z}(s)=s_E^{!}([X])\in 
{\cal A}_{m-e}(Z(s))$. It is called the localized top Chern class of 
$E$ with respect to $s$. The most important property of ${\bf Z}(s)$ is
 $i_{\ast}{\bf Z}(s)=c_{e}(E)\cap [X]\in {\cal A}_{m-e}(X)$ under the inclusion
map $i:Z(s)\mapsto X$. Namely it is mapped to the (global) top Chern class of $E$
 under the inclusion map $i$.

  If one goes through the definition of $s_E^{!}$, one may write 
 $s_E^{!}([X])$ in an alternative way.

\begin{prop}\label{prop; local}
 Let ${\bf C}_{Z(s)}$ denote the normal cone of $Z(s)$ in $X$ and let
 $s(Z(s), X)=s({\bf C}_{Z(s)})$ denote the total Segre class of the normal cone. 
Then the localized top Chern class ${\bf Z}(s)$ is equal to
 $\{c(i^{\ast}E)\cap s(Z(s), X)\}_{m-e}$.

Moreover, for all $r\in {\bf N}$ the cycle class 
$\{c(i^{\ast}E)\cap s(Z(s), X)\}_{m-e-r}=0$.
\end{prop}

\medskip

\noindent Proof: This is the consequence of the commutative diagram on
 [F] page 244 and proposition 6.1.(a) [F] page 94.
The object $\{c(i^{\ast}E)\cap s(Z(s), X)\}_{m-e-r}$ vanishes in the
 cycle class group because in the proof of proposition 6.1.(a) [F]
 $c_{d+r}(\xi)=0$ for $d=rank_{\bf C}\xi$, $r\in {\bf N}$. 
 $\Box$

The constraint on the grading will be used extensively in the following
 discussion.

\medskip

 Insteading of taking $Z(s)$, we consider $Y\subset Z(s)$ to be
 a closed sub-scheme and the denote inclusion $Y\subset X$ by $i_Y$.
Then the expression
 $\{c(i_Y^{\ast}E)\cap s(Y, X)\}_{m-e}\in {\cal A}_{m-e}(Y)$ defines
a cycle class localized in $Y$.

\medskip

\begin{defin}\label{defin; contribution}
 Define the cycle class 
${\bf Z}_Y(s)=\{c(i_Y^{\ast}E)\cap s(Y, X)\}_{m-e}$ to be
 the localized top Chern class contribution of $Y\subset Z(s)$.
\end{defin}

\medskip

 The geometric meaning of this definition is clarified by the following
 proposition.

\medskip

\begin{prop}\label{prop; sum}
Suppose that the zero locus $Z(s)$ can be decomposed into $\coprod_q Y_q$
 such that $Y_{q_1}\cap Y_{q_2}=\emptyset$, whenever $q_1\not=q_2$ and let 
$j_{Y_q}:Y_q\mapsto Z(s)$ denote the inclusion map.

Then $$\sum_q{j_{Y_q}}_{\ast}{\bf Z}_{Y_q}(s)={\bf Z}(s).$$  
\end{prop}

\medskip

\noindent Proof: When $Z(s)=\coprod_q Y_q$, ${\bf C}_{Z(s)}=\coprod_q
 {\bf C}_{Y_q}$.  Thus ${\cal A}_{\ast}(Z(s))\ni s(Z(s), X)=
 \sum_q j_{Y_q\ast}s(Y_q, X)$.

 Then by projection formula (theorem 3.2.(c) on [F] page 50)
 $${\bf Z}(s)=\{c(i^{\ast}E)\cap s(Z(s), X)\}_{m-e}
=\sum_q\{c(i^{\ast}E)\cap \{j_{Y_q}\}_{\ast}s(Y_q, X)\}_{m-e}$$
$$=\sum_q{j_{Y_q}\ast}\bigl(c(j_{Y_q}^{\ast}i^{\ast}E)\cap s(Y_q, X)\bigr)
\}_{m-e}
=\sum_q\{{j_{Y_q}}_{\ast}\bigl(c(i_{Y_q}^{\ast}E)\cap s(Y_q, X)\bigr)
\}_{m-e}=\sum_q{j_{Y_q}}_{\ast}{\bf Z}_{Y_q}(s).$$ $\Box$

   We apply this set up to the family invariant of $C-{\bf M}(E)E$,
 ${\cal AFSW}_{M_{n+1}\times T(M)\mapsto M_n\times T(M)}(1, C-{\bf M}(E)E)$.

\medskip

\subsection{\bf An Identification Upon the Local Contribution}
\label{subsection; localcont}

\bigskip

 In our setting, we take $X={\bf P}_{M_n\times T(M)}({\bf V}_{canon})$ and
 $E={\bf H}\otimes \pi_{X}^{\ast}{\bf W}_{canon}$. Then the algebraic family
 moduli space ${\cal M}_{C-{\bf M}(E)E}=Z(s_{canon})$ and the
 family invariant ${\cal AFSW}_{M_{n+1}\times T(M)
\mapsto M_n\times T(M)}(1, C-{\bf M}(E)E)$
 is defined to be 

$${\bf Z}\cong 
{\cal A}_0(X)\ni c_{top}({\bf H}\otimes \pi_{X}^{\ast}{\bf W}_{canon})\cap 
 c_1({\bf H})^{dim_{\bf C}M_n+
{C^2-C\cdot c_1({\bf K}_M)\over 2}+p_g-\sum_i {m_i^2+m_i\over 2}}.$$

As in the earlier sections
 let $\Gamma\not=\gamma_n$ be an $n$-vertex admissible graph such that 
all the type $I$ exceptional classes $e_i, 1\leq i\leq n$ over $Y(\Gamma)$
satisfy the following condition, 

\noindent {\bf Special Condition}: \label{Scondition}

either 

\noindent (i). $(C-{\bf M}(E)E)\cdot e_i<0$, i.e. ${\bf M}(E)E\cdot e_i>0$. 

or 

\noindent (ii). the condition 
$e_i^2=-1$, i.e. $e_i$ is a $-1$ type $I$ exceptional class.

Let $k_i, 1\leq i\leq p$ be the
 subscripts in $\{1, 2, \cdots, n\}$ such that $(C-{\bf M}(E)E)\cdot e_{k_i}<0$.
 By permuting the indexes we may assume $k_i=i$ for $1\leq i\leq p$.
 From now on we adopt this simplified notation. 

Because $Y(\Gamma)\subset M_n$ is a closed inclusion, the restriction
 ${\cal M}_{C-{\bf M}(E)E}\times_{M_n}Y(\Gamma)$ is a closed sub-scheme of
 ${\cal M}_{C-{\bf M}(E)E}=Z(s_{canon})$. Then we may take 
 $Y={\cal M}_{C-{\bf M}(E)E}\times_{M_n}Y(\Gamma)$ and definition
 \ref{defin; contribution} determines a localized top Chern class contribution
 of $Y\subset {\cal M}_{C-{\bf M}(E)E}$.

 Then $$\hskip -1.2in {\bf Z}_Y(s)\cap c_1(i_Y^{\ast}{\bf H})^{{C^2-C\cdot 
c_1({\bf K}_M)\over 2}+p_g-\sum_i {m_i^2+m_i\over 2}}
=\{c(i_Y^{\ast}{\bf H}\otimes \pi_{X}^{\ast}{\bf W}_{canon})\cap
 s(Y, X)\}_{dim_{\bf C}M_n+{C^2-C\cdot 
c_1({\bf K}_M)\over 2}+p_g-\sum_i {m_i^2+m_i\over 2}}$$
$$\cap
 c_1(i_Y^{\ast}{\bf H})^{dim_{\bf C}M_n+{C^2-C\cdot 
c_1({\bf K}_M)\over 2}+p_g-\sum_{1\leq i\leq n} {m_i^2+m_i\over 2}}
\in 
{\cal A}_0({\cal M}_{C-{\bf M}(E)E}\times_{M_n}Y(\Gamma))\cong {\bf Z}$$

 gives a localized contribution of the algebraic
 family invariant over $Y(\Gamma)$.

\medskip

\noindent {\bf Question 1}: 
Can we express (enumerate) the localized contribution of
the (algebraic) family invariant of $C-{\bf M}(E)E$ 
over $Y(\Gamma)$ in terms of the family invariant of some other classes?

\medskip

  The complete answer of this question in terms of differential topology has been
 presented in [Liu1] by using the concept of modified family invariants,
 the complete solution by a purely algebraic approach is beyond the
scope of the present paper. Instead, we try to motivate the readers by
 providing a light-weighted version which answers the following questions.

\medskip

\noindent {\bf Question 2}: What is the explicit 
form of the typical family invariant that we express
 the localized contributions of the family
 invariant over $Y(\Gamma)$?

\medskip

\noindent {\bf Question 3}: Why does the procedure of enumerating 
the local contributions of the algebraic family invariants involve the so-called
 higher level admissible decomposition classes (defined in [Liu1]), 
i.e. the local contributions from the $Y(\Gamma')$, $\Gamma'<\Gamma$?
\label{question3}

\medskip

 The conceptual understanding of both {\bf Question 2} and {\bf Question 3} are
 essential to understand the solution of {\bf Question 1}. 

\medskip

 In corollary \ref{cor; decom} on page \pageref{cor; decom}, we have shown that 
${\cal M}_{C-{\bf M}(E)E}\times_{M_n} Y(\Gamma)$
 can be decomposed into ${\cal M}_{C-{\bf M}(E)E-\sum_i e_{k_i}}
\times_{M_n}Y(\Gamma)$ and 
 $\cup_{i\not=0}\pi_{{\bf H}\otimes \pi^{\ast}_{{\bf P}({\bf V}_{canon}^{\circ})}
{\bf W}_{canon}}({\bf C}_{\rho_i} \cap {\bf s}_{canon}^{\circ})$.

\bigskip

 As our goal is to illustrate the patterns and difficulties involved,  
 we make an additional assumption to simplify the discussion while
 the general situation without imposing this assumption will be treated
 elsewhere [Liu5] by using the residual intersection theory.

\medskip

\noindent {\bf Simplifying Assumption}: 
The space ${\cal C}_{C-{\bf M}(E)E-\sum_i e_{k_i}}$
 and $\cup_{i\not=0}\pi_{{\bf H}\otimes 
\pi^{\ast}_{{\bf P}({\bf V}_{canon}^{\circ})}
{\bf W}_{canon}}({\bf C}_{\rho_i} \cap {\bf s}_{canon}^{\circ})$
are disjoint, i.e. ${\cal C}_{C-{\bf M}(E)E-\sum_i e_{k_i}}\cap 
\cup_{i\not=0}\pi_{{\bf H}\otimes 
\pi^{\ast}_{{\bf P}({\bf V}_{canon}^{\circ})}
{\bf W}_{canon}}({\bf C}_{\rho_i} \cap {\bf s}_{canon}^{\circ})=\emptyset$.
\label{assum}
\medskip

 As a direct consequence of this assumption we may name
 $Y_1={\cal C}_{C-{\bf M}(E)E-\sum_i e_{k_i}}\times_{M_n}Y(\Gamma)$, 
$Y_2=\cup_{i\not=0}\pi_{{\bf H}\otimes 
\pi^{\ast}_{{\bf P}({\bf V}_{canon}^{\circ})}
{\bf W}_{canon}}({\bf C}_{\rho_i} \cap {\bf s}_{canon}^{\circ})$, and
 $Y={\cal M}_{C-{\bf M}(E)E}=Y_1\coprod Y_2$.

 Then the argument of proposition \ref{prop; sum} implies that (after
replacing $Z(s)$ 
 by $Z(s)\times_X Y$)

$${\bf Z}_Y(s_{canon})=
{j_{Y_1}}_{\ast}{\bf Z}_{Y_1}(s_{canon})+{j_{Y_2}}_{\ast}
{\bf Z}_{Y_2}(s_{canon}).$$

 This implies that the localized contribution of the algebraic family 
invariant of $C-{\bf M}(E)E$ over $Y(\Gamma)$ can be decomposed into
 two parts, one from \label{yone}
${j_{Y_1}}_{\ast}{\bf Z}_{Y_1}(s_{canon})\cap 
 c_1(i_Y^{\ast}{\bf H})^{dim_{\bf C}M_n+{C^2-C\cdot 
c_1({\bf K}_M)\over 2}+p_g-\sum_{i\leq n}{m_i^2+m_i\over 2}}\in {\bf Z}$ 
and another from 
 ${j_{Y_2}}_{\ast}{\bf Z}_{Y_2}(s_{canon})\cap c_1(i_Y^{\ast}{\bf H})^{
{C^2-C\cdot c_1({\bf K}_M)\over 2}+p_g-\sum_{i\leq n} {m_i^2+m_i\over 2}}
\in {\bf Z}$.

 To answer {\bf Question 1}, we have the following theorem:

\medskip

\begin{theo}\label{theo; question1}
 Under the {\bf Simplifying Assumption} above, the integer
 ${j_{Y_1}}_{\ast}{\bf Z}_{Y_1}(s_{canon})\cap 
c_1(i_Y^{\ast}{\bf H})^{dim_{\bf C}M_n+{C^2-C\cdot 
c_1({\bf K}_M)\over 2}+p_g-\sum_{i\leq n} {m_i^2+m_i\over 2}}$
 can be identified with the mixed family invariant

$${\cal AFSW}_{M_{n+1}\times_{M_n}Y(\Gamma)\times T(M)\mapsto Y(\Gamma)\times 
T(M)}(c_{total}(\tau), C-{\bf M}(E)E-\sum_{i\leq p} e_{k_i})$$

 for some $\tau\in K_0(Y(\Gamma)\times T(M))$ represented by a locally
 free sheaf on some zariski open subset of $Y(\Gamma)\times T(M)$. 

The mixed family invariant is identically zero
 when there exists an type $I$ exceptional class $e_{i}$ 
with $0>(C-{\bf M}(E)E)\cdot e_{i}>e_{i}^2$.
\end{theo}

\medskip

\begin{rem}\label{rem; question1}
The statement in the theorem still holds without the {\bf Simplifying Assumption},
but it is beyond the reach of the current discussion. 
 The element $\tau\in K_0(Y(\Gamma)\times T(M))$ admits a locally free 
representative and its associated vector bundle
 was called $\kappa$, the residual relative
 obstruction bundle on page 443 of [Liu1]. 
\end{rem}

\medskip

\noindent Proof of theorem \ref{theo; question1}: 
Because $Y_1={\cal M}_{C-{\bf M}(E)E-\sum_{i\leq p}e_{k_i}}\times_{M_n}Y(\Gamma)
\equiv 
Z(s_{canon}^{\circ})\times_{M_n}Y(\Gamma)$, our goal is to identify 
${j_{Y_1}}_{\ast}{\bf Z}_{Y_1}(s_{canon})\cap c_1(i_Y^{\ast}{\bf H})^{{C^2-C\cdot 
c_1({\bf K}_M)\over 2}+p_g-\sum_{i\leq n} {m_i^2+m_i\over 2}}$ with 
some mixed family invariant of $C-{\bf M}(E)E-\sum_{i\leq p}e_{k_i}$ through a
 nine steps process.

\medskip

\noindent Step I: By definitions of 
$m=dim_{\bf C}{\bf P}({\cal V}_{canon}^{\circ})
=dim_{\bf C}B+q(M)+rank_{\bf C}{\bf V}_{canon}^{\circ}-1$ and
 $e=rank_{\bf C}{\bf W}_{canon}$, 
${\bf Z}_{Y_1}(s_{canon})=\{c(i_{Y_1}^{\ast}({\bf H}\otimes
\pi^{\ast}_{{\bf P}({\bf V}_{canon}^{\circ})}{\bf W}_{canon}))\cap 
s(Y_1, {\bf P}({\cal V}_{canon}^{\circ}))\}_{m-e}$, and we desire to
 transform it into a more recognized form.

 Denote the projection map $Y_1\mapsto Y(\Gamma)$ by $\pi_{Y_1}$.
 Firstly, the inclusion $Y_1\subset X={\bf P}({\bf V}_{canon}^{\circ})$
 factors through $Y_1\subset X\times_{M_n}Y(\Gamma)\subset X$
 and it induces the short exact sequence of cones (see page 72, example
 4.1.6. in [F] for the definition)

$$0\mapsto\pi_{Y_1}^{\ast}
{\bf N}_{Y(\Gamma)}M_n\mapsto {\bf C}_{Y_1}X\mapsto {\bf C}_{Y_1}X
\times_{M_n}Y(\Gamma)\mapsto 0$$

and it implies the equality

 $$s(Y_1, X)=s(\pi_{Y_1}^{\ast}
{\bf N}_{Y(\Gamma)}M_n)\cap s(Y_1, X\times_{M_n}Y(\Gamma)).$$

 We plug this identity into the defining formula of ${\bf Z}_{Y_1}(s_{canon})$
 and get

$${\bf Z}_{Y_1}(s_{canon})=\{c(i_{Y_1}^{\ast}({\bf H}\otimes
\pi^{\ast}_{{\bf P}({\bf V}_{canon}^{\circ})}{\bf W}_{canon}))\cap 
s(\pi_{Y_1}^{\ast}{\bf N}_{Y(\Gamma)}M_n)\cap 
s(Y_1, X\times_{M_n}Y(\Gamma))\}_{m-e}.$$

\medskip

\noindent Step II:  By lemma \ref{lemm; cone} 
the cone $\cup_{i\not=0}{\bf C}_{\rho_i}$ (excluding the
irreducible component ${\bf O}^{\circ}$) is the locus over which the
 map $g$ from the total space of 
${\bf H}\otimes \pi^{\ast}_{{\bf P}({\bf V}_{canon}^{\circ})}
{\bf W}_{canon}^{\circ}|_{Y(\Gamma)\times T(M)}$ to the total space of 
${\bf H}\otimes \pi^{\ast}_{{\bf P}({\bf V}_{canon}^{\circ})}
{\bf W}_{canon}|_{Y(\Gamma)\times T(M)}$ 
 fails to be injective. 

 Denote the projection map from $Y_1$ to $Y(\Gamma)\times T(M)\subset M_n\times 
 T(M)$ by $k_{Y_1}$.

The assumption that $Y_1\cap Y_2=\emptyset$ implies the injection, 

$$0\mapsto i_{Y_1}^{\ast}{\bf H}\otimes 
k_{Y_1}^{\ast}{\bf W}_{canon}^{\circ}
\mapsto i_{Y_1}^{\ast}{\bf H}\otimes k_{Y_1}^{\ast}{\bf W}_{canon}$$

because $Y_2\subset {\bf s}_{canon}^{\circ}$ the locus over which the
 map $g|_{{\bf s}_{canon}^{\circ}}$ fails to be injective, is disjoint from
 $Y_1$.

 By using $\pi_{{\bf P}({\bf V}_{canon}^{\circ})}\circ i_{Y_1}=k_{Y_1}$, we
may rewrite $c(i_{Y_1}^{\ast}({\bf H}\otimes 
\pi_{{\bf P}({\bf V}_{canon}^{\circ})}^{\ast}
{\bf W}_{canon}/))$ as
 $c(i_{Y_1}^{\ast}{\bf H}\otimes k_{Y_1}^{\ast}{\bf W}_{canon}^{\circ})
\cap c(i_{Y_1}^{\ast}{\bf H}\otimes k_{Y_1}^{\ast}
{\bf W}_{canon}/i_{Y_1}^{\ast}{\bf H}\otimes 
k_{Y_1}^{\ast}{\bf W}_{canon}^{\circ})$.

\medskip

\noindent Step III: As a direct consequence of our {\bf Simplifying Assumption} 
on page \pageref{assum} that the cone $\cup_{i\not=0} {\bf C}_{\rho_i}$ is
 disjoint from $Y_1$, the sheaf 
$i_{Y_1}^{\ast}{\cal R}^0\pi_{\ast}\bigl(
{\cal O}_{\sum_{i\leq p} \Xi_i}(-{\bf M}(E)E)\otimes {\cal E}_C\bigr)$
 is trivial and the sheaf 
 $i_{Y_1}^{\ast}{\cal R}^1\pi_{\ast}\bigl(
{\cal O}_{\sum_{i\leq p} \Xi_i}(-{\bf M}(E)E)\otimes {\cal E}_C\bigr)$ is 
 locally free on 
$Y_1={\cal M}_{C-{\bf M}(E)E-\sum_{i\leq p}e_{k_i}}\times_{M_n}Y(\Gamma)$.
 We may denote the algebraic vector bundle associated with the first derived
 image sheaf by ${\bf G}\mapsto Y_1$.
 By the commutative diagram in the statement of proposition \ref{prop; com},
 $i_{Y_1}^{\ast}
({\bf H}\otimes {\bf G}) \cong i_{Y_1}^{\ast}{\bf H}\otimes k_{Y_1}^{\ast}
{\bf W}_{canon}/i_{Y_1}^{\ast}{\bf H}\otimes 
k_{Y_1}^{\ast}{\bf W}_{canon}^{\circ}$.

We notice that $i_{Y_1}:Y_1\subset X$ factors through 
 $\bar{i}_{Y_1}:Y_1\subset X\times_{M_n}Y(\Gamma)=X'$.

  Combining the conclusion of step II and the above identification, we
 may rewrite 

 $$\hskip -1.2in {\bar{i}_{Y_1\ast}}{\bf Z}_{Y_1}(s_{canon})
=\bar{i}_{Y_1\ast}
\{c(i_{Y_1}^{\ast}{\bf H}\otimes k_{Y_1}^{\ast}{\bf W}_{canon}^{\circ})
\cap s(Y_1, X')
\cap c(i_{Y_1}^{\ast}{\bf H}\otimes {\bf G})
\cap s(\pi_{Y_1}^{\ast}{\bf N}_{Y(\Gamma)}M_n)\}_{m-e}$$
$$\hskip -1.2in =\sum_{-\infty<r\leq dim_{\bf C}X'}\bar{i}_{Y_1\ast}\{
c(i_{Y_1}^{\ast}{\bf H}\otimes k_{Y_1}^{\ast}{\bf W}_{canon}^{\circ})
\cap s(Y_1, X\times_{M_n}Y(\Gamma))\}_{m-e+r}\cap \bar{i}_{Y_1\ast}
\{c(i_{Y_1}^{\ast}({\bf H}\otimes {\bf G}))\cap 
s(\pi_{Y_1}^{\ast}{\bf N}_{Y(\Gamma)}M_n)\}_{dim_{\bf C}X'-r}.$$

 \medskip

\noindent Step IV: We recall that 
$Z(s_{canon}^{\circ})\times_{M_n}Y(\Gamma)=Y_1$. Then by taking
 $X'=X\times_{M_n}Y(\Gamma)$, $i_{X'}:X'\subset X$, and 
$E'=i_{X'}^{\ast}{\bf H}\otimes 
\pi_X^{\ast}{\bf W}_{canon}^{\circ}$
we find that ($m'=dim_{\bf C}X'$, $e'=rank_{\bf C}{\bf W}_{canon}^{\circ}$)

 $${\bf Z}(s_{canon}^{\circ}|_{X'})=
\{c(i_{Y_1}^{\ast}({\bf H}\otimes \pi_{X}^{\ast}{\bf W}_{canon}^{\circ}))
\cap s(Y_1, X')\}_{m'-e'}$$
$$=\{c(i_{Y_1}^{\ast}({\bf H}\otimes \pi_{X}^{\ast}{\bf W}_{canon}^{\circ}))
\cap s(Y_1, X')\}_{m-e-rank_{\bf C}{\bf N}_{Y(\Gamma)}M_n+rank_{\bf C}{\bf G}}$$

 by using $m'=dim_{\bf C}X-(dim_{\bf C}M_n-dim_{\bf C}
Y(\Gamma))=m-codim_{\bf C}Y(\Gamma)=m-rank_{\bf C}{\bf N}_{Y(\Gamma)}M_n$ and 
 $e'=rank_{\bf C}{\bf W}_{canon}-rank_{\bf C}{\bf G}=e-rank_{\bf C}{\bf G}$. 

\medskip

Then the summation and its range of $r$ in the final expression in step III 
can be changed into

$$\hskip -1.2in =\sum_{rank_{\bf C}{\bf G}-rank_{\bf C}{\bf N}_{Y(\Gamma)}M_n
\leq r\leq dim_{\bf C}X'}\bar{i}_{Y_1\ast}\{
c(i_{Y_1}^{\ast}{\bf H}\otimes k_{Y_1}^{\ast}{\bf W}_{canon}^{\circ})
\cap s(Y_1, X')\}_{m-e+r}\cap \bar{i}_{Y_1\ast}
\{c(i_{Y_1}^{\ast}{\bf H}\otimes {\bf G})
\cap s(\pi_{Y_1}^{\ast}{\bf N}_{Y(\Gamma)}M_n)\}_{dim_{\bf C}X'-r},$$

by applying the grading constraint from proposition \ref{prop; local} to $X'$, 
$E'=i_{X'}^{\ast}\bigl({\bf H}\otimes {\bf W}_{canon}^{\circ}\bigr)$ and
 $s_{canon}^{\circ}|_{X'}$.

\medskip

 To simply this final expression, some preparation in step V, VI, VII are
 necessary.

\medskip

\subsection{\bf The Canonical Algebraic Kuranishi Model of Type I Exceptional 
Curves}\label{subsection; typeI}

\bigskip

\noindent Step V: We construct the canonical algebraic Kuranishi model of
 a type I exceptional class $e_i=E_i-\sum_{j_i}E_{j_i}$ as the following.

 Let $\Gamma_{e_i}$ denote the $n$-vertex admissible graph such that 

\noindent (i). the direct descendents of the $i$-th vertex
 are exactly all the $j_i$-th vertexes. 

\noindent (ii). the $i$-th vertex is the unique vertex among the $n$ vertexes which
 has any direct descendent.

\medskip

\begin{figure}
\centerline{\epsfig{file=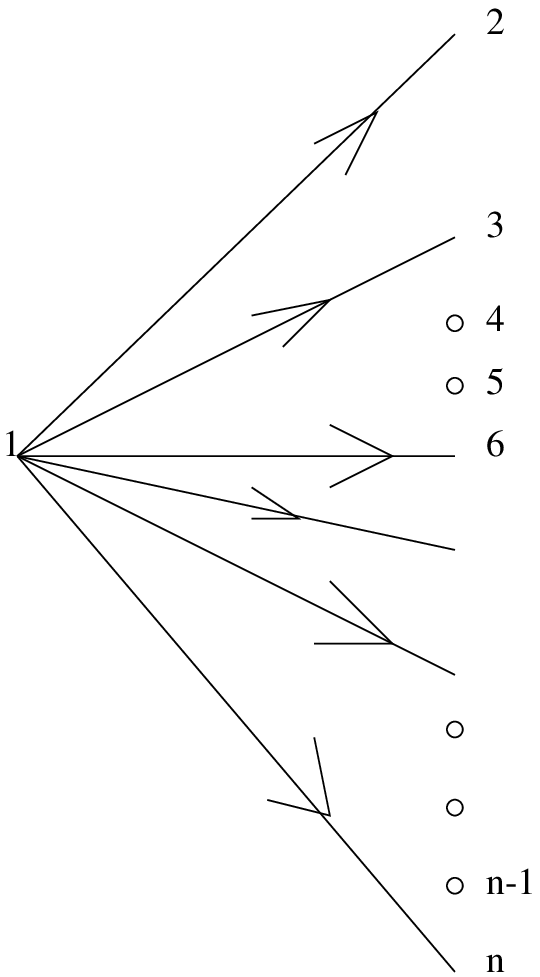,height=4cm}}
\centerline{fig.1}
\end{figure}

 See figure 1 for an example. Then $Y(\Gamma_{e_i})\subset M_n$ is the locus
 over which the class $e_i$ becomes effective and $codim_{\bf C}Y(\Gamma_{e_i})$
 is equal to the number of $1$-edges in the graph $\Gamma_{e_i}$.  

\medskip

\begin{lemm}\label{lemm; normal}
 There exists a canonical sheaf isomorphism on the
 normal sheaf ${\cal N}_{Y(\Gamma_{e_i})}M_n\cong
 {\cal R}^1\pi_{\ast}\bigl({\cal O}_{\Xi_i}(E_i-\sum_{j_i}E_{j_i})\bigr)$.
\end{lemm}

\medskip

\noindent Proof: The sheaf short exact sequence 

$$0\mapsto {\cal O}_{M_{n+1}\times_{M_n}Y(\Gamma_{e_i})}\mapsto 
{\cal O}_{M_{n+1}\times_{M_n}Y(\Gamma_{e_i})}(E_i-\sum_{j_i}E_{j_i})
\mapsto {\cal O}_{\Xi_i}(E_i-\sum_{j_i}E_{j_i})\mapsto 0$$

 implies the following short exact sequence

$$\hskip -1.2in 
0\mapsto {\cal R}^1\pi_{\ast}\bigl(
{\cal O}_{M_{n+1}\times_{M_n}Y(\Gamma_{e_i})}\bigr)\mapsto 
{\cal R}^1\pi_{\ast}\bigl({\cal O}_{M_{n+1}\times_{M_n}Y(\Gamma_{e_i})}(E_i-
\sum_{j_i}E_{j_i})\bigr)\mapsto 
{\cal R}^1\pi_{\ast}\bigl({\cal O}_{\Xi_i}(E_i-\sum_{j_i}E_{j_i})\bigr)\mapsto 0.$$

We notice ${\cal R}^1\pi_{\ast}\bigl(
{\cal O}_{M_{n+1}\times_{M_n}Y(\Gamma_{e_i})}\bigr)\cong
 H^1(M, {\cal O}_M)\otimes {\cal O}_{Y(\Gamma_{e_i})}$.

\medskip

 Consider the sheaf short exact sequence 

$$\hskip -1in 0\mapsto {\cal O}_{M_{n+1}}(E_i-\sum_{j_i}E_{j_i})\mapsto 
{\cal O}_{M_{n+1}}(E_i)\mapsto
 {\cal O}_{\sum_{j_i}E_{j_i}}(E_i)\mapsto 0$$

and take the right derived long exact sequence along $\pi=f_n:M_{n+1}\mapsto M_n$,
we get the following five-term sheaf exact sequence over $M_n$,

$$\hskip -1in 0\mapsto {\cal R}^0\pi_{\ast}\bigl({\cal O}_{M_{n+1}}(E_i-\sum_{j_i}
E_{j_i})\bigr)\mapsto {\cal O}_{M_n}\mapsto {\cal R}^0\pi_{\ast}\bigl(
{\cal O}_{\sum_{j_i}E_{j_i}}(E_i)\bigr)\mapsto 
{\cal R}^1\pi_{\ast}\bigl({\cal O}_{M_{n+1}}(E_i-
\sum_{j_i}E_{j_i})\bigr)$$

$$\hskip -1in 
\mapsto {\cal R}^1\pi_{\ast}\bigl({\cal O}_{M_{n+1}}(E_i)\bigr)\mapsto 0.$$

\medskip

 Similar to the discussion in proposition 5.3 of [Liu3], the sheaf 
${\cal R}^0\pi_{\ast}\bigl(
{\cal O}_{\sum_{j_i}E_{j_i}}(E_i)\bigr)$ is locally free with its rank equal to
 the number of $1$-edges in $\Gamma_{e_i}$. It is the canonical 
obstruction bundle of the type I class $e_i$.

By using 
$$\hskip -.3in {\cal O}_{Y(\Gamma_{e_i})}\otimes H^1(M, {\cal O}_M)\cong 
{\cal R}^1\pi_{\ast}\bigl(
{\cal O}_{M_{n+1}\times_{M_n}Y(\Gamma_{e_i})}\bigr)\cong
 {\cal R}^1\pi_{\ast}\bigl({\cal O}_{M_{n+1}}(E_i)\bigr)|_{Y(\Gamma_{e_i})},$$
 
it is easy to modify the above five-term exact sequence to the
following canoncial algebraic Kuranishi model of $e_i=E_i-\sum_{j_i}E_{j_i}$,

$$\hskip -1in 0\mapsto 
{\cal R}^0\pi_{\ast}\bigl({\cal O}_{M_{n+1}}(E_i-\sum_{j_i}
E_{j_i})\bigr)\mapsto {\cal O}_{M_n}\mapsto {\cal R}^0\pi_{\ast}\bigl(
{\cal O}_{\sum_{j_i}E_{j_i}}(E_i)\bigr)\mapsto 
{\cal R}^1\pi_{\ast}\bigl({\cal O}_{M_{n+1}}(E_i-
\sum_{j_i}E_{j_i})\bigr)/{\cal R}^1\pi_{\ast}\bigl(
{\cal O}_{M_{n+1}}\bigr)\mapsto 0.$$

The
 morphism ${\cal O}_{M_n}\cong 
{\cal R}^0\pi_{\ast}\bigl({\cal O}_{M_{n+1}}(E_i)\bigr)
\mapsto {\cal R}^0\pi_{\ast}\bigl(
{\cal O}_{\sum_{j_i}E_{j_i}}(E_i)\bigr)$ defines
 a canonical section whose zero locus
 $=Y(\Gamma_{e_i})$. By proposition 4.3.
 in [Liu1], the space $Y(\Gamma_{e_i})\subset M_n$
 is smooth and its codimension in $M_n$ matches with the rank of
 the locally free sheaf ${\cal R}^0\pi_{\ast}\bigl(
{\cal O}_{\sum_{j_i}E_{j_i}}(E_i)\bigr)$. 
Thus the algebraic section defining $Y(\Gamma_{e_i})$
 is regular.  An immediate consequence of the
 regularity of the section is the following 
 identification of the normal sheaf 

$$\hskip -1in {\cal N}_{Y(\Gamma_{e_i})}M_n\cong  {\cal R}^0\pi_{\ast}\bigl(
{\cal O}_{\sum_{j_i}E_{j_i}}(E_i)\bigr)|_{Y(\Gamma_{e_i})} \cong
{\cal R}^1\pi_{\ast}\bigl({\cal O}_{M_{n+1}}(E_i-
\sum_{j_i}E_{j_i})\bigr)/{\cal R}^1\pi_{\ast}\bigl(
{\cal O}_{M_{n+1}}\bigr)|_{Y(\Gamma_{e_i})}\cong 
{\cal R}^1\pi_{\ast}\bigl({\cal O}_{\Xi_i}(E_i-\sum_{j_i}E_{j_i})\bigr).$$

The lemma is proved. $\Box$

\medskip

\subsection{\bf Short Exact Sequences on ${\cal N}_{Y(\Gamma_{e_i})}M_n$}
\label{subsection; short}

\bigskip

\noindent Step VI: Consider the ${\bf P}^1$ fibration $\Xi_i\mapsto 
 Y(\Gamma_{e_i})$ and the projection map $\pi:T(M)\times \Xi_i\mapsto 
 T(M)\times Y(\Gamma_{e_i})$ induced from 
$\pi:M_{n+1}\times_{M_n}Y(\Gamma_{e_i})\mapsto Y(\Gamma_{e_i})$.
 Then we have the following proposition,

\medskip

\begin{prop}\label{prop; residue}
 Let $e_i$ be a type $I$ exceptional class satisfying $e_i\cdot (C-{\bf M}(E)E)<0$.
 Suppose that the intersection pairing 
$e_i\cdot (C-{\bf M}(E)E)\leq e_i^2<0$, then there exists a
 relative effective divisor $\Delta_i\subset \Xi_i\mapsto Y(\Gamma_{e_i})$
 of relative degree $({\bf M}(E)E+e_i)\cdot e_i\geq 0$, an invertible
 sheaf $\pi^{\ast}{\cal Q}_i$ on $\Xi_i$ 
pulled-back from $Y(\Gamma_{e_i})$ and a short 
exact sequence of locally free sheaves,

$$\hskip -1in 0\mapsto  {\cal R}^0\pi_{\ast}\bigl({\cal O}_{\Delta_i}(E_i-
\sum_{j_i}E_{j_i})\bigr)\otimes {\cal Q}_i 
\mapsto {\cal R}^1\pi_{\ast}\bigl({\cal O}_{\Xi_i}(-{\bf M}(E)E)\bigr)\mapsto 
{\cal R}^1\pi_{\ast}\bigl({\cal O}_{\Xi_i}(E_i-\sum_{j_i}E_{j_i})
\bigr)\otimes {\cal Q}_i \mapsto 0$$
 
exact on a zariski open subset of $Y(\Gamma_{e_i})$.

Suppose that the intersection pairing satisfies $0>e_i\cdot (C-{\bf M}(E)E)>e_i^2$,
 then there exists a relative effective divisor $\Delta_i\subset \Xi_i\mapsto
 Y(\Gamma_{e_i})$ of relative degree $-({\bf M}(E)E+e_i)\cdot e_i>0$, an 
invertible sheaf $\pi^{\ast}{\cal Q}_i$ on $\Xi_i$ pulled-back from $Y(\Gamma_{e_i})$ and 
 a short exact sequence of locally free sheaves,

$$\hskip -1in 0\mapsto {\cal R}^0\pi_{\ast}\bigl({\cal O}_{\Delta_i}(E_i
-\sum_{j_i}E_{j_i})\bigr)\otimes {\cal Q}_i 
\mapsto {\cal R}^1\pi_{\ast}\bigl({\cal O}_{\Xi_i}(E_i-\sum_{j_i}E_{j_i})
\bigr) \otimes {\cal Q}_i  \mapsto 
{\cal R}^1\pi_{\ast}\bigl({\cal O}_{\Xi_i}(-{\bf M}(E)E)\bigr)  \mapsto 0$$

 exact on a Zariski open subset of $Y(\Gamma_{e_i})$.
\end{prop}

\medskip

The reader should notice the similarity between proposition \ref{prop; residue}
 and proposition \ref{prop; Q} in subsection \ref{subsection; GT}.
 By lemma \ref{lemm; normal} in the previous subsection, these exact sequences
 can be viewed as exact sequences about ${\cal N}_{Y(\Gamma_{e_i})}M_n$.

\medskip

\noindent Proof of the proposition: Denote the set of all descendent indexes
 of $i$ by $J_i$.
Firstly, notice that the effective
 divisors $E_{j_i}\subset M_{n+1}$, $j_i\in J_i$
 restrict to cross sections on $\Xi_i\mapsto Y(\Gamma_{e_i})$. As 
 $\Xi_i$ is a ${\bf P}^1$ fibration, the invertible 
 sheaves ${\cal O}_{\Xi_i}(E{j_i}|_{\Xi_i})$ for different $j_i\in J_i$ are
 equivalent after tensoring invertible sheaves pulled back from the base 
$Y(\Gamma_{e_i})$. 
Define $q=({\bf M}(E)E+e_i)\cdot e_i\in {\bf Z}$.
Suppose that $q\geq 0$, fix one $l\in J_i$ and define 
$\Delta_i=qE_l|_{\Xi_i}+\sum_{j\not\in J_i\cup \{i\}}m_jE_j$.
 Then ${\cal O}_{\Xi_i}(-{\bf M}(E)E)$ and
 ${\cal O}_{\Xi_i}(-\Delta_i+(E_i-\sum_{j_i\in J_i}E_{j_i}))$ have the same 
relative degrees along $\Xi_i\mapsto Y(\Gamma_{e_i})$ 
and the former is equivalent to the latter after tensoring the
 latter sheaf by 
 some invertible sheaf $\pi^{\ast}{\cal Q}_i$ 
pulled-back from the base $Y(\Gamma_{e_i})$.

 Then by tensoring the defining short exact sequence of $\Delta_i$ by 
${\cal O}_{\Xi_i}(E_i-\sum_{j_i}E_{j_i})\otimes \pi^{\ast}{\cal Q}_i$,

 one gets 

$$\hskip -1in 0\mapsto {\cal O}_{\Xi_i}(-{\bf M}(E)E) 
\mapsto {\cal O}_{\Xi_i}(E_i-\sum_{j_i\in J_i}E_{j_i})\otimes \pi^{\ast}{\cal Q}_i 
\mapsto {\cal O}_{\Delta_i}(E_i-\sum_{j_i\in J_i}E_{j_i})\otimes
 \pi^{\ast}{\cal Q}_i 
\mapsto 0.$$ 

 By taking the derived long exact sequence along $
\pi:\Xi_i\mapsto Y(\Gamma_{e_i})$ and restrict to the complement of the
 support of ${\cal R}^0\pi_{\ast}\bigl({\cal O}_{\Xi_i}(-{\bf M}(E)E)
\bigr)$, we get
 the desired sheaf short exact sequence stated in the proposition.

 Notice that the choices of $\Delta_i$, ${\cal Q}_i$ are not unique.

 The proof of the $q<0$ case is rather 
similar and we leave it to the reader. $\Box$

\medskip

\subsection{\bf Some Vanishing Results about the Local Contribution}
\label{subsection; vani}

\bigskip

\noindent Step VII: In this step, we derive a vanishing lemma which will be used
 in the next few steps.

\medskip

\begin{lemm}\label{lemm; equal}
Let $V\mapsto X$ be a vector bundle of rank $v$ and let $s:X\mapsto V$ be
 a regular section of $V$ with $i:Z(s)\subset X$, $codim_{\bf C}Z(s)=v$.
 Let ${\bf Q}$ be a line bundle over $X$, then
 $0=c_r(i^{\ast}({\bf Q}\otimes V)-i^{\ast}V)=\{c(i^{\ast}({\bf Q}\otimes V))
\cap  s(i^{\ast})V\}_{dim_{\bf C}Z(s)-r}\in {\cal A}_{dim_{\bf C}Z(s)-r}(Z(s))$ 
for all $r\in {\bf N}$.
\end{lemm}

\medskip

\noindent A sketch of the proof: 
Firstly, assume that ${\bf Q}$ is effective and consider
 a fixed section $s_{D_{\bf Q}}$ defining the
 effective divisor $D_{\bf Q}$. Then $s\otimes s_{D_{\bf Q}}$
 is a section of $V\otimes {\bf Q}$ and $Z(s)\subset Z(s\otimes s_{D_{\bf Q}})$.

 By the regularity condition on $s$, ${\bf N}_{Z(s)}X\cong i^{\ast}V$.
 Blowing up $Z(s)\subset X$ into a smooth divisor $D$.

By a direct computation following [F] page 161-162, equations (1), (2), (3),
 one finds that

$$c(i^{\ast}({\bf Q}\otimes V)\cap s(Z(s), X)=
c({\bf Q}\otimes V)\cap s(Z(s\otimes s_{D_{\bf Q}}), X)-
c({\cal O}(-D)\otimes {\bf Q}\otimes V)\cap s(D_{\bf Q}, X).$$

 By the grading constraint from 
proposition \ref{prop; local}, the degree 
 $dim_{\bf C}X-v-r=dim_{\bf C}Z(s)-r$ pieces of both terms on the right hand side
 vanish for all $r\in {\bf N}$.

 Thus $c_r(i^{\ast}({\bf Q}\otimes V)-i^{\ast}V)=0$ for all $r\in {\bf N}$.

\medskip

 Secondly when ${\bf Q}$ is not effective, 
write ${\bf Q}={\bf Q}_1\otimes {\bf Q}_2^{-1}$,
 where ${\bf Q}_1, {\bf Q}_2$ are both effective. 
 This is always possible as we can twist ${\bf Q}$ by a high power of
 ample line bundle ${\bf D}$ to make both ${\bf Q}_1={\bf Q}\otimes {\bf D}^l$ and
 ${\bf Q}_2={\bf D}^l$ effective for large enough $l\gg 0$. 

We know that 
$f_r(m, n)=c_r(i^{\ast}({\bf Q}_1^m\otimes{\bf Q}_2^n\otimes V)-i^{\ast}V)=0$
 for all $m, n, r\in {\bf N}$ because ${\bf Q}_1^m\otimes {\bf Q}_2^n$ is
 effective.

 On the other hand, $c_1({\bf Q}_1^m\otimes {\bf Q}_2^n)=m\cdot c_1({\bf Q}_1)
+n\cdot c_1({\bf Q}_2)$ and $f_r(m, n)$, being an algebraic 
combination of Chern classes
 of $i^{\ast}{\bf Q}$ and $i^{\ast}V$,
 must be a polynomial expression in terms
  of $m$ and  $n$. Then the polynomial in $n$, $f_r(1, n)$
 has an infinite number of roots, 
$f_r(1, n)\equiv 0$ for all $n\in {\bf Z}$. In particular, $f_r(1, -1)=0$ and thus

$c_r(i^{\ast}({\bf Q}\otimes V)-i^{\ast}V)=
c_r(i^{\ast}({\bf Q}_1\otimes {\bf Q}_2^{-1}\otimes V)-i^{\ast}V)=0$. $\Box$

\medskip

\noindent Step VIII: After the preparation in step VI and step VII, 
we continue the discussion from step $IV$ and show that

\begin{prop}\label{prop; zero}
 Suppose that $e_i\cdot (C-{\bf M}(E)E)>e_i^2$ for some $i$, $1\leq i\leq p$,
 then local contribution of ${\cal M}_{C-{\bf M}(E)E-\sum_{i\leq p}e_i}
\times_{M_n}Y(\Gamma)$ 
 to the family invariant 
${\cal AFSW}_{M_{n+1}\times T(M)\mapsto M_n\times T(M)}(1, C-{\bf M}(E)E)$ 
defined by (consult page \pageref{yone}

$${j_{Y_1}}_{\ast}{\bf Z}_{Y_1}(s_{canon})\cap 
 c_1(i_Y^{\ast}{\bf H})^{dim_{\bf C}M_n+{C^2-C\cdot 
c_1({\bf K}_M)\over 2}+p_g-\sum_{i\leq n}{m_i^2+m_i\over 2}}\in {\bf Z}$$ 

 vanishes.
\end{prop}

\medskip

\noindent Proof: For notational simplicity, we may assume that
 $e_1\cdot (C-{\bf M}(E)E)>e_1^2$ by permuting the indexes 
$\{1, 2, \cdots, p\}$ (if it is necessary).

 Firstly, we notice that for $\underline{C}=C-{\bf M}(E)E$, we have an
 equality,  
 $$ 0<e_1\cdot (C-{\bf M}(E)E-e_1)=e_1\cdot (\underline{C}-e_i)
 =\bigl({\underline{C}^2-c_1({\bf K}_{M_{n+1}/M_n})\cdot \underline{C}\over 2}
+p_g\bigr)$$
$$-\bigl({(\underline{C}-e_1)^2-c_1({\bf K}_{M_{n+1}/M_n})\cdot 
(\underline{C}-e_1)\over 2}+p_g+
{e_1^2-c_1({\bf K}_{M_{n+1}/M_n})\cdot e_1\over 2}\bigr).$$

Thus the family dimension of the class $\underline{C}$ is strictly larger than the
 family dimension of the splitting $\underline{C}-e_1$ and $e_1$.
 We will argue the vanishing of the intersection number by the negativity of the
 dimension count.
 
\medskip

 Firstly, take ${\bf W}$ to be the algebraic vector bundle associated to
 ${\cal R}^0\pi_{\ast}\bigl({\cal O}_{{\bf M}(E)E+E_1-\sum_{j_1}E_{j_1}}
\otimes {\cal E}_C\bigr)$.

 Then we argue 

$$\hskip -1.2in 
\{c(i_{Y_1}^{\ast}{\bf H}\otimes \pi^{\ast}_{{\bf P}({\bf V}_{canon}^{\circ})}
{\bf W})\cap 
s(Y_1, X'\}_{dim_{\bf C}X-
codim_{\bf C}Y(\Gamma)-rank_{\bf C}{\bf W}}\cap c_1(i_{Y_1}^{\ast}{\bf H})^{
{\underline{C}^2-c_1({\bf K}_{M_{n+1}/M_n})\cdot \underline{C}\over 2}+
p_g+dim_{\bf C}M_n}$$

 vanishes.

 It is because 

$$\hskip -1.2in dim_{\bf C}X-codim_{\bf C}Y(\Gamma)-rank_{\bf C}{\bf W}
={(\underline{C}-e_1)^2-c_1({\bf K}_{M_{n+1}/M_n})\cdot (\underline{C}-e_1)\over 2}
+p_g+dim_{\bf C}M_n+
\sum_{i\leq p}{e_i^2-c_1({\bf K}_{M_{n+1}/M_n})\cdot e_i\over 2}$$
$$\hskip -1.2in 
<{\underline{C}^2-c_1({\bf K}_{M_{n+1}/M_n})\cdot \underline{C}\over 2}+
p_g+dim_{\bf C}M_n.$$

Secondly, we would like to transform $\alpha=\{c(i_{Y_1}^{\ast}{\bf H}\otimes 
{\bf W})\cap s(Y_1, X')\}_{dim_{\bf C}X-
codim_{\bf C}Y(\Gamma)-rank_{\bf C}{\bf W}}\in 
{\cal A}_{dim_{\bf C}X-codim_{\bf C}Y(\Gamma)-rank_{\bf C}{\bf W}}
(Y_1)$ to the original local
expression $\{c(i_{Y_1}^{\ast}
\pi^{\ast}_{{\bf P}({\bf V}_{canon}^{\circ})}{\bf W}_{canon})\cap 
s(\pi_{Y_1}^{\ast}{\bf N}_{Y(\Gamma)}M_n)\cap 
s(Y_1, X\times_{M_n}Y(\Gamma))\}_{m-e}$.

\medskip

 By proposition \ref{prop; residue} and the fact
 $e_1\cdot (C-{\bf M}(E)E)>e_1^2$, there exists a short exact sequence
 on ${\cal R}^1\pi_{\ast}\bigl({\cal O}_{\Xi_1}(E_1-\sum_{j_1}E_{j_1})\bigr)$.
 We may tensor it with ${\cal E}_C$ and get

$$\hskip -1in 0\mapsto {\cal R}^0\pi_{\ast}\bigl({\cal O}_{\Delta_2}(E_1
-\sum_{j_1}E_{j_1})\bigr)\otimes {\cal Q}_1\otimes {\cal E}_C 
\mapsto {\cal R}^1\pi_{\ast}\bigl({\cal O}_{\Xi_1}(E_1-\sum_{j_1}E_{j_1})
\bigr) \otimes {\cal Q}_1 \otimes {\cal E}_C \mapsto 
{\cal R}^1\pi_{\ast}\bigl({\cal O}_{\Xi_1}(-{\bf M}(E)E)\otimes {\cal E}_C
\bigr)  \mapsto 0.$$

 Let us set a few notations before moving forward.

 By proposition \ref{prop; residue} the sheaf 
${\cal R}^0\pi_{\ast}\bigl({\cal O}_{\Delta_i}(E_i
-\sum_{j_i}E_{j_i})\bigr)\otimes {\cal Q}_i\otimes {\cal E}_C$
 is locally free on $U_i$, the complement of the support of
 ${\cal R}^0\pi_{\ast}\bigl({\cal O}_{\Xi_i}(-{\bf M}(E)E)\bigr)$. \label{ui}
Denote  
 the corresponding vector bundle by $\nu_i\mapsto U_i$. Denote the normal bundle 
 of $Y(\Gamma_{e_i})\subset M_n$ by ${\bf N}_i, 1\leq i\leq p$. 
Denote the inclusion $Y(\Gamma)\subset Y(\Gamma_{e_i})$ by $h_i$.

 Because $Y(\Gamma)$ is the locus over which all $e_i, 1\leq i\leq p$
 are effective, $Y(\Gamma)=\cap_{1\leq i\leq p} Y(\Gamma_{e_i})$ and
 ${\bf N}_{Y(\Gamma)}=\oplus_{1\leq i\leq p}h_i^{\ast}{\bf N}_i$.

  Define $U=U_1\cap Y(\Gamma)\subset Y(\Gamma)$ and set $\pi_{\bf F}:
{\bf F}=\pi_{X'}^{\ast}h_1^{\ast}\nu_1\otimes {\bf H}
\pi_{X'}^{\ast}\oplus_{2\leq i\leq p}h_i^{\ast}{\bf N}_i|_{X'\times_{Y(\Gamma)}U}
\mapsto X'\times_{Y(\Gamma)}U$
 to be the projection map of the vector bundle ${\bf F}$. Set
 $s_{\bf F}:X'\times_{Y(\Gamma)}U\mapsto {\bf F}$ to be the zero section map.

 Because $Y_1\subset X'\times_{Y(\Gamma)}U$, 
$\alpha\in {\cal A}_{dim_{\bf C}X-codim_{\bf C}
Y(\Gamma)-rank_{\bf C}{\bf W}}(Y_1)$
 defines a class in ${\cal A}_{dim_{\bf C}X-codim_{\bf C}
Y(\Gamma)-rank_{\bf C}{\bf W}}(X'\times_{Y(\Gamma)}U)$ and we abuse the
 notation slightly and denote it by the same symbol.
 Then by page 67, example 3.3.2. of [F], 
$c({\bf F})\cap \alpha=s_{\bf F}^{\ast}s_{{\bf F}\ast}\alpha$ and then
$$\alpha=\{\alpha\}_{dim_{\bf C}X-codim_{\bf C}
Y(\Gamma)-rank_{\bf C}{\bf W}}=s_{\bf F}^{\ast}s_{{\bf F}\ast}
\{s({\bf F})\cap \alpha\}_{dim_{\bf C}X-codim_{\bf C}
Y(\Gamma)-rank_{\bf C}{\bf W}+rank_{\bf C}{\bf F}}.$$ 

By the short exact sequence of $\nu_1$ in 
proposition \ref{prop; residue} and by lemma \ref{lemm; normal}, 
we have 

$$s(\pi_{X'}^{\ast}\nu_1\otimes {\bf H})\cap s\bigl(\pi_{X'}^{\ast}
(h_1^{\ast}{\cal R}^1\pi_{\ast}\bigl(
{\cal O}_{\Xi_i}(-{\bf M}(E)E)\otimes {\cal E}_C\bigr))\otimes 
{\bf H}|_{X'\times_{Y(\Gamma)}U}\bigr)=
s(\pi_{X'}^{\ast}(
h_1^{\ast}({\bf N}_1\otimes{\bf Q}_1\otimes {\bf E}_C))\otimes {\bf H}|_{X'
\times_{Y(\Gamma)}U})$$

$$=s(\pi_{X'}^{\ast}
h_1^{\ast}{\bf N}_1|_{X'\times_{Y(\Gamma)}U})
\cap c(\pi_{X'}^{\ast}h_1^{\ast}{\bf N}_1|_{X'\times_{Y(\Gamma)}U}-
\pi_{X'}^{\ast}h_1^{\ast}{\bf N}_1\otimes {\bf Q}_1\otimes
 {\bf E}_C\otimes {\bf H}|_{X'\times_{Y(\Gamma)}U})=s(\pi_{X'}^{\ast}
h_1^{\ast}{\bf N}_1|_{X'\times_{Y(\Gamma)}U}).$$

$$c(\pi_{X'}^{\ast}h_1^{\ast}{\bf N}_1|_{X'\times_{Y(\Gamma)}U}-
\pi_{X'}^{\ast}h_1^{\ast}{\bf N}_1\otimes {\bf Q}_1\otimes
 {\bf E}_C\otimes {\bf H}|_{X'\times_{Y(\Gamma)}U})=s(\pi_{X'}^{\ast}
h_1^{\ast}{\bf N}_1|_{X'\times_{Y(\Gamma)}U})\equiv 1$$

because $X\times_{M_n}Y(\Gamma_{e_1})$
 is the regular zero locus of the global reqular section of
 $\pi_{X}^{\ast}
{\cal R}^0\pi_{\ast}\bigl({\cal O}_{\sum_{j_1} E_{j_1}}(E_1)\bigr)$ on 
$X={\bf P}_{M_n\times T(M)}({\bf V}_{canon})$,
 lemma \ref{lemm; equal} 
(we take ${\bf Q}=
\pi_{X}^{\ast}(h_1^{\ast}{\bf Q}_1\otimes {\bf E}_C)\otimes {\bf H}$) and
 the fact $U\subset Y(\Gamma)\subset Y(\Gamma_{e_1})$.

Thus, 

$$\hskip -1.2in 
\{s(i_{Y_1}^{\ast}{\bf F})\cap \alpha\}_{dim_{\bf C}X+codim_{\bf C}Y(\Gamma)-
rank_{\bf C}{\bf W}+rank_{\bf C}{\bf F}}=
\{s(i_{Y_1}^{\ast}\pi_{X'}^{\ast}{\bf N}_{Y(\Gamma)})\cap 
 c(i_{Y_1}^{\ast}{\bf H}\otimes {\cal R}^1\pi_{\ast}
\bigl({\cal O}_{\Xi_1}(-{\bf M}(E)E)\otimes {\cal E}_C\bigr))$$ 
$$\cap c(i_{Y_1}^{\ast}({\bf H}\otimes \pi_{X}^{\ast}
{\bf W}))\cap s(Y_1, X')
 \}_{dim_{\bf C}X+codim_{\bf C}Y(\Gamma)-
rank_{\bf C}{\bf W}+rank_{\bf C}{\bf F}}$$

$$=\{s({\bf N}_{Y(\Gamma)})\cap c(i_{Y_1}^{\ast} {\bf H}\otimes 
\pi_{{\bf P}({\bf V}_{canon}^{\circ})}^{\ast}{\bf W}_{canon})\cap
 s(Y_1, X')\}_{dim_{\bf C}X-codim_{\bf C}Y(\Gamma)-
rank_{\bf C}{\bf W}+rank_{\bf C}{\bf F}}\in {\cal A}_{\ast}(Y_1),$$

due to the Chern classes identity 
$$c({\bf W}_{canon})=c({\bf W})\cap c({\cal R}^1\pi_{\ast}
\bigl({\cal O}_{\Xi_1}(-{\bf M}(E)E)\otimes {\cal E}_C\bigr))$$
 derived from a short exact sequence similar to the one in
 proposition \ref{prop; com} relating ${\cal W}_{canon}^{\circ}$ and
 ${\cal W}_{canon}$. 

\medskip

Thirdly, the integral grading of the last expression above is equal to

$$\hskip -1.2in dim_{\bf C}X+codim_{\bf C}Y(\Gamma)-
rank_{\bf C}{\bf W}+rank_{\bf C}{\bf F}=dim_{\bf C}X+\sum_{1\leq i\leq p}
{e_i^2-c_1({\bf K}_{M_{n+1}/M_n})\cdot e_i\over 2}-rank_{\bf C}{\bf W}-
\sum_{2\leq i\leq p}{e_i^2-c_1({\bf K}_{M_{n+1}/M_n})\cdot e_i\over 2}$$

$$\hskip -1.2in+rank_{\bf C}\nu_1 
=dim_{\bf C}X+{e_1^2-c_1({\bf K}_{M_{n+1}/M_n})\cdot e_1\over 2}
-rank_{\bf C}{\bf W}+rank_{\bf C}
\nu_1=dim_{\bf C}X-rank_{\bf C}{\bf W}_{canon}$$

$$\hskip -1.2in 
+\{{e_1^2-c_1({\bf K}_{M_{n+1}/M_n})\cdot e_1\over 2}+rank_{\bf C}
\nu_1+rank_{\bf C}{\cal R}^1\pi_{\ast}
\bigl({\cal O}_{\Xi_1}(-{\bf M}(E)E)\otimes {\cal E}_C\bigr)\}
=rank_{\bf C}X-rank_{\bf C}{\bf W}_{canon}=m-e,$$

 and the final expression

$$\{c(i_{Y_1}^{\ast}{\bf H}\otimes \pi_{{\bf P}({\bf V}_{canon}^{\circ})}^{\ast}
{\bf W}_{canon})\cap s(Y_1, X\times_{M_n}Y(\Gamma))\cap 
s(\pi_{Y_1}^{\ast}{\bf N}_{Y(\Gamma)})\}_{m-e}$$

 matches with the final expression in step I. Therefore, its cap product with 
$c_1(i_{Y_1}{\bf H})^{{\underline{C}^2-c_1({\bf K}_{M_{n+1}/M_n})\cdot 
\underline{C}\over 2}+p_g+dim_{\bf C}M_n}$ 
must be zero as well. $\Box$

\bigskip

\noindent Step IX: In this final step, we work with the situation that
 $0>e_i^2\geq e_i\cdot (C-{\bf M}(E)E)$, $1\leq i\leq p$. Define
 the class $\tau$ and identify the local family invariant contribution on $Y_1$
 with ${\cal AFSW}_{M_{n+1}\times_{M_n}Y(\Gamma)\times T(M)\mapsto Y(\Gamma)\times
 T(M)}(c_{total}(\tau), C-{\bf M}(E)E-\sum_{1\leq i\leq p}e_i)$.

\medskip

 Define $\tau$ to be the equivalence class represented by the sheaf 

$$\oplus_{1\leq i\leq p}h_i^{\ast}
\bigl({\cal R}^0\pi_{\ast}\bigl({\cal O}_{\Xi_i\cap 
\cup_{0\leq a<i}\Xi_a}(-{\bf M}(E)E)\otimes {\cal E}_C\bigr)\oplus 
 {\cal R}^0\pi_{\ast}\bigl({\cal O}_{\Delta_i}(-{\bf M}(E)E)\bigr)\otimes 
{\bf Q}_i\otimes {\cal E}_C\bigr),$$

locally free
 on $\cap_{1\leq i\leq p} U_i\times T(M)$, where $U_i\subset Y(\Gamma_{e_i})$
 stands for the Zariski open subset defined on page \pageref{ui} in step VIII.
 Because $\pi_{Y_1}:Y_1\mapsto Y(\Gamma)$ factors through $U\subset Y(\Gamma)$, 
 $c_{rank_{\bf C}\tau+r}(\pi_{Y_1}^{\ast}\tau)=0$ for $r\in {\bf N}$.

 The family moduli space 
${\cal M}_{C-{\bf M}(E)E-\sum e_i}\times_{M_n}Y(\Gamma)=Y_1\subset X'$ 
above $Y(\Gamma)$ is
 defined to be the zero locus $Z(s_{canon}^{\circ})\times_{M_n}Y(\Gamma)$. 
 By applying proposition \ref{prop; local}, the mixed family invariant 
${\cal AFSW}_{M_{n+1}\times_{M_n}Y(\Gamma)\times T(M)\mapsto Y(\Gamma)\times
 T(M)}(c_{total}(\tau), C-{\bf M}(E)E-\sum_{1\leq i\leq p}e_i)$ can be
 identified with the following expression involving localized top Chern class,

$$\hskip -1.2in
\bar{i}_{Y_1\ast}\bigl(\sum_{0\leq r\leq rank_{\bf C}\tau}
c_r(i_{Y_1}^{\ast}\tau)\cap \{c(i_{Y_1}^{\ast}{\bf H}\otimes 
\pi_{{\bf P}({\bf V}_{canon}^{\circ})}^{\ast}{\bf W}_{canon}^{\circ})\cap 
s(Y_1, X')\}_{dim_{\bf C}X'-e'}$$
$$\cap c_1(i_{Y_1}^{\ast}{\bf
 H})^{{(\underline{C}-\sum e_i)^2-c_1({\bf K}_{M_{n+1}/M_n})\cdot 
(\underline{C}-\sum e_i)\over 2}+p_g+dim_{\bf C}Y(\Gamma)+rank_{\bf C}\tau-r}
\bigr)$$

$$\hskip -1.2in
=\bar{i}_{Y_1\ast}\bigl(
c_{rank_{\bf C}\tau}(i_{Y_1}^{\ast}{\bf H}\otimes \pi_{Y_1}^{\ast}\tau)\cap \{
 c(i_{Y_1}^{\ast}{\bf H}\otimes 
\pi_{{\bf P}({\bf V}_{canon}^{\circ})}^{\ast}{\bf W}_{canon}^{\circ})\cap 
s(Y_1, X')\}_{dim_{\bf C}X'-e'}$$
$$\cap c_1(i_{Y_1}^{\ast}{\bf H})^{
{(\underline{C}-\sum e_i)^2-c_1({\bf K}_{M_{n+1}/M_n})\cdot 
(\underline{C}-\sum e_i)\over 2}+p_g+dim_{\bf C}Y(\Gamma)}\bigr).$$

Thus, it suffices to identify  

$$\bar{i}_{Y_1\ast}\bigl(c_{rank_{\bf C}\tau}(i_{Y_1}^{\ast}{\bf H}\otimes
 \pi_{Y_1}^{\ast}\tau)\cap \{
 c(i_{Y_1}^{\ast}{\bf H}\otimes 
\pi_{{\bf P}({\bf V}_{canon}^{\circ})}^{\ast}{\bf W}_{canon}^{\circ})\cap 
s(Y_1, X\times_{M_n}Y(\Gamma))\}_{dim_{\bf C}X'-e'}\bigr)$$

 with the final expression derived in step IV.

By a direct comparision it suffices to identify  
$c_k(i_{Y_1}^{\ast}{\bf H}\otimes \pi_{Y_1}^{\ast}\tau)$ with 
$\{c(i_{Y_1}^{\ast}{\bf H}\otimes {\bf G})
\cap s(\pi_{Y_1}^{\ast}{\bf N}_{Y(\Gamma)}M_n)\}_{dim_{\bf C}X'-k}$ 
for all $k\in \{0\}\cup {\bf N}$ and argue the vanishing of
 $\{c(i_{Y_1}^{\ast}{\bf H}\otimes 
{\bf G})\cap s(\pi_{Y_1}^{\ast}
{\bf N}_{Y(\Gamma)}M_n)\}_{dim_{\bf C}X'-rank_{\bf C}\tau
 -k}$, $k\in {\bf N}$ by the vanishing of
$c_{rank_{\bf C}\tau+k}({\bf H}\otimes \pi_{Y_1}^{\ast}\tau)=0$, 
$k>0$ for any locally free 
 $\tau$ on $U\times T(M)$.   

\medskip

 Firstly, by using ${\bf N}_{Y(\Gamma)}M_n=\oplus_{i=1}^p h_i^{\ast}{\bf N}_i$,

$$c(i_{Y_1}^{\ast}{\bf H}\otimes 
{\bf G})\cap s(\pi_{Y_1}^{\ast}{\bf N}_{Y(\Gamma)}M_n)
=c(i_{Y_1}^{\ast}{\bf H}\otimes 
{\bf G}-\pi_{Y_1}^{\ast}{\bf N}_{Y(\Gamma)}M_n)$$
$$=c(i_{Y_1}^{\ast}{\bf H}\otimes 
{\bf G}-i_{Y_1}^{\ast}{\bf H}\otimes (\oplus_{1\leq i\leq p}
 \pi_{Y_1}^{\ast}h_i^{\ast}{\bf Q}_i\otimes {\bf E}_C\otimes {\bf N}_i))
\cap_{i=1}^p c(i_{Y_1}^{\ast}{\bf H}\otimes 
\pi_{Y_1}^{\ast}h_i^{\ast}{\bf Q}_i\otimes {\bf E}_C\otimes 
{\bf N}_i-\pi_{Y_1}^{\ast}
h_i^{\ast}{\bf N}_i).$$

 Please consult proposition \ref{prop; residue} in
 step VI for the definitions of ${\bf Q}_i$.

\medskip

By lemma \ref{lemm; normal}, the locus
 $X\times_{M_n}Y(\Gamma_{e_i})\subset X$ is the zero locus of some 
regular global section of the locally free sheaf
 $\pi_{X}^{\ast}
{\cal R}^0\pi_{\ast}\bigl({\cal O}_{\sum_{j_i}E_{j_i}}(E_i)\bigr)$, 
which implies (by lemma \ref{lemm; equal}) the equality 
 $c({\bf H}\otimes {\bf Q}_i\otimes {\bf E}_C\otimes 
\pi_{X'}^{\ast}{\bf N}_i-\pi_{X'}^{\ast}{\bf N}_i)
\equiv 1$ for each $1\leq i\leq p$.

 Thus, we may drop $\cap_{i=1}^p c(i_{Y_1}^{\ast}{\bf H}\otimes 
\pi_{Y_1}^{\ast}h_i^{\ast}{\bf Q}_i\otimes {\bf E}_C\otimes 
{\bf N}_i-\pi_{Y_1}^{\ast}
h_i^{\ast}{\bf N}_i)$ from the above expression and prove that 

$$[{\cal G}]-[\oplus_{1\leq i\leq p}h_i^{\ast}({\cal Q}_i\otimes {\cal E}_C 
\otimes {\cal N}_i)]
\equiv \tau,$$

in $K_0(U\times T(M))$. 
As usual, we have used the calligraphic characters 
${\cal G}$, ${\cal N}_i$,
 ${\cal Q}_i$ and ${\cal E}_C$ to denote the locally free sheaf associated 
 with ${\bf G}$, ${\bf N}_i$, ${\bf Q}_i$ and ${\bf E}_C$.

\medskip

 This equality follows from the definitions of $\tau$ and ${\bf G}$, 
the short exact sequences 
of ${\bf N}_i$ stated in proposition \ref{prop; residue}
 and the following short exact sequences 

$$\hskip -1in 0\mapsto {\cal R}^0\pi_{\ast}\bigl({\cal O}_{\Xi_i\cap 
\cup_{a<i}\Xi_a}(-{\bf M}(E)E)\otimes {\cal E}_C\bigr) 
\mapsto {\cal R}^1\pi_{\ast}\bigl({\cal O}_{\Xi_i}(-{\bf M}(E)E)\otimes 
{\cal E}_{C-\sum_{a=1}^ie_a}\bigr)
\mapsto {\cal R}^1\pi_{\ast}\bigl({\cal O}_{\Xi_i}(-{\bf M}(E)E)\otimes 
{\cal E}_C\bigr)\mapsto 0$$

on $U_i\times T(M)$.
This ends the proof of theorem \ref{theo; question1}. $\Box$

\medskip

\subsection{\bf Some Partial Orderings Among 
$(\Gamma, \sum_{e_i\cdot (C-{\bf M}(E)E)<0}e_i)$}\label{subsection; order}

\bigskip

 At the end, let us address {\bf Question 3} on page \pageref{question3} briefly.
 Theorem \ref{theo; question1} has answered the question of identifying 
 the local contributions of the family invariant under the {\bf
 Simplifying Assumption}. 

 On the other hand, it is clear from proposition \ref{prop; intersect}, the 
 local contributions over the whole $Y(\Gamma)$ (not only from $Y_1$)
 involve more terms yet to
 be identified and are related to the local contributions upon 
$({\bf C}_{\rho}-{\bf O}^{\circ})\cap 
{\bf s}_{canon}^{\circ}$. In principle we may repeat our discussion
 upon $Y(\Gamma)$ to some other admissible $Y(\Gamma'), \Gamma'<\Gamma$. (compare 
with $\sqsupset$ on page \pageref{sqsupset}).

 Our discussion in theorem \ref{theo; question1} makes us believe 
that the mixed family invariants
 similar to the expression in theorem \ref{theo; question1} could have
 appeared while we enumerate these unknown contributions. 

 On the other hand, if we consider the local contributions of
 the family invariants upon the closures of two admissible strata 
$Y(\Gamma_1), Y(\Gamma_2)$, $Y(\Gamma_1)\cap
 Y(\Gamma_2)\not=\emptyset$ simutaneously, it leads to potential 
over-counting
 as the local contributions from $Y(\Gamma_1)\cap Y(\Gamma_2)$ are counted
 twice altegother. 
 The phenomenon becomes more complicated when more than two different $Y(\Gamma)$
 are involved and some combinatorial partial ordering (see page \pageref{partial} 
below)
 upon these admissible graphs 
$\Gamma\in adm(n)$ have to be imposed in order to get a 
 consistent enumeration on ${\cal FASW}$ without over-counting.
 
A few partial orderings among the admissible strata $Y_{\Gamma}$, 
$\Gamma\in adm(n)$
can been introduced as below:

Given a fixed ${\bf M}(E)E=\sum_{1\leq i\leq n}m_i E_i$ encoding the
 singular multiplicities of curve singularities, 
consider all the admissible strata $Y_{\Gamma}, \Gamma\in adm(n)$ 
satisfying the special condition on page \pageref{Scondition}.
 
\medskip  

 One may introduce three partial orderings $>$, $\sqsupset$, $\gg$ 
among all such $(\Gamma, \sum_{e_i\cdot (C-{\bf M}(E)E)<0}e_i)$ by the
 following conditions based on
 the degenerations of type $I$ exceptional classes: \label{partial}

 Let $\Gamma$ and $\Gamma'$ be two admissible graphs.
If the following condition (i) holds, 
 
\noindent (i). $\Gamma>\Gamma'$, i.e. $Y_{\Gamma'}\subset Y(\Gamma)-Y_{\Gamma}$

\medskip

then $(\Gamma, \sum_{e_i\cdot {\bf M}(E)E>0} e_i)$
 is said to be larger than $(\Gamma', \sum_{e_i'\cdot {\bf M}(E)E>0} e_i')$ under
 a partial ordering $>$, \label{>}

 $$(\Gamma, \sum_{e_i\cdot {\bf M}(E)E>0} e_i)
 > (\Gamma', \sum_{e_i'\cdot {\bf M}(E)E>0} e_i').$$

\medskip

 We say that 
$(\Gamma, \sum_{e_i\cdot {\bf M}(E)E>0} e_i)
 \sqsupset (\Gamma', \sum_{e_i'\cdot {\bf M}(E)E>0} e_i'))$ if additional to
 (i)., \label{sqsupset} the following condition (ii). is satisfied,

\noindent (ii). The class $\sum_{e_i\cdot (C-{\bf M}(E)E)<0}e_i-\sum_{
e_j'\cdot (C-{\bf M}(E)E<0}e_j'$ is effective over $Y_{\Gamma'}$. In 
other words for
 all $b\in Y_{\Gamma'}$, the exceptional curve above $b$ dual to each $e_j'$ with
 $e_j'\cdot (C-{\bf M}(E)E)<0$ is an irreducible component of the
 tree of ${\bf P}^1s$ dual to an $e_i$ with $e_i\cdot (C-{\bf M}(E)E<0$.
 
\medskip

We say that $(\Gamma, \sum_{e_i\cdot {\bf M}(E)E>0} e_i)
 \gg (\Gamma', \sum_{e_i'\cdot {\bf M}(E)E>0} e_i')$ under
 $\gg$ if $(i).$ and
 the following condition $(iii).$ hold,

\noindent (iii). For all $e_i$, with $e_i\cdot (C-{\bf M}(E)E)<0$,
 the corresponding $e_i'=e_i$. There exists at least one $1\leq j\leq n$
 such that $e_j^2=-1$ but $e_j'\cdot (C-{\bf M}(E)E)<0$. 
(compare with page 409-410 definition 4.5 of [Liu1])

\medskip

 If the equality may hold, we replace the symbols $>$, $\sqsupset$, $\gg$ by
 $\geq$, $\sqsupseteq$, $\underline{\gg}$.
\medskip

 We end the current paper by the following observation,

\medskip

\begin{prop}\label{prop; middle}
 Let $\Gamma$ satisfies the {\bf special condition} on page \pageref{Scondition}
 let 
 $(\Gamma, \sum_{e_i\cdot (C-{\bf M}(E)E<0}e_i)> (\Gamma', \sum_{e_i'\cdot
 (C-{\bf M}(E)E<0}e_i')$, then there exists an intermidiate
 pair $(\Gamma'', \sum_{e_i''\cdot
 (C-{\bf M}(E)E<0}e_i'')$ such that

$$(\Gamma, \sum_{e_i\cdot (C-{\bf M}(E)E<0}e_i)\sqsupseteq   
(\Gamma'', \sum_{e_i''\cdot (C-{\bf M}(E)E<0}e_i'')\underline{\gg}
(\Gamma', \sum_{e_i'\cdot (C-{\bf M}(E)E<0}e_i').$$
\end{prop}

\medskip

\noindent Proof: Consider the cohomology class $\sum_{e_i\cdot (C-{\bf M}(E)E)<0}
e_i-\sum_{e_j'\cdot (C-{\bf M}(E)E<0} e_j'$. If it has been effective on 
the whole $Y_{\Gamma'}$, then we may take $\Gamma"=\Gamma$ and the statement holds 
trivially. 

 If the class 
$\sum_{e_i\cdot (C-{\bf M}(E)E)<0}e_i-\sum_{e_j'\cdot (C-{\bf M}(E)E)<0} e_j'$ 
is not effective on $Y_{\Gamma'}$, one defines the index sets 
$I=\{ i|1\leq i\leq n,
e_i\cdot (C-{\bf M}(E)E<0\}$ and $J=\{ j|1\leq j\leq n, e_j'\cdot 
(C-{\bf M}(E)E)<0\}$ as subsets of the universal set 
 $\{1, 2, \cdots, n\}$. Then there exists an non-empty subset $J_0\subset J-I$ 
such that

\medskip

\noindent (a). $\sum_{i\in I}e_i-\sum_{j\in J-J_0}e_j'$ 
is effective on $Y_{\Gamma'}$.

\medskip 

\noindent (b). Consider any proper subset $J_1\subset J_0$, then
 $\sum_{i\in I}e_i-\sum_{j\in J-J_1}e_j'$ is non-effective on $Y_{\Gamma'}$.

\medskip

 By a direct calculation
 the intersection pairing $e_i\cdot e_i'\in {\bf Z}$ is always negative.
Thus for all points in $Y_{\Gamma'}$, 
the irreducible ${\bf P}^1$ dual to $e_i'$ is always an irreducible
 component of the tree of ${\bf P}^1$ dual to $e_i$.  Therefore, $e_i-e_i'$, 
 $i\in I$, is effective over $Y_{\Gamma'}$. 
So is their sum $\sum_{i\in I}(e_i-e_i')=
\sum_{i\in I}e_i-\sum_{j\in I}e_j'$.  

As $\sum_{i\in I}e_i-\sum_{j\in I}e_j'$ is
 effective but $\sum_{i\in I}e_i-\sum_{j\in J}e_j'$ is not, there exists
 at least one minimal subset $J_0\subset J-I$ satisfying (a). and (b).

 Then the co-existence of the type $I$ exceptional classes 
$e_j', j\in (J-J_0)\cup J^c$ as 
irreducible rational curves defines an 
admissible stratum encoded by some admissible graph $\Gamma''$ with the
properties
 $e_j''\equiv e_j$ for $j\in (J-J_0)\cup J^c$, 
$e_j''=E_j$, i.e. $(e_j'')^2=-1$, for $j \in J_0$.

 Recall that the closed set 
$Y(\Gamma)$ is the locus in $M_n$ over which the type $I$ 
exceptional classes $e_i, 1\leq i\leq n$
 are effective and it is apparent that all $e_i$, as effective 
combinations of $e_j''=e_j'$, $j \in (J-J_0)\cup J^c$ 
and some $-1$ classes $e_j'', j\in J_0$ (effective over the whole
 $M_n$), are effective over 
 $Y_{\Gamma''}$. Thus $Y(\Gamma)\supset Y_{\Gamma''}$ and $\Gamma>\Gamma''$
 accordingly.
 For a similar reason $\Gamma''>\Gamma'$ as well.

 The non-emptyness of $J_0$ implies 
 $(\Gamma'', \sum_{e_j''\cdot (C-{\bf M}(E)E)<0}e_j'')\gg 
(\Gamma', \sum_{e_j'\cdot (C-{\bf M}(E)E)<0}e_j')$.

 On the other hand, to show that $\sum_{e_i\cdot (C-{\bf M}(E)E)<0}e_i-\sum_{
e_j''\cdot (C-{\bf M}(E)E)<0}e_j''$ is effective over $Y_{\Gamma''}$ (as
 requried by the partial ordering $\sqsupset$), we show
that the difference must be an effective combination of $e_j'', j\in (J-J_0)
\cup J^c$ and
 $e_j'', j\in J_0$.

 Firstly, we know that  
$$\sum_{e_i\cdot (C-{\bf M}(E)E)<0}e_i-\sum_{
e_j''\cdot (C-{\bf M}(E)E)<0}e_j''=\sum_{i\in I}e_i-\sum_{j\in J-J_0}e_j'$$

is effective over $Y_{\Gamma'}$. Thus, it must be an effective combination of
 $e_j', 1\leq j\leq n$, say $\sum_{1\leq j\leq n}c_j e_j$, $c_j\geq 0$.

 We argue that for all $j\in J_0$, $c_j=0$. Otherwise, we may take
 $J_1=\{j|j\in J_0, c_j=0\}\subset J_0$ and make
 $\sum_{i\in I}e_i-\sum_{j\in J-J_1}e_j'$ effective over $Y_{\Gamma'}$,
 violating the minimality condition (b). that $J_0$ satisfies.

 Thus, $\sum_{j\leq n}c_je_j'=\sum_{j\in J- J_0}c_je_j'+\sum_{j\not\in J}c_je_j'$.

 By the defining properties of $\Gamma'$ and $\Gamma''$, 
$e_j'=e_j'', j\in (J-J_0)\cup J^c$. 
Thus all these type $I$ classes $e_j', j\in (J-J_0)\cup J^c$ 
are effective over $Y_{\Gamma''}$ and so is 
 their effective combination 

$$\sum_{j\not\in J_0}c_je_j'=\sum_{e_i\cdot (C-{\bf M}(E)E)<0}e_i-\sum_{
e_j''\cdot (C-{\bf M}(E)E)<0}e_j''.$$  

The proposition is proved. $\Box$

\bigskip

 We point out the geometric origins of these partial orderings.

\medskip

\begin{rem}\label{rem; partial}
 The partial ordering $>$ indicates that the admissible stratum $Y_{\Gamma'}$
 can be degenerated from $Y_{\Gamma}$. It implies that 
 ${\cal M}_{C-{\bf M}(E)E}\times_{M_n}
Y(\Gamma')$ is contained inside ${\cal M}_{C-{\bf M}(E)E}\times_{M_n}Y(\Gamma)$.

The partial ordering $\sqsupset$ indicates that every type $I$ classes $e_j'$
 with $e_j'\cdot (C-{\bf M}(E)E)<0$ are degenerated from some $e_i$ as one of its
 irreducible components. The $\sqsupset$ implies that
family moduli space 
${\cal M}_{C-{\bf M}(E)E-\sum_{e_i\cdot (C-{\bf M}(E)E)<0}e_i}\times_{M_n}
Y(\Gamma')$ is contained inside 
${\cal M}_{C-{\bf M}(E)E-\sum_{e_i'\cdot (C-{\bf M}(E)E)<0}e_i'}\times_{M_n}
Y(\Gamma')$.

The partial ordering $\gg$ indicates that when one degenerates from $Y_{\Gamma}$
 to $Y_{\Gamma'}$, some new type $I$ class $e_j'$
 with $e_j'\cdot (C-{\bf M}(E)E)<0$ appears while the remaining
 $e_i$, $e_i\cdot (C-{\bf M}(E)E)<0$ are preserved. The $\gg$ implies
 that ${\cal M}_{C-{\bf M}(E)E-\sum_{e_i'\cdot (C-{\bf M}(E)E)<0}e_i'}\times_{M_n}
Y(\Gamma')$ is contained in 
${\cal M}_{C-{\bf M}(E)E-\sum_{e_i\cdot (C-{\bf M}(E)E)<0}e_i}
\times_{M_n}Y(\Gamma)$.

Notice that the $\Gamma'$ or $\Gamma''$ in proposition \ref{prop; middle} may
 not always satisfies the {\bf special condition} on page \pageref{Scondition}. 
 This indicates that there is some other admissible strata $Y_{\dot{\Gamma}}$ 
in $M_n$ 
satisfying the {\bf special condition} and $Y_{\Gamma'}$ or $Y_{\Gamma''}$
 is in the intersection of $Y(\Gamma)$ and $Y(dot{\Gamma})$.

\end{rem}

{}


\bigskip

\begin{thebibliography}{}

\bibitem[Be]{Be} K. Behrend,
 {\em Gromov-Witten Invariants in Algebraic Geometry}
 Inventione Math., 
 {\bf 127} pp601-617 (1997)

\bibitem[BGV]{BGV} N. Berline, E Getzler, M. Vergne
 {\em Heat Kernels and Dirac Operators} Grundlehren der mathematischen
 Wissenschaften  {\bf 298} New York: Springer-Verlag (1991)



\bibitem[BPV]{BPV} W. Barth, C. Peters, A. Van De Ven,
 {\em  Compact Complex Surfaces} Ergebnisse der Math. {\bf 4}
 New York: Springer-Verlag (1984)


\bibitem[BT]{BT} R. Bott and L. Tu,
 {\em  Differential Forms in Algebraic Topology} Graduate texts in
 mathematics; {\bf 82}

\bibitem[F]{F} W. Fulton.
{\em Intersection Theory}. A series of Modern Surveys in Mathematics, 
Springer-Verlag, (1984).

\bibitem[F2]{F3} W. Fulton.
{\em Introduction to Toric Varieties}. Annals of Mathematics Studies, 
{\bf 131} Princeton University Press, (1984).



\bibitem[Fr]{Fr} Robert Friedman.
{\em Algebraic Surfaces and Holomorphic Vector Bundles} Universitext, 
Springer-Verlag, (1998).


\bibitem[FM]{FM} R. Friedman and J. Morgan
{\em Algebraic Surfaces and Seiberg-Witten Invariants}.
Journal of Algebraic Geometry, {\bf 6}, no.6, pp445-479, (1997).

\bibitem[FM2]{FM2} R. Friedman and J. Morgan
{\em Obstruction Bundles, Semiregularity, and Seiberg-Witten Invariants}.
Communications of Analysis and Geometry, {\bf 7}, no.3, pp451-495, (1999).

\bibitem[FS]{FS} R. Fintushel, R. Stern
{\em Rational Blowdowns of Smooth $4$-Manifolds}.
Journal of Differential Geometry, {\bf 46}, pp.181-235, (1997).

\bibitem[FS2]{FS2} R. Fintushel, R. Stern
{\em Blowup Formula For Donaldson Invariants}.
Annals of Mathematics (2), {\bf 143}, pp.529-546, (1996).


\bibitem[GH]{GH} P. Griffiths and J. Harris,
 {\em  Principles of Algebraic Geometry}. New York: John Wiley and Sons, (1978).

\bibitem[Got]{Got}  L. G$\ddot o$ttsche.
{\em A Conjectural Generating Function for Numbers of Curves on
Surfaces}. Preprint. {\bf alg-geom 9711012} (1997).



\bibitem[Ha]{Ha} R. Hartshorne,
 {\em   Algebraic Geometry}.
 Graduate texts in Mathematics, {\bf 52}

\bibitem[HM1]{HM1} J. Harvey and G. Moore,
 {\em   Algebras, BPS states, and Strings}.
 hep-th 9510182, {\bf} (1995).

\bibitem[K]{K} M. Karoubi
{\em  K-Theory, an introduction}. Grundlehren der mathematischen Wissenschaften
 {\bf 226}, New York: Springer-Verlag (1970).


\bibitem[KM]{KM} P. Kronheimer. and T. Mrowka. 
{\em The Genus of Imbedded Surfaces
in the Projective Spaces}. Math. Research Letters.
 {\em 1}, pp. 797-808 (1994).

\bibitem[LiT1]{LiT1} J. Li,  G. Tian
 {\em Virtual Moduli Spaces and Gromov-Witten Invariants of 
Algebraic Varieties}. Journal of Amer. Math. Soc. {\bf 11}, 
 pp119-174 (1998)

\bibitem[LiT2]{LiT2} J. Li,  G. Tian
 {\em Virtual Moduli Spaces and Gromov-Witten Invariants of 
 General Symplectic Manifolds}. Preprint alg-geom/9608032


\bibitem[Liu]{Liu} A. K. Liu,
 {\em  Ph.D thesis} Harvard University, (1996). 


\bibitem[Liu1]{Liu1} A. K. Liu,
 {\em  Family Blowup Formula, Admissible Graphs and the 
Enumeration of Singular Curves (I)}. 
  Journal of Differential Geometry {\bf 56} pp381-579 (2001)

\bibitem[Liu2]{Liu2} A. K. Liu,
 {\em  Cosmic String and Family Seiberg-Witten Theory} 
 preprint, (2003)

\bibitem[Liu3]{Liu3} A. K. Liu,
 {\em  The Family Blowup Formula of the Family Seiberg-Witten Invariants}
 preprint, DG0305294, (2003)

\bibitem[Liu4]{Liu4} A. K. Liu,
 {\em  A Note on Curve Counting Scheme in an Algebraic Family and the
 Admissible Decomposition Classes}
 preprint DG/0308196, (2003)

\bibitem[Liu5]{Liu5} A. K. Liu,
 {\em  The Algebraic Proof of the Universality Theorem}
 in preparation, (2003)



\bibitem[LL1]{LL1} T. J. Li and A. K. Liu,
 {\em  Family Seiberg-Witten Invariant and Wall Crossing Formula}. 
 to appear in Communications of Analysis and Geometry.

\bibitem[LL2]{LL2} T. J. Li and A. K. Liu. 
 {\em  General Wall Crossing Formula}. Mathematical Research Letters.
 {\bf 2} pp. 97-118, (1995)

\bibitem[LL3]{LL3} T. J. Li and A. K. Liu. 
 {\em The Symplectic Structures of Rational and Ruled Surfaces and 
the generalized Adjunction inequality}. 
 {\bf 2} pp. 453-471, (1995)




\bibitem[Law]{Law} B.H. Lawson Jr. and M-L. Michelson. 
 {\em Spin Geometry} Princeton Mathematical Series 
 {\bf 38} Princeton Univ. Press, Princeton (1989)



\bibitem[Mc]{Mc} D. Mcduff.
 {\em Remark on the Uniqueness of Symplectic Blowing Ups}, 
London Math. Soc. Lecture Note Ser.{\bf 192} Symplectic Geometry, pp157-167, (1993)

\bibitem[Mor]{Mor} J. W. Morgan,
 {\em The Seiberg-Witten Equations and Applications to the Topology of 
Smooth Four-Manifolds}, Mathematical Notes (Princeton University Press)
{\bf 44}  (1996)



\bibitem[R]{R} Y. Ruan.
 {\em Virtual Neighborhood and Pseudo-Holomorphic Curves},
{\bf alg-geom/9611021} preprint, (1996)





\bibitem[RT1]{RT1} Y. Ruan and G. Tian. 
 {\em The Mathematical Theory of Quantum Cohomology}, Journal of
Differential Geometry. {\bf 42} no. 2. Sep. pp 259-367, (1995).

\bibitem[RT2]{RT2} Y. Ruan and G. Tian. 
 {\em Higher Genus Symplectic Invariants and Sigma Models coupled with
 Gravity}, Inventione Mathematicae. {\bf 130} no. 3. pp. 455-516, (1997).


\bibitem[S]{S} B. Siebert 
 {\em Gromov-Witten Invariants for General Symplectic Manifolds}. 
{\bf dg-ga/9608005}, preprint, (1996)





\bibitem[T1]{T1} C.H. Taubes.
 {\em SW $\rightarrow$ Gr, from the Seiberg-Witten Equations to 
 Pseudo-holomorphic Curves}. Journal of American Mathematical Society. {\bf 9} 
no. 3. (1996).


\bibitem[T2]{T2} C.H. Taubes,
 {\em  Gr $\rightarrow$ SW, from Pseudo-holomorphic Curves to
 Seiberg-Witten Solutions}, Journal of Differential Geometry, {\bf 51}, 
pp203-334. (1999).

\bibitem[T3]{T3} C.H. Taubes,
{\em  The Seiberg Witten Invariants and Gromov Invariant},
Mathematical Research Letters. {\bf 2} pp221-238, (1995).

\bibitem[V]{V} Israel Vainsencher.
 {\em Enumeration of n-fold tangent hypersurfaces to a surface},
 Journal of ALgebraic Geometry. {\bf 4} pp503-526. (1995).


\bibitem[W]{W} Edward Witten.
 {\em Monopoles and Four Manifolds}.
  Math. Research Letters. {\bf 1}, pp. 769-796. (1994).


\bibitem[YZ]{YZ} S. T. Yau and E. Zaslow.
 {\em  BPS states, string duality, and Nodal Curves on $K3$}. Nuclear Physics B.
  {\bf 471} no. 3. pp503-512. (1996).



\end{thebibliography}
\end{document}